\documentclass[a4paper,11pt,reqno]{amsart}
\usepackage[english]{babel}
\usepackage{amsmath}
\usepackage{amscd}
\usepackage{amsgen}
\usepackage{amsfonts}
\usepackage{amssymb}
\usepackage{amsthm}
\usepackage{comment}
\usepackage{stmaryrd}
\usepackage{amsmath}
\usepackage[all]{xy}
\usepackage{latexsym}
\usepackage{graphics}
\usepackage{verbatim}
\usepackage{url}
\usepackage{slashed}
\usepackage[mathscr]{euscript}

\newtheorem{theorem}{Theorem}
\newtheorem{prop}[theorem]{Proposition}
\newtheorem{prop*}{Proposition}
\newtheorem{lem}[theorem]{Lemma}
\newtheorem{cor}[theorem]{Corollary}
\newtheorem{cor*}{Corollary}
\newtheorem{conj}[theorem]{Conjecture}
\newtheorem{defn}[theorem]{Definition}
\newtheorem{rem}[theorem]{Remark}

\newtheorem{thm}{Theorem}

\newcommand{\lo}{\mathring{L}}                      
\newcommand{\JF}{\mathcal{F}}                      
\newcommand{\F}{\mathring{\mathcal{F}}}      
\newcommand{\Hess}{\operatorname{Hess}}   
\newcommand{\tr}{\operatorname{tr}}              
\newcommand{\id}{\operatorname{Id}}            
\newcommand{\Div}{\operatorname{div}}         
\newcommand{\W}{\overline{\mathcal{W}}}     
\newcommand{\CI}{{\mathcal C}}                      
\newcommand{\Ric}{\operatorname{Ric}}          
\newcommand{\scal}{\operatorname{scal}}        
\newcommand{\grad}{\operatorname{grad}}      
\newcommand{\SC}{\mathcal{S}}                      
\newcommand{\Sc}{\mathscr{S}}                       
\newcommand{\res}{\operatorname{Res}}        
\newcommand{\G}{{\mathcal G}}                      

\newcommand{\N}{{\mathbb N}}

\newcommand{\Pf}{\operatorname{Pf}}             
\newcommand{\End}{\operatorname{End}}       
\newcommand{\LOP}{\mathcal {D}}                  

\newcommand{\Iv}{\mathcal{I}}                        
\newcommand{\Jv}{\mathcal{J}}                       

\def\I{{\sf I}}                                                   
\def\H{{\mathcal H}}
\def\T{{\mathcal T}}                                        
\def\R{{\mathbb R}}                                        
\def\Rho{{\sf P}}                                             
\def\J{{\sf J}}                                                 
\def\B{{\mathcal B}}                                        
\def\st{\stackrel{\text{def}}{=}}

\numberwithin{theorem}{section} \numberwithin{equation}{section}

\title{Extrinsic Paneitz operators and $Q$-curvatures for hypersurfaces}

\author[Andreas Juhl]{Andreas Juhl}

\address{Humboldt-Universit\"at, Institut f\"ur Mathematik, Unter den Linden 6, 10099 Berlin, Germany}

\email{juhl.andreas@googlemail.com}


\begin{document}

\begin{abstract}
For any hypersurface $M$ of a Riemannian manifold $X$, recent works introduced the
notions of extrinsic conformal Laplacians and extrinsic $Q$-curvatures. Here we
derive explicit formulas for the extrinsic version ${\bf P}_4$ of the Paneitz operator and 
the corresponding extrinsic fourth-order $Q$-curvature ${\bf Q}_4$ in general
dimensions. In the critical dimension $n=4$, this result yields a closed formula for the
global conformal invariant $\int_M {\bf Q}_4 dvol$ (for closed $M$) and various decompositions 
of ${\bf Q}_4$, which are analogs of the Alexakis/Deser-Schwimmer type decompositions of global
conformal invariants. These results involve a series of obvious local conformal
invariants of the embedding $M^4 \hookrightarrow X^5$ (defined in terms of the Weyl
tensor and the trace-free second fundamental form) and a non-trivial local conformal
invariant $\mathcal{C}$. In turn, we identify $\mathcal{C}$ as a linear combination of two
local conformal invariants $\Jv_1$ and $\Jv_2$. We also observe that these are special
cases of local conformal invariants for hypersurfaces in backgrounds of general
dimension. Moreover, in the critical dimension $n=4$, a linear combination of $\Jv_1$ and
$\Jv_2$ can be expressed in terms of obvious local conformal invariants of the
embedding $M \hookrightarrow X$. This finally reduces the non-trivial part of the
structure of ${\bf Q}_4$ to the non-trivial invariant $\Jv_1$. For totally umbilic $M$, 
the invariants $\Jv_i$ vanish, and the formula for ${\bf P}_4$ substantially simplifies. 
For closed $M^4 \hookrightarrow \R^5$, we relate the integrals of $\Jv_i$ to functionals 
of Guven and Graham-Reichert. Moreover, we establish a Deser-Schwimmer type decomposition 
of the Graham-Reichert functional of a hypersurface $M^4 \hookrightarrow X^5$ in general 
backgrounds. In this context, we find one further local conformal invariant $\Jv_3$. Finally, we 
derive an explicit formula for the singular Yamabe energy of a closed $M$. The resulting 
explicit formulas show that it is proportional to the total extrinsic fourth-order $Q$-curvature. 
This observation confirms a special case of a general fact and serves as an additional cross-check
of our main result. Throughout, we carefully discuss the relation of our formulas to the recent literature. 
\end{abstract}

\maketitle

\begin{center} \today \end{center}

\tableofcontents

\section{Introduction and formulation of the main results}\label{results}

The significance of the Yamabe operator
$$
    P_2 = \Delta - \left(\frac{n}{2}-1\right) \J
$$
and the Paneitz operator
$$
    P_4 = \Delta^2 - \delta ((n-2) \J h - 4\Rho) d + \left(\frac{n}{2}-2\right)
    \left(\frac{n}{2} \J^2 - 2 |\Rho|^2 - \Delta (\J) \right)
$$
in geometric analysis is well-known (\cite{sharp, CY, CGY, CBull, IMC, DGH, J1, BJ} 
and references therein).
These differential operators are defined on any Riemannian manifold $(M,g)$ of
respective dimension $n \ge 2$ and $n \ge 3$.  Paneitz \cite{pan} discovered the
operator $P_4$ in general dimensions in $1983$. Around the same time, it appeared in
dimension $4$ in several other contexts in \cite{FT, Rieg, ES}.
Here we use the following conventions. The dimension of $M$ is denoted by $n$,
$\delta$ is the divergence operator on $1$-forms, $\Delta = \delta d$ the
non-positive Laplacian, $2(n-1) \J = \scal$ and $(n-2)\Rho = \Ric - \J h$. $\Rho$ is
the Schouten tensor of $h$. It naturally acts on $1$-forms. The operators $P_2$ and
$P_4$ are the first two elements in the sequence of so-called GJMS-operators
$P_{2N}$ (with $2N \le n$ for even $n$ and $N \ge 1$ for odd $n$) \cite{GJMS}. These
self-adjoint geometric differential operators have as leading term a power of the
Laplacian $\Delta$ and are covariant
\begin{equation}\label{CTL-GJMS}
    e^{(\frac{n}{2}+N)\varphi} P_{2N} (\hat{h})(f) = P_{2N}(h)(e^{(\frac{n}{2}-N) \varphi} f)
\end{equation}
under conformal changes $h \mapsto \hat{h} = e^{2\varphi} h$, $\varphi \in
C^\infty(M)$, of the metric. The original definition of the GJMS-operators rests on
the ambient metric of Fefferman and Graham \cite{FG-final}. 

The quantities $Q_2 = \J$ and $Q_4  = \frac{n}{2} \J^2 - |\Rho|^2 - \Delta (\J)$ are known as 
Branson's $Q$-curvature (of respective order $2$ and $4$). The general Branson's $Q$-curvatures 
are defined by 
$$
   P_{2N}(1) = \left(\frac{n}{2}-N\right) (-1)^N Q_{2N}
$$
if $2N < n$, and by a continuation in dimension argument for $2N=n$. For even $n$, the 
case $2N =n$ will be referred to as the critical case. A remarkable property of the critical 
$Q$-curvature $Q_n$ is its transformation law
\begin{equation}\label{CTL-Q}
    e^{n\varphi} Q_n(\hat{h}) = Q_n(h) + (-1)^{\frac{n}{2}} P_n(h) (\varphi), \; \varphi \in C^\infty(M).
\end{equation}
This property may be derived from the conformal covariance \eqref{CTL-GJMS} in the
non-critical case using a continuation in dimension argument \cite{sharp}. Combining
\eqref{CTL-Q} with the self-adjointness of $P_n$ and $P_n(h)(1)=0$ it follows that
for a closed manifold $M$ the integral
$$
   \int_M Q_n(h) dvol_h
$$
is a (global) conformal invariant. This invariant for the conformal class $[h]$ is
related to the conformal anomaly of the renormalized volume of a Poincar\'e-Einstein
metric associated to $h$ \cite{GZ}. In dimension $n=4$, the integrand of this
anomaly is just a constant multiple of $\J^2 - |\Rho|^2$. The complexity of formulas
for $P_{2N}$ and $Q_{2N}$ dramatically increases with $N$. However, a recursive
structure enables one to derive explicit formulas at least for $P_6$ and $P_8$. For more
details, we refer to \cite{J-ex, J-heat, FG-J}.

In the recent works \cite{GW-LNY, GW-RV, JO1}, generalizations ${\bf P}_{N}$ of the GJMS-operators 
$P_N$ and ${\bf Q}_N$ of Branson's $Q$-curvatures $Q_N$ were introduced in the context of the 
singular Yamabe problem for hypersurfaces. These are called extrinsic conformal Laplacians and extrinsic
$Q$-curvatures. In \cite{GW-LNY}, the operators ${\bf P}_N$ were defined in terms of
conformal tractor calculus. The alternative approach in \cite{JO1} rests on an
extension of the notion of residue families as developed on \cite{J1}. The latter
method proves the self-adjointness of all extrinsic conformal Laplacians
by connecting them to scattering theory (see Theorem \ref{ECL-S}).


We assume that the closed manifold $M$ is the boundary of a compact manifold $(X,g)$
with $h$ being induced by $g$. Let $\iota: M \hookrightarrow X$ denote the
embedding. The operators ${\bf P}_N(g)$ still act on $C^\infty(M)$ and, for even
$N$, have leading term a power of the Laplacian of $M$.\footnote{In \cite{JO1}, the leading 
term of ${\bf P}_{2N}$ is a constant multiple of $\Delta^N$.} But their lower-order 
terms depend on the embedding $\iota$. For odd $N$, their leading terms depend on the
trace-free part $\lo = L - H h$ of the second fundamental form $L$. Here $H$ is the
mean curvature, i.e., $\tr(L) = nH$. The embedding $M \hookrightarrow X$ is called
totally umbilic if $\lo=0$. Again, the operators ${\bf P}_N$ are self-adjoint, and
they are covariant
\begin{equation}\label{CTL-CL}
    e^{\frac{n+N}{2} \iota^*(\varphi)} {\bf P}_{N} (\hat{g})(f) 
   = {\bf P}_{N}(g)(e^{\frac{n-N}{2} \iota^*(\varphi)} f)
\end{equation}
under conformal changes $g \mapsto \hat{g} = e^{2\varphi} g$, $\varphi \in
C^\infty(X)$. We recall that, in contrast to \eqref{CTL-GJMS}, the property \eqref{CTL-CL}
concerns conformal changes of the metric $g$ on the ambient space $X$.

In the following, extrinsic conformal Laplacians and extrinsic $Q$-curvatures will
be denoted by boldface letters. For simplicity, we often omit their dependence on
the metric. The operator ${\bf P}_1$ vanishes, and the first two non-trivial
extrinsic conformal Laplacians are given by (see \cite[Proposition 8.5]{GW-LNY},
\cite[Sections 13.10-13.11]{JO1})

\begin{prop*}\label{P23} It holds
\begin{equation*}\label{ECL-2}
    {\bf P}_2(g) = P_2(h) + \frac{n-2}{4(n-1)} |\lo|^2, \quad n \ge 2
\end{equation*}
and
\begin{equation*}\label{ECL-3}
    {\bf P}_3(g) = 8 \delta (\lo d) + \frac{n-3}{2} \frac{4}{n-2} (\delta \delta (\lo) - (n-3) (\lo,\Rho) + (n-1) (\lo,\JF)),
   \quad n \ge 3.
\end{equation*}
Here $\JF$ is the conformally invariant Fialkow tensor (see \eqref{FT}). The
corresponding extrinsic $Q$-curvatures are
\begin{equation*}\label{EQ2}
   {\bf Q}_2(g) = Q_2(h) - \frac{1}{2(n-1)}|\lo|^2 \quad \mbox{with} \quad Q_2(h) = \J^h
\end{equation*}
and
\begin{equation*}\label{EQ3}
   {\bf Q}_3(g) = \frac{4}{n-2} (\delta \delta (\lo) - (n-3) (\lo,\Rho) + (n-1) (\lo,\JF)).
\end{equation*}
These satisfy the fundamental transformation laws
$$
    e^{2 \iota^*(\varphi)} {\bf Q}_2(\hat{g}) = {\bf Q}_2(g) - {\bf P}_2(g)(\varphi)
$$
if $n=2$ and
$$
   e^{3 \iota^*(\varphi)} {\bf Q}_3(\hat{g}) = {\bf Q}_3(g) + {\bf P}_3(g)(\varphi)
$$
if $n=3$.
\end{prop*}

Here we used the general convention $(-1)^{N} {\bf P}_{2N}(1) = (\frac{n}{2}-N) {\bf Q}_{2N}$
and ${\bf P}_N (1) = \frac{n-N}{2} {\bf Q}_N$ for odd $N$\footnote{In \cite{JO1}, we used 
the convention ${\bf P}_N (1) = \frac{n-N}{2} {\bf Q}_N$ for all $N$.} (see also the comments after 
Theorem \ref{ECL-S}). 

In the critical dimension $n=2$, the extrinsic $Q$-curvature $\bf{Q}_2$ is a linear combination of 
$\J$ and $|\lo|^2$. The integrals of both terms are conformally invariant. This implies the conformal 
invariance of $\int_M {\bf Q}_2 dvol_h$ for a closed surface.

In the critical dimension $n=3$, the extrinsic $Q$-curvature ${\bf Q}_3$ is proportional to $(\lo,\JF) $, 
up to a total divergence. This implies the conformal invariance of $\int_M {\bf Q}_3 dvol_h$ for a closed 
$M$. In connection with the study of conformal anomalies for CFT on manifolds with boundaries, it has been
argued in \cite{Fu, Sol} that the boundary terms of anomalies are linear
combinations of $(\lo,\W)$ and $\tr(\lo^3)$. The Fialkow equation \eqref{Fial} shows
that $(\lo,\JF)$ can be written as such a linear combination. Note also the formula
or ${\bf Q}_3$ in general dimensions has a simple pole at $n=2$ and its residue at
$n=2$ is a multiple of \footnote{For a justification of the limit, we refer to
\cite[Section 13.10]{JO1}.}
$$
   \delta \delta (\lo) + (\lo,\bar{\Rho}) + H |\lo|^2.
$$
The fact that this quantity is a constant multiple of the singular Yamabe
obstruction $\B_2$ is a special case of \cite[Theorem 11.6]{JO1}. We shall see
another instance of it in connection with ${\bf Q}_4$.

Now we state the main result of this paper. The following theorem displays a formula
for the extrinsic Paneitz operator ${\bf P}_4$ for general background metrics in
general dimensions. The formulation requires some more notation. In the following,
we shall use a bar to distinguish curvature quantities of the background metric $g$
on $X$ from curvature quantities of the induced metric $h$ on $M$. Accordingly, it
will be convenient to denote the background metric $g$ also by $\bar{g}$. Let
$\overline{W}$ be the Weyl tensor of the background metric $\bar{g}$ and let
$\W_{ij} = \overline{W}_{0ij0}$ be defined by inserting a unit normal vector $N =
\partial_0$ of $M$ into the first and the last slot of $\overline{W}$. Then $\W$ is
a trace-free conformally invariant symmetric bilinear form on $M$. It naturally acts
on $1$-forms on $M$. The component $\bar{\Rho}_{00}$ is defined by inserting two
unit normal vectors into the Schouten tensor $\bar{\Rho}$ of $\bar{g}$. Let
$\bar{\nabla}$ be the Levi-Civita derivative of the background metric $\bar{g}$. The
symbol $\delta$ also will be used for the divergence operator on symmetric bilinear
forms on $M$.

\begin{thm}\label{main} Assume that $n \ge 4$. Then
\begin{align}\label{P4-gen-dim}
   {\bf P}_4 & = \Delta^2 - \delta((n-2)\J h - 4 \Rho) d \notag \\
   & + \delta \left( 4 \tfrac{3n-5}{n-2} \lo^2 + \tfrac{n^2-12n+16}{2(n-1)(n-2)} |\lo|^2 h + \tfrac{4(n-1)}{n-2} \W \right) d
   + \left(\tfrac{n}{2}-2\right)  {\bf Q}_4.
\end{align}
Here ${\bf Q}_4$ is the sum of
\begin{itemize}
\item the intrinsic $Q_4$ of $(M,h)$,
\item the four divergence terms
\begin{align}\label{div-terms}
  \tfrac{2(n-1)}{(n-3)(n-2)} \delta \delta (\W) + \tfrac{2(n-1)}{(n-3)(n-2)} \delta \delta (\lo^2)
   + \tfrac{4}{n-3} \delta (\lo \delta(\lo)) + \tfrac{3n-4}{2(n-1)(n-2)} \Delta (|\lo|^2),
\end{align}
\item the derivative terms
\begin{equation}\label{d-terms}
    2 \lo^{ij} \bar{\nabla}_0(\bar{\Rho})_{ij}
   - \tfrac{4}{n-3} \lo^{ij} \bar{\nabla}_0 (\overline{W})_{0ij0}+ 2 (\lo,\Hess(H)),
\end{equation}
\item the four $\W$-terms
\begin{align}\label{W-terms}
    \tfrac{2 (n-1)^2}{(n-3)(n-2)^2} |\W|^2 - \tfrac{2(n-4)(n-1)}{(n-3)(n-2)} (\Rho,\W)
   + \tfrac{4 (3n-5)(n-1)}{(n-3)(n-2)^2} (\lo^2,\W) - \tfrac{2(n-1)^2}{(n-3)(n-2)} H (\lo,\W),
\end{align}
\item the four Schouten tensor terms
\begin{align}\label{P-terms}
    - \tfrac{2(n^2-9n+12)}{(n-3)(n-2)} (\lo^2,\Rho)
    - \tfrac{n^3-5n^2+18n-20}{2(n-3)(n-2)(n-1)} \J |\lo|^2  + 2 H(\lo,\Rho) - 2 |\lo|^2 \bar{\Rho}_{00}
\end{align}
\item and the four quartic $L$-terms
\begin{align}\label{L-terms}
    & - 3 H^2 |\lo|^2 - \tfrac{2(n-3)}{n-2} H \tr(\lo^3) + \tfrac{2(5n^2-14n+9)}{(n-3)(n-2)^2} \tr(\lo^4)
    - \tfrac{15n^4-49n^3+36n^2+24n-32}{8(n-3)(n-2)^2(n-1)^2} |\lo|^4.
\end{align}
\end{itemize}
Here all scalar products, norms, traces, and divergences are defined by the metric $h$.
\end{thm}

The following comments are in order.

As noted above, in \cite{JO1} a different normalization of ${\bf P}_4$ has been used. In fact, in \cite{JO1} 
we set ${\bf P}_4(g) = 9 P_4(h)$ if $g$ is the conformal compactification of a Poincar\'e-Einstein metric with
conformal infinity $[h]$. Here we adopt the convention that
${\bf P}_4(g) = P_4(h)$ in the Poincar\'e-Einstein case. This implies that
$$
    \int_{M^4} {\bf Q}_4 dvol_h = 3! 4! /9 \int_{M^4} V_4 dvol_h = 16 \int_{M^4} V_4 dvol_h,
$$
where $V_4$ is the fourth singular Yamabe renormalized volume coefficient. For a discussion of the Yamabe 
energy $\int_M V_4 dvol_h$, we refer to Section \ref{SYE}.

Formula \eqref{P4-gen-dim} makes the self-adjointness of ${\bf P}_4$ obvious.

Of course, the four total divergence terms in \eqref{div-terms} vanish by integration on a closed $M$. Since 
$(\lo^2,\W)$, $|\W|^2$ and the quartic terms $\tr (\lo^4)$, $|\lo|^4$ in \eqref{L-terms} are conformally invariant, 
one obtains another conformally covariant operator by removing these terms.

The sum in Theorem \ref{main} describing ${\bf Q}_4$ has a simple pole in $n=3$. The
residue of ${\bf Q}_4$ at $n=3$ equals the sum of the total divergence terms
$$
    4 \delta \delta (\W) + 4 \delta \delta(\lo^2) + 4 \delta (\lo \delta(\lo))
$$
the normal derivative term
$
   - 4 \lo^{ij} \bar{\nabla}_0(\overline{W})_{0ij0},
$
the Weyl tensor terms
$$
    8 |\W|^2 + 4  (\Rho,\W) + 32 (\lo^2,\W) - 8 H (\lo,\W),
$$
the Schouten tensor terms
$$
     12 (\lo^2,\Rho) - 4 \J |\lo|^2
$$
and
$$
    24 \tr(\lo^4) - 8  |\lo|^4.
$$
The sum of these terms actually coincides with $24\B_3$, where $\B_3$ is the singular Yamabe 
obstruction of the embedding $M^4 \hookrightarrow X^5$ \cite[Proposition 1.1]{GGHW}, \cite[Theorem 1]{JO2}. 
This relation is another special case of \cite[Theorem 11.6]{JO1}, and it gives another proof of the conformal
invariance of $\B_3$. In \eqref{d-terms}, the normal derivative of the Weyl tensor has a coefficient that is 
singular at $n=3$. The connection between ${\bf Q}_4$ and $\B_3$ actually explains the appearance of this term 
in ${\bf Q}_4$ by its appearance in $\B_3$.

The conformal covariance of the operator displayed in Theorem \ref{main} can be
confirmed by direct calculations - for an outline of the arguments, we refer to
Section \ref{covariance-direct}.

The formula in Theorem \ref{main} is written in terms of the tensors $\W$ and
$\bar{\Rho}$ of the background metric, their first-order normal derivatives as well
as the intrinsic Schouten tensor $\Rho$ (and its trace $\J$) and the second
fundamental form $L$ (and its trace $nH$) of the embedding $M \hookrightarrow X$.
Alternative formulas can be obtained using the trace-free part of the Fialkow tensor
$\JF$.

In view of the particular significance of the result in the critical dimension
$n=4$, we separately formulate this special case. For $n=4$, Theorem \ref{main}
reduces to

\begin{cor*}\label{P4-crit} In the critical dimension $n=4$, the extrinsic Paneitz operator ${\bf P}_4$ is given by
\begin{equation}\label{EP4-crit}
    {\bf P}_4  = \Delta^2 - \delta (2 \J h - 4 \Rho) d + \delta \left(14 \lo^2 - \frac{4}{3} |\lo|^2 h + 6 \W \right) d.
\end{equation}
\end{cor*}

Note that the right-hand side of \eqref{EP4-crit} is a sum of $P_4$ and three
individually conformally covariant operators.

Next, we consider the extrinsic $Q$-curvature ${\bf Q}_4$ in the critical dimension
$n=4$ more closely. Let $M^4$ be closed. Like the total integral of $Q_4$, the total
integral of ${\bf Q}_4$ over $M$ is a global conformal invariant. In fact, combining
the conformal transformation property \cite[Section 10]{JO1}
\begin{equation}\label{ECTL}
    e^{4\iota^*(\varphi)} {\bf Q}_4 (\hat{g}) = {\bf Q}_4(g) + {\bf P}_4(g)(\varphi)
\end{equation}
with ${\bf P}_4(g)(1)=0$ and the self-adjointness of ${\bf P}_4$, shows that the
total integral of ${\bf Q}_4$ is an invariant of the conformal class $[g]$. The
following result describes this global conformal invariant.

\begin{cor*}\label{Q4-g-int} For a closed $4$-manifold $M$, it holds
\begin{align}\label{Q4-g-invariant}
     & \int_M {\bf Q}_4 dvol_h =  \int_M \left( 2 \J^2 - 2|\Rho|^2 + \frac{9}{2} |\W|^2\right) dvol_h \notag \\
     & + \int_M \left(2 (\lo,\bar{\nabla}_0 (\bar{\Rho}))
     - 4 \lo^{ij} \bar{\nabla}_0 (\overline{W})_{0ij0} + 2 (\lo,\Hess(H)) + 2 H (\lo,\Rho) - 9 H (\lo,\W) \right) dvol_h \notag \\
     & + \int_M \left( 8 (\lo^2,\Rho) -2 \bar{\Rho}_{00} |\lo|^2 - 3 \J |\lo|^2  - 3 H^2 |\lo|^2
    + 21 (\lo^2,\W) - H \tr(\lo^3) \right) dvol_h \notag \\
     & + \int_M \left( \frac{33}{2} \tr (\lo^4) - \frac{14}{3} |\lo|^4 \right) dvol_h.
\end{align}
\end{cor*}

Several comments are in order.

The first integral on the right-hand side of \eqref{Q4-g-invariant} does not depend
on $L$, the second integral is linear in $\lo$, all terms except the last one in the
third integral are quadratic in $\lo$, and the terms in the last integral are
quartic in $\lo$.

The right-hand side of \eqref{Q4-g-invariant} can be decomposed as a sum of a series
of global conformal invariants. More precisely, the first integral is the sum of the
global conformal invariant
$$
    \int_M Q_4(h) dvol_h = \int_M (2\J^2 - 2 |\Rho|^2) dvol_h
$$
and the integral of a constant multiple of the local conformal invariant $|\W|^2$.
The quartic terms $|\lo|^4$, $\tr(\lo^4)$ and $(\lo^2,\W)$ are local conformal
invariants. The conformal invariance of $\int_M {\bf Q}_4 dvol_h$ implies that the
remaining terms define a global conformal invariant. In fact, this sum is the
integral of another local conformal invariant. We set
\begin{align}\label{new invariant}
      \CI & \st 2 (\lo,\bar{\nabla}_0 (\bar{\Rho}))
     - 4 \lo^{ij} \bar{\nabla}_0 (\overline{W})_{0ij0} + 2 (\lo,\Hess(H)) + 2 H (\lo,\Rho) - 9 H (\lo,\W) \notag \\
     & + 8 (\lo^2,\Rho) -2 \bar{\Rho}_{00} |\lo|^2 - 3 \J |\lo|^2  - 3 H^2 |\lo|^2
    - H \tr(\lo^3) + 2 \delta \delta (\lo^2) + \frac{1}{2} \Delta(|\lo|^2).
\end{align}
All terms in the latter sum except the last two divergence terms were taken from
\eqref{Q4-g-invariant}. Note that $\mathcal{C} = 0$ if $\lo=0$. The advantage of
adding these two divergence terms becomes clear in the following result.

\begin{thm}\label{LCI} Let $n=4$. Then $e^{4 \iota^*(\varphi)} \hat{\mathcal{C}} = \mathcal{C}$ for all
$\varphi \in C^\infty(X)$, i.e., $\mathcal{C}$ is a local invariant of the conformal class $[g]$.
\end{thm}

The local conformal invariant $\mathcal{C}$ contains two terms with normal
derivatives of the curvature tensor of the background metric: $(\lo,\bar{\nabla}_0
(\bar{\Rho}))$ and $\lo^{ij} \bar{\nabla}_0 (\overline{W})_{0ij0}$. It turns out
that $\mathcal{C}$ is a linear combination
\begin{equation}\label{C-deco}
   \mathcal{C} = - 4 \Jv_1+2 \Jv_2
\end{equation}
of {\em two} local conformal invariants $\Jv_i$ containing these two normal
derivative terms, respectively. These local invariants are given by\footnote{The integrated invariant $\Jv_1$ also 
appears in \cite{AS} in the context of anomalies of CFT's on manifolds with boundary. Generalizations 
for embeddings $M^4 \hookrightarrow X^n$ with $n \ge 5$ have been found in \cite{CHBRS}.}
\begin{equation}\label{CI-1}
   \Jv_1 \st \lo^{ij} \bar{\nabla}_0(\overline{W})_{0ij0} + 2 H (\lo,\W) + \frac{2}{9} |\delta(\lo)|^2
   - 2 (\lo^2,\Rho) + \J |\lo|^2 - \delta \delta (\lo^2)
\end{equation}
and
\begin{align}\label{CI-2}
   \Jv_2 \st & (\lo, \bar{\nabla}_0(\bar{\Rho})) + H (\lo,\Rho) - \frac{1}{2} H (\lo,\W) + (\lo,\Hess(H)) \notag \\
   & + \frac{4}{9} |\delta(\lo)|^2
   - \bar{\Rho}_{00} |\lo|^2 + \frac{1}{2} \J |\lo|^2 - \frac{3}{2} H^2 |\lo|^2 - \frac{1}{2} H \tr(\lo^3)
   - \delta \delta (\lo^2) + \frac{1}{4} \Delta (|\lo|^2).
\end{align}
For a closed hypersurface $M^4 \hookrightarrow X^5$, the integrals of these
invariants define generalizations of the conformally invariant Willmore functional
of hypersurfaces $M^2 \hookrightarrow X^3$ in the sense that the leading terms of
their Euler-Lagrange equations (for variations of the embedding) are constant
multiples of $\Delta^2(H)$.


It seems to be of independent interest that both invariants $\Jv_1$ and $\Jv_2$
allow generalizations for general dimensions. In fact, the quantities
\begin{align*}
    \Jv_1 & \st  \lo^{ij} \bar{\nabla}_0(\overline{W})_{0ij0} + 2 H (\lo,\W) + \tfrac{n-2}{(n-1)^2} |\delta(\lo)|^2 \notag \\
    & - \tfrac{n-2}{n-3} (\lo^2,\Rho) - \tfrac{n-2}{(n-3)(n-6)} \J |\lo|^2 + \tfrac{n-4}{(n-3)(n-6)} \Delta (|\lo|^2) -
    \tfrac{1}{n-3} \delta \delta (\lo^2)
\end{align*}
and
\begin{align*}
    \Jv_2 & \st (\lo,\bar{\nabla}_0(\bar{\Rho})) + (\lo,\Hess(H)) + H(\lo,\Rho) - \tfrac{n-3}{n-2} H (\lo,\W) \notag \\
    & + \tfrac{n}{(n-1)^2} |\delta(\lo)|^2 - \bar{\Rho}_{00} |\lo|^2 - \tfrac{1}{(n-3)(n-6)} \J |\lo|^2
    - \tfrac{3}{2} H^2 |\lo|^2 - \tfrac{n-3}{n-2} H \tr (\lo^3) + \tfrac{n-4}{n-3} (\lo^2,\Rho) \notag \\
    & - \tfrac{1}{n-3} \delta \delta (\lo^2) + \tfrac{n-5}{2(n-3)(n-6)} \Delta (|\lo|^2)
\end{align*}
are local conformal invariants of weight $-4$ of an embedding $M^n \hookrightarrow
X^{n+1}$, i.e., it holds 
$$
   e^{4 \iota^*(\varphi)} \hat{\Jv}_j = \Jv_j
$$ 
for $j=1,2$ (Proposition \ref{J-12}). Both invariants have a simple formal pole at $n=3$ with
residue $-\mathcal{D}((\lo^2))_\circ)$. Here $\mathcal{D}: b \mapsto \delta \delta
(b) + (\Rho,b)$ is a conformally covariant operator $S_0^2(M) \to C^\infty(M)$ on
trace-free symmetric $2$-tensors on $M^3$ and $(\lo^2)_\circ$ denotes the trace-free
part of $\lo^2$. Note that the term $\mathcal{D}((\lo^2)_\circ)$ contributes to the
singular Yamabe obstruction $\B_3$ of $M^3 \hookrightarrow X^4$ \cite{GGHW},
\cite{JO2}. We also note that both invariants $\Jv_i$ have a simple formal pole at
$n=6$ with residues being proportional to the local invariant $P_2 (|\lo|^2)$ of
weight $-4$.

Now we return to the critical dimension $n=4$. The above results imply the following
decomposition of the critical extrinsic $Q$-curvature of order $4$ in terms of local
conformal invariants.

\begin{thm}\label{alex} In the critical dimension $n=4$, the extrinsic $Q$-curvature ${\bf Q}_4$ admits the
decomposition
\begin{align}\label{Q4-ex1}
    {\bf Q}_4 & = Q_4 + \frac{9}{2} \Iv_4 + \mathcal{C} + 21 \Iv_6 + \frac{33}{2} \Iv_2 - \frac{14}{3} \Iv_1 \notag \\
    & + \delta \delta (\lo^2) + \frac{1}{6} \Delta (|\lo|^2) + 4 \delta (\lo \delta(\lo)) + 3 \delta \delta (\W),
\end{align}
where the local conformal invariants $\Iv_j$ are defined in Section \ref{ECI}. In particular, ${\bf Q}_4(g)$ is a
linear combination of the Pfaffian of $(M,h)$, local conformal invariants of the embedding $M \hookrightarrow X$
and a divergence term, i.e., it holds
\begin{equation}\label{DS-Q4}
    {\bf Q}_4 = a \Pf_4 + \sum_j b_j I_j + \mbox{total divergence}
\end{equation}
with the Pfaffian density $\Pf_4$ and local conformal invariants $I_j$ of the embedding $M \hookrightarrow X$.
\end{thm}


Some further comments are in order.

In the first line of \eqref{Q4-ex1}, all terms except the intrinsic $Q_4$ of $(M,h)$
are local conformal invariants. Likewise, all terms in the second line are total
divergences. The local conformal invariants in \eqref{DS-Q4} are intrinsic and
extrinsic.

The transformation law \eqref{ECTL} shows that the conformal variation of the sum of
the total divergences in the second lone of \eqref{Q4-ex1} is given by the second-order part in
\eqref{EP4-crit}. In fact, a direct calculation confirms that this conformal
variation equals
$$
   14 \delta (\lo^2 d\varphi) - \frac{4}{3} \delta (|\lo|^2 d\varphi) + 6 \delta (\W d\varphi).
$$
In other words, it is natural to view the pair $({\bf P}_4,{\bf Q}_4)$ as the sum of the pairs
\begin{align*}
   (P_4,Q_4) \quad \mbox{and} \quad ( P_4^e, Q_4^e)
\end{align*}
with
$$
    P_4^e \st \delta \left( 14 \lo^2 - \frac{4}{3} |\lo|^2 h + 6 \W\right)  d \quad \mbox{and} \quad
    Q_4^e \st \delta \delta (\lo^2) + \frac{1}{6} \Delta (|\lo|^2) + 4 \delta (\lo \delta(\lo)) + 3 \delta \delta (\W) ,
$$
and a linear combination of the local conformal invariants $\mathcal{C}, \; \Iv_1,
\; \Iv_2, \; \Iv_4, \; \Iv_6$. Both pairs satisfy the {\em same} conformal
transformation law. It is also worth noting that for any linear combination
$\tilde{Q}_4^e$ of $\delta \delta (\lo^2)$, $\Delta(|\lo|^2)$, $\delta(\lo
\delta(\lo))$ and $\delta \delta (\W)$ there is a second-order operator
$\tilde{P}_4^e$ so that the pair $(\tilde{P}_4^e,\tilde{Q}_4^e)$ satisfies the same
conformal transformation law as $(P_4,Q_4)$. So the most interesting and most
complex part of the structure of the pair $({\bf P}_4,{\bf Q}_4)$ is the local
conformal invariant $\mathcal{C}$ or, equivalently, the local conformal invariants
$\Jv_1$ and $\Jv_2$.

Moreover, it turns out that
\begin{equation}\label{JJD}
   \Jv_1 - 2 \Jv_2 = - \frac{4}{3} |\lo|^4 + 3  \tr(\lo^4)  + \lo^{kl} \lo^{ij}  W_{kijl} + 3 (\lo^2,\W) 
   - \frac{1}{2} |\overline{W}_0|^2,
\end{equation}
where the right-hand side is a linear combination of obvious local conformal invariants 
(Corollary \ref{J12-LC}). This identity leads to the following equivalent decomposition of the critical ${\bf Q}_4$.

\begin{cor*}\label{Q4-crit} The critical extrinsic $Q$-curvature ${\bf Q}_4$ admits the decomposition
\begin{align}\label{Q4-ex2}
    {\bf Q}_4 & = Q_4 - 3\Jv_1 + \frac{9}{2} \Iv_4 - \Iv_5 + 18 \Iv_6 + \frac{1}{2} \Iv_7 
    - \frac{10}{3} \Iv_1 + \frac{27}{2} \Iv_2 \notag \\
    & + \delta \delta (\lo^2) +  \frac{1}{6} \Delta (|\lo|^2) + 4 \delta (\lo \delta(\lo)) + 3 \delta \delta (\W).
\end{align}
\end{cor*}

This result finally describes ${\bf Q}_4$ in terms of trivial conformal invariants $\Iv_j$, the 
non-trivial conformal invariant $\Jv_1$, and some divergence terms.

In general dimensions, the extrinsic ${\bf Q}_4$ can be written as the sum of $Q_4$,
a linear combination of the local conformal invariants $\Iv_1$, $\Iv_2$, $\Iv_4$, 
$\Iv_6$, $-4 \Jv_1 + 2 \Jv_2$, the product of $n-4$ with a curvature term ${\bf
E}_4$ and a divergence term (Theorem \ref{Q-deco-I}). ${\bf E}_4$ admits a
continuation to the critical dimension $n=4$ and the conformal variation of 
$\int_M {\bf E}_4 dvol$ equals the divergence term in \eqref{Q4-ex1} or \eqref{Q4-ex2} (Remark \ref{E-var}).
This generalizes the observation that in the explicit formula \eqref{EQ3} for ${\bf
Q}_3$ the conformal variation of $\int_M (\lo,\Rho) dvol$ is given by $\delta \delta
(\lo)$.

The analog of the decomposition \eqref{DS-Q4} for ${\bf Q}_2(g)$ for a surface $M^2
\hookrightarrow X^3$ is obviously true since it is true for $Q_2(h)$. In this case
and for closed $M$, the total integral of ${\bf Q}_2$ is a linear combination of the
Euler characteristic of $M$ and the total integral of the local conformal invariant
$|\lo|^2$. The latter integral is the Willmore energy of $M \hookrightarrow X$.
Similarly, ${\bf Q}_3$ for $M^3 \hookrightarrow X^4$ is a linear combination of the
local conformal invariant $(\lo,\JF)$ and a divergence term.

Formulas \eqref{DS-Q4} and \eqref{Q4-ex2} are analogs of the Deser-Schwimmer
decomposition of global conformal invariants \cite{DS} (established by Alexakis in
the monograph \cite{alex} and a series of papers). 

In \cite{GR}, Graham and Reichert studied the renormalized volume of minimal hypersurfaces in a 
Poincar\'e-Einstein background. This led to new conformally invariant energies. Moreover, 
an explicit formula for such an energy was derived if the boundary of the hypersurface is a four-manifold. 
In Section \ref{GR-functional}, we shall examine the Graham-Reichert energy functional of four-dimensional 
hypersurfaces from the perspective of an analog of the Deser-Schwimmer decomposition. This also 
reveals one further local conformal invariant of $M^4 \hookrightarrow X^5$. It is given by
\begin{equation}\label{CI-3}
     \Jv_3 \st \frac{1}{2} \lo^{ij} \bar{\nabla}_0(\overline{W})_{0ij0} + H (\lo,\W) + (\Rho,\W) - \bar{B}_{00}
     + \frac{1}{2} \delta \delta (\W),
\end{equation} 
where $\bar{B}$ is the Bach tensor of the background metric (Corollary \ref{J3}). Note that $\Jv_3 = 0$ 
if $\overline{W}=0$.\footnote{The authors of \cite{AS} informed us that in a forthcoming paper 
they prove that the integrated invariant $\Jv_3$ is a linear combination of the integrals of the invariants
$\Iv_5$, $\Iv_6$ and $\Iv_7$ (see \eqref{J3-AS}).} 

For other results on extrinsic analogs of the Deser-Schwimmer classification, we refer to \cite{mondino}.

Finally, we return to the operator ${\bf P}_4$ in general dimensions and take a
closer look at its structure for a totally umbilic hypersurface, i.e., if $\lo = 0$.
In this case, Theorem \ref{main} reduces to the following result.

\begin{cor*}\label{PQ4-gen} Assume that $\lo=0$ and $n \ge 4$. Then
\begin{align}\label{P4-final}
   {\bf P}_4  & = P_4 (f) + 4 \frac{n-1}{n-2} \delta(\W d)  \notag\\
   & + \left(\frac{n}{2}-2\right) \frac{2(n-1)}{(n-2)(n-3)} \left(\frac{n-1}{n-2} |\W|^2 - (n-4) (\Rho,\W)
   + \delta \delta (\W) \right).
\end{align}
In particular, it holds ${\bf P}_4 = P_4$ iff $\W=0$. As a consequence,
\begin{equation}\label{Q4-final}
   {\bf Q}_4 = Q_4 + \frac{2(n-1)}{(n-2)(n-3)} \left(\frac{n-1}{n-2} |\W|^2 - (n-4) (\Rho,\W)
   + \delta \delta (\W) \right).
\end{equation}
Thus, in the critical dimension $n=4$, it holds
\begin{equation}\label{Q4-crit-0}
   {\bf P}_4(f) = P_4(f) + 6 \delta(\W df) \quad \mbox{and}
   \quad  {\bf Q}_4 = Q_4 + \frac{9}{2} |\W|^2 + 3 \delta \delta (\W).
\end{equation}
\end{cor*}

The assumption $\lo=0$ is a conformally invariant condition. The conformally
invariant condition $\W=0$ yields another interesting special case. If $\bar{g} =
dr^2 + h_r$ so that $g_+ = r^{-2} \bar{g}$ satisfies $\Ric (g_+) + n (g_+) = 0$ (we
shall refer to this case as the Poincar\'e-Einstein case), then it holds $\lo=0$
{\em and} $\W = 0$. Hence, in this case, formula \eqref{P4-final} shows that ${\bf
P}_4$ reduces to the intrinsic Paneitz operator $P_4$ of $M$. Formula
\eqref{Q4-crit-0} implies that the total integral of the critical ${\bf Q}_4$ equals
$$
     \int_{M^4} \left(Q_4 + \frac{9}{2} |\W|^2 \right) dvol_h
$$
if $\lo=0$. Now the Chern-Gauss-Bonnet formula\footnote{Here $|W|^2 = W_{ijkl} W^{ijkl}$.}
\begin{equation}\label{CGB-4}
   8 \pi^2 \chi(M) = \frac{1}{4} \int_{M^4} |W|^2 dvol_h + \int_{M^4} Q_4 dvol_h
\end{equation}
immediately implies that the total integral of the critical ${\bf Q}_4$ is conformally invariant
if $\lo=0$. In this special case, one easily sees that $f \mapsto P_4(f) + 6 \delta(\W df)$
is a conformally covariant operator: both $P_4$ and $\delta(\W d)$ are conformally
covariant. Also one can directly verify the fundamental transformation law
\eqref{ECTL} for the critical ${\bf Q}_4$-curvature if $\lo=0$. The residue at $n=3$
of the right-hand side of \eqref{Q4-final} is a constant multiple of $2 |\W|^2 +
(\Rho,\W) + \delta \delta (\W)$ and this quantity is a multiple of the singular Yamabe
obstruction $\B_3$ (if $\lo=0$).

We finish this section with a detailed review of the paper.

In Section \ref{notation}, we fix notation and collect basic identities for later
references. Section \ref{Laplace-scattering} briefly recalls the relation between
extrinsic conformal Laplacians and scattering theory for the singular Yamabe metric.
It implies a description of the conformal Laplacians in terms of the asymptotic
expansions of eigenfunctions of the singular Yamabe metric (see \eqref{P-sol}). This
is followed up in Section \ref{sol-inter} by a discussion of the first few terms in
the asymptotic expansion of such eigenfunctions, which suffice to prove Proposition
\ref{P23}. Moreover, Theorem \ref{P4-inter} gives a preliminary formula for ${\bf
P}_4$. This formula describes ${\bf P}_4$ in terms of data written in coordinates for which 
the singular Yamabe metric has a simple normal form. These
coordinates have been used in \cite{JO1} and recently in \cite{CMY}; they differ 
from the adapted coordinates in \cite{JO1}. The displayed formula resembles the formula for the usual
Paneitz operator in terms of a Poincar\'e-Einstein metric. However, to get
the desired formula for ${\bf P}_4$, it remains to evaluate it further and 
perform a conformal change. These delicate tasks are realized in the following sections. 
The most subtle part concerns the study of two normal derivatives of the scalar curvature 
of the background metric. In Section \ref{critical}, we derive a formula for the total integral 
of ${\bf Q}_4$ in the critical dimension $n=4$ (Corollary \ref{Q4-g-int}). Section \ref{TUH} 
treats the operator ${\bf P}_4$ in the totally umbilic case in general dimensions, proving 
Corollary \ref{PQ4-gen}, and Section \ref{general} contains a proof of the formula for 
${\bf P}_4$ in the general case (Theorem \ref{main}). As a corollary, we rederive a 
formula for the singular Yamabe obstruction $\B_3$ of $M^3 \hookrightarrow X^4$ 
by taking a formal residue at  $n=3$ (Corollary \ref{QB}). In the following sections, we look closely
at the structure of ${\bf Q}_4$. In particular, we prove that $\Jv_1$, 
$\Jv_2$ and $\mathcal{C}$ are local conformal invariants (Theorem \ref{LCI}) and 
establish the decomposition in Theorem \ref{alex}. The main technical result is Lemma \ref{V-1}. 
We round off this section with the formulation of an analogous conjectural decomposition of 
${\bf Q}_n$ in general dimensions. Section \ref{GR-functional} contains
similar results related to the conformally invariant functional $\mathcal{E}_{GR}$ introduced in
\cite{GR}. In this context, we find the local conformal invariant $\Jv_3$ (Corollary \ref{J3}). 
Corollary expresses $\mathcal{E}_{GR}$ in terms of the Euler characteristic of $M$ and the integrals 
of the local invariants $\Jv_1$ and $\Iv_j$. For a flat background, we relate in Section \ref{flat-bg} the 
energy functionals defined by $\mathcal{C}$ to a functional of Guven. In Section \ref{SYE},
we derive a formula for the singular Yamabe energy of $M^4 \hookrightarrow X^5$
(Theorem \ref{V4-final}). It requires solving the singular Yamabe problem to
sufficiently high order. The result is proportional to the total integral of ${\bf Q}_4$. 
Since this relation is a special case of a general result, the calculation serves as an 
additional cross-check. In the Appendix, we collect various technical results. 
The first section outlines a direct proof of the conformal covariance of the operator displayed 
in Theorem \ref{main}. Section \ref{ECI} reviews the known trivial and non-trivial extrinsic 
conformal invariants of hypersurfaces $M^4 \hookrightarrow X^5$. In Section \ref{deco-general}, 
we briefly review the roles of Deser-Schwimmer type decompositions for extrinsic conformal invariants 
of hypersurfaces in other mathematical and physical contexts. Section \ref{invariant} 
provides a new proof of the conformal invariance of the basic conformal invariant ${\text Wm}$ 
introduced in \cite{BGW-1}. Section \ref{diff-rel} contains a proof of 
the relation \eqref{JJD}.

Some of the results of this paper were announced in \cite{J-announce}. The results
in the paper \cite{BGW-1} overlap with the current results. Among other
things, the authors of \cite{BGW-1} derive an explicit formula for the extrinsic Paneitz operator 
${\bf P}_4$ in general dimensions.\footnote{In contrast to \cite{BGW-1}, the arguments in the present paper 
do not rely on computer calculations.}  In Section \ref{BGW-1}, we (almost) prove that this formula is equivalent 
to Theorem \ref{main} (see also Remark \ref{BGW-1-compare}).\footnote{However, we do not verify the 
equivalence of the terms which are quartic in $L$ (see \eqref{L-terms}).}

{\em Acknowledgment.} The author is grateful to B. {\O}rsted for numerous discussions.

\section{Notation and basic identities}\label{notation}

All manifolds $X$ are smooth. For a manifold $X$, $C^\infty(X)$ is the space of
smooth functions on $X$. Metrics on $X$ are denoted by $g$. $dvol_g$ is the
Riemannian volume element defined by $g$. Let $\mathfrak{X}(X)$ be the space of
smooth vector fields on $X$. The Levi-Civita connection of $g$ is denoted by
$\nabla_X^g$ or simply $\nabla_X$ for $X \in \mathfrak{X}(X)$ if $g$ is understood.
In these terms, the curvature tensor $R$ of the Riemannian manifold $(X,g)$ is
defined by $R(X,Y)Z =\nabla_X \nabla_Y (Z) - \nabla_Y \nabla_X (Z) -
\nabla_{[X.Y]}(Z)$ for vector fields $X,Y,Z \in \mathfrak{X}(X)$. The components of
$R$ are defined by $R(\partial_i,\partial_j)(\partial_k) = {R_{ijk}}^l \partial_l$.
$\Ric$ and $\scal$ are the Ricci tensor and the scalar curvature of $g$. On a
manifold $(X,g)$ of dimension $n$, we set $2(n-1) \J = \scal$ and define the
Schouten tensor $\Rho$ of $g$ by 
$$
   (n-2)\Rho = \Ric - \J g
$$ 
(if $n \ge 3$). Let $W$ be the Weyl tensor. Then the curvature tensor admits the 
decomposition $R = W - \Rho \owedge g$. We recall that $W$ vanishes in dimension $3$. 
The Cotton tensor $C$ is defined by 
$$
   C_{ijk} = \nabla_k (\Rho)_{ij} - \nabla_j(\Rho)_{ik}.
$$
Then
$$
    (n-3) C_{ijk} = \nabla^l(W)_{lijk}.
$$
Finally, let 
$$
   B_{ij} = \nabla^k(C)_{ijk} + \Rho^{kl} W_{iklj}
$$
be the Bach tensor. These conventions are as in \cite{J1}. 

For a metric $g$ on $X$ and $u \in C^\infty(X)$, let $\grad_g(u)$ be the gradient of
$u$ with respect to $g$ so that $g(\grad_g(u),V) = \langle du,V \rangle$ for all $V
\in \mathfrak{X}(X)$. $g$ defines pointwise scalar products $(\cdot,\cdot)$ and
norms $|\cdot|$ on $\mathfrak{X}(X)$, on forms and on general tensors. $\delta^g$ is
the divergence operator on differential forms or symmetric bilinear forms. On forms,
it coincides with the negative adjoint $-d^*$ of the exterior differential $d$ with
respect to the Hodge scalar product defined by $g$. Let $\Delta_g = \delta^g d$ be
the non-positive Laplacian on $C^\infty(X)$. On the Euclidean space $\R^n$, it
equals $\sum_i \partial_i^2$. In addition, $\Delta$ will also denote the
Bochner-Laplacian (when acting on $L$, say).

A metric $g$ on a manifold $X$ with boundary $M$ induces a metric $h$ on $M$. In
such a setting, we distinguish the curvature quantities of $g$ and $h$ by adding a
bar to those of $g$. The covariant derivative, the curvature tensor and the Weyl
tensor of $(X,g)$ are $\bar{\nabla}$, $\bar{R}$ and $\overline{W}$. Similarly,
$\overline{\Ric}$ and $\overline{\scal}$ are the Ricci tensor and the scalar
curvature of $g$.

The following conventions coincide with those in \cite{JO1, JO2}.

A hypersurface is given by an embedding $\iota: M \hookrightarrow X$. Accordingly,
tensors on $X$ are pulled back by $\iota^*$ to $M$. In practice, we often omit this
pull-back. For a hypersurface $\iota: M \hookrightarrow X$ with the induced metric
$h = \iota^*(g)$ on $M$, the second fundamental form $L$ is defined by $L(X,Y)= - h
(\nabla^g_X(Y), N)$ for vector fields $X, Y \in \mathfrak{X}(M)$ and a unit normal
vector field $\partial_0 = N$. We set $n H = \tr_h(L)$ if $M$ has dimension $n$.
Then $H$ is the mean curvature of $M$. Let $\lo = L - H h$ be the trace-free part of
$L$. Sometimes we identify $L$ with the shape operator $S$ being defined by
$h(X,S(Y)) = L(X,Y)$.

We use metrics, as usual, to raise and lower indices. In particular, we set
$(L^2)_{ij} = L_i^k L_{kj} = h^{lk} L_{il} L_{kj}$ and similarly for higher powers
of $L$. We always apply the Einstein summation convention, i.e., we sum over
repeated indices.

The $1$-form $\overline{\Ric}_{0} \in \Omega^1(M)$ is defined by
$\overline{\Ric}_{0}(X) = \overline{\Ric}(X,\partial_0)$ for $X \in
\mathfrak{X}(M)$. Similarly, we write $b_0$ for the analogous $1$-form defined by a
bilinear form $b$, and we let $\overline{W}_0$ be the $3$-tensor on $M$ with
components $\overline{W}_{ijk0}$, i.e., we always insert the normal vector
$\partial_0$ into the last slot. Moreover, we set $\W_{ij} = \overline{W}_{0ij0}$.
We define $(\lo,\overline{W}_0) \in \Omega^1(M)$ by $\lo^{ij} \overline{W}_{\cdot
ij0}$.

The curvatures of the background metric $g$ on $X$ and the induced metric $h$ on the
hypersurface $M$ are connected through the Gauss equations
\begin{align}
    \iota^* \overline{R} & = R + \frac{1}{2} L \owedge L, \label{Gauss-R}\\
    \iota^* \overline{\Ric} & = \Ric + L^2 - n H L + \bar{\G}, \label{Gauss-Ric} \\
    \iota^* \bar{\J} & = \J + \frac{1}{2(n-1)} |\lo|^2 - \frac{n}{2} H^2 + \bar{\Rho}_{00} \label{Gauss-scalar}
\end{align}
with $\bar{\G}_{ij} \st \overline{R}_{0ij0}$ and the Codazzi-Mainardi equation
\begin{equation}\label{CM-general}
   \nabla_j (L)_{ik} - \nabla_i(L)_{jk} = \overline{R}_{ijk0}.
\end{equation}
Taking traces gives
\begin{equation}\label{CM-trace}
    \delta(L) - n dH = (n-1) \bar{\Rho}_0 \quad \mbox{and} \quad \delta(\lo) - (n-1) dH = (n-1) \bar{\Rho}_0.
\end{equation}
The trace-free part of the Codazzi-Mainardi equation states that
\begin{equation}\label{CM-TF-3}
   \nabla_j(\lo)_{ik} - \nabla_i(\lo)_{jk} + \frac{1}{n-1} \delta(\lo)_j h_{ik} - \frac{1}{n-1} \delta(\lo)_{i} h_{jk}
   = \overline{W}_{ijk0}.
\end{equation}
For any embedding $\iota: M^n \hookrightarrow X^{n+1}$ ($n \ge 3$), the tensor
\begin{equation}\label{FT}
   \JF \st  \iota^* \bar{\Rho} - \Rho + H \lo + \frac{1}{2} H^2 h
\end{equation}
is conformally invariant: $\hat{\JF} = \JF$. The invariance of $\JF$ also follows from the identity
\begin{equation}\label{Fial}
   (n-2) \left( \iota^* \bar{\Rho} - \Rho + H \lo + \frac{1}{2} H^2 h\right) = \lo^2 - \frac{|\lo|^2}{2(n-1)}h + \W
\end{equation}
(\cite[Lemma 6.23.3]{J1}). Following \cite{V} and \cite{GW-LNY}, we refer to $\JF$ as the Fialkow tensor and to 
the equation \eqref{Fial} as the Fialkow equation. Taking the trace in \eqref{Fial}, yields the Gauss equation
\eqref{Gauss-scalar}. The trace-free part $\mathring{\JF}$ of $\JF$ is given by 
\begin{equation}\label{TF-F}
   \mathring{\JF} = \frac{1}{n-2} \lo^2 - \frac{1}{n(n-2)} |\lo|^2 h + \frac{1}{n-2} \W.
\end{equation}
Finally, we show that $\JF$ naturally appears in the Gauss equation for the Weyl tensor. We calculate 
\begin{align*}
    \iota^* \overline{W} & = \iota^* \overline{R} + \iota^* \bar{\Rho} \owedge h \\
    & = R + \frac{1}{2} L \owedge L +  \iota^* \bar{\Rho} \owedge h \\
    & = W - \Rho \owedge h +  \frac{1}{2} L \owedge L +  \iota^* \bar{\Rho} \owedge h \\
    & = W + \frac{1}{2} L \owedge L  + (\JF - H\lo - \frac{1}{2} H^2 h) \owedge h \\
    & = W + \frac{1}{2} \lo \owedge \lo  + \JF \owedge h
\end{align*}
using the Gauss equation \eqref{Gauss-R} and the decompositions of the curvature tensors. 
This proves the Gauss equation
\begin{equation}\label{Gauss-Weyl}
   \iota^* \overline{W} = W + \frac{1}{2} \lo \owedge \lo + \JF \owedge h
\end{equation}
for the Weyl tensor. 


\section{Extrinsic conformal Laplacians and the scattering operator}\label{Laplace-scattering}

There are two different methods to define extrinsic conformal Laplacians ${\bf
P}_N$. In \cite{GW-LNY}, Gover and Waldron defined these operators in terms of
compositions of so-called Laplace-Robin operators. The latter notion has its origin
in conformal tractor calculus. But Laplace-Robin operators are also linked to
representation theory \cite{JO0} and scattering theory \cite{JO1}. From the point of
view developed in \cite{JO1}, the extrinsic conformal Laplacians appear in terms of
so-called residue families as introduced in \cite{J1} in the setting of
Poincar\'e-Einstein metrics. This also provides a natural definition of extrinsic
$Q$-curvatures. Roughly speaking, residue families may be viewed as curved versions
of symmetry breaking operators in representation theory \cite{KS}. One of the main
results in \cite{JO1} states that both approaches define the same operators.

Moreover, the following result states that the extrinsic conformal Laplacians can be
identified with residues of the geometric scattering operator of the singular metric
$\sigma^{-2} \bar{g}$, where $\sigma \in C^\infty(X)$ satisfies the condition
\begin{equation}\label{CSC}
   \scal_{\sigma^{-2} \bar{g}} = - n(n+1)
\end{equation}
(at least asymptotically). In other words, $\sigma$ is a solution of a singular Yamabe problem. 
For more details on the structure of $\sigma$, we refer to Section \ref{SYE}.

\begin{theorem}[{\cite[Theorem 4]{JO1}}]\label{ECL-S} If \eqref{CSC} and $N \in \N$ satisfies
$1 \le N \le n$ and $(n/2)^2-(N/2)^2$ avoids the discrete spectrum of $-\Delta_{\sigma^{-2} \bar{g}}$,
then
$$
    {\bf P}_N \sim \res_{\frac{n-N}{2}} (\Sc(\lambda)).
$$
The operator $\Sc(\lambda)$ is the scattering operator of the Laplacian of the singular 
Yamabe metric $\sigma^{-2} \bar{g}$.
\end{theorem}


Here a comment on the convention of the normalization of ${\bf P}_{2N}$ is in order. In contrast 
to the normalization ${\bf P}_{2N} = (2N-1)!!^2 \Delta^N + LOT$ used in \cite{JO1}, we 
normalize ${\bf P}_{2N}$ so that its leading term is $\Delta^N$. Accordingly, we adapt the 
normalization of the extrinsic $Q$-curvatures. Then ${\bf Q}_{2N}$ reduces to $(-1)^N Q_{2N}$ 
in the Poincar\'e-Einstein case.

Theorem \ref{ECL-S} extends a result of \cite{GZ} for GJMS-operators. A different
perspective was taken in \cite{CMY} by {\em defining} extrinsic conformal Laplacians
through the residues of the scattering operator. However, this paper did not clarify
the relation of these residues to the operators defined in the works of Gover and
Waldron. Since $|d\sigma|^2 = 1$ on $M$, the singular metric $\sigma^{-2}g$ is
asymptotically hyperbolic, and the proof of the above result again rests on results
in scattering theory as developed in \cite{GZ}. The definition of the scattering
operator $\SC(\lambda)$ combines the existence of local asymptotic expansions of
eigenfunctions of the Laplacian $\Delta_{\sigma^{-2} g}$ with the meromorphic
continuation of the global resolvent. However, the residues of the scattering
operator, which are of interest here, can be described only in terms of the local
asymptotic expansion of eigenfunctions.

In order to analyze these expansions, one has to choose suitable coordinates. In
\cite{JO1}, we utilized so-called {\em adapted} coordinates (which are best suited for the
study of residue families). In these coordinates, the metric $\sigma^{-2} \bar{g}$
takes the form $s^{-2} (a(s) ds^2 + h_s)$ with some coefficient $a \in C^\infty(X)$.
On the other hand, one also may use coordinates so that the metric $\sigma^{-2} \bar{g}$ takes
the form $\hat{r}^{-2}(d\hat{r}^2 + h_{\hat{r}})$. Then the metric $\hat{\bar{g}} \st 
d\hat{r}^2 + h_{\hat{r}}$ is conformally related to the original metric $g$, i.e., 
it holds
\begin{equation}\label{hat-change}
   \hat{\bar{g}} = e^{2\omega} \bar{g}
\end{equation}
with some $\omega \in C^\infty(X)$ so that $\iota^*(\omega) = 0$.\footnote{These coordinates 
have been used in \cite[Section 7]{JO1} and also in \cite{CMY}. The latter paper derived formulas for 
${\bf P}_2$ and ${\bf P}_3$ which are equivalent to those in Proposition \ref{P23}.} In these terms,
assume that $u$ satisfies
$$
    -\Delta_{\hat{r}^{-2} (d\hat{r}^2 + h_{\hat{r}})} u = \lambda(n-\lambda) u
$$
and has $f \in C^\infty(M)$ as boundary value. Then $u$ has an asymptotic expansion of the form
$$
    u \sim \sum_{j \ge 0} \hat{r}^{\lambda+j} \T_j (\lambda)(f) 
   + \sum_{j \ge 0} \hat{r}^{n-\lambda+j} \T_j(n-\lambda) \Sc(\lambda)(f)
$$
with meromorphic one-parameter families $\T_j(\lambda)$ of differential operators on $M$. For more 
details, we refer to \cite[Section 7]{JO1}. There are analogous expansions in terms of adapted coordinates. But 
since $\iota^*(\omega) = 0$, the scattering operator does not depend on these coordinates. The meromorphic 
families $\T_N(\lambda)$ has a simple poles at $\lambda = \frac{n-N}{2}$, and it holds
\begin{equation}\label{P-sol}
    {\bf P}_N \sim \res_{\lambda=\frac{n-N}{2}} (\T_N(\lambda))
\end{equation}
(see \cite[Theorem 9.3]{JO1}); recall that ${\bf P}_N$ is self-adjoint.


\section{Solution operators and the proof of Proposition \ref{P23}}\label{sol-inter}

We assume that the metric $g_+ = r^{-2} (dr^2 + h_r)$ has constant scalar curvature $-n(n+1)$. 

In the present section, we derive formulas for the solution operators $\T_j(\lambda)$ for $j \le 4$. 
We apply the results for the metric $\hat{g}_+ = \hat{r}^{-2} \hat{\bar{g}}$ to prove Proposition 
\ref{P23}. Finally, we provide a preliminary formula for ${\bf P}_4$ in terms of the metric $\hat{\bar{g}}$ 
(see \eqref{hat-change}).

For generic $\lambda$, we analyze the sum
$$
     \sum_{j \ge 0}  r^{\lambda+j}  \T_j(\lambda) (f) = \sum_{j \ge 0} r^{\lambda+j} f_j
$$
contributing to the formal asymptotic expansion of an eigenfunction $u$ so that
$$
     \Delta_{g_+}(u) + \lambda(n\!-\!\lambda) u = 0.
$$
We expand the Laplacian in the form
\begin{align*}
   \Delta_{g_+} & = r^2 \partial_r^2 + (1-n) r \partial_r + r^2 \frac{1}{2} \tr (h_r^{-1} h_r') \partial_r + r^2 \Delta_{h_r} \\
   & =  r^2 \partial_r^2 + (1-n) r \partial_r + r^2 \frac{v'}{v}(r) \partial_r + r^2 \Delta_{h_r},
\end{align*}
where
\begin{equation}\label{v-coeff}
   v(r) = dvol_{h_r}/dvol_h = \frac{1}{2} \tr(h_r^{-1} h_r') = 1+ r v_1 + r^2 v_2  + \cdots.
\end{equation}
Here the prime $'$ denotes the derivative in $r$. Let $\bar{\J}$ be defined for the metric $\bar{g} = dr^2 + h_r$.

\begin{lem}\label{key-form} Assume that $g_+ = r^{-2}(dr^2 + h_r)$ is a metric of constant
scalar curvature $-n(n+1)$. Then
$$
    \bar{\J} = - \frac{1}{2r} \tr (h_r^{-1} h_r')
$$
or, equivalently,
\begin{equation}\label{basic}
    r \bar{\J} = - \frac{v'}{v}.
\end{equation}
\end{lem}

\begin{proof} We recall that the transformation law for scalar curvature under the conformal change
$e^{2\varphi} g$ of the metric $g$ on a manifold $M$ of dimension $n$ reads
$$
   e^{2\varphi} \hat{\J} = \J - \Delta_g (\varphi) - \frac{n-2}{2} |d\varphi|_g^2.
$$
We apply this law for $\bar{g} = r^2 g_+$ on $X$. Then
$$
    r^2 \bar{\J} = \J_{g_+} - \Delta_{g_+}(\log r) - \frac{n-1}{2} |d\log r|^2_{g_+}.
$$
Now $\J_{g_+} = - \frac{n+1}{2}$,
$$
    \Delta_{g_+}(\log r) = -n + \frac{1}{2} r \tr(h_r^{-1}h_r')
$$
and $|d \log r|^2 = 1$ imply  the assertion.
\end{proof}

Lemma \ref{key-form} implies
\begin{align*}
   \Delta_{g_+} (r^\lambda f)
   = -\lambda(n-\lambda) r^{\lambda} f - \lambda r^{\lambda+2} \bar{\J} f + r^{\lambda+2} \Delta_{h_r}(f)
\end{align*}
for $f \in C^\infty(M)$. Hence
$$
   (\Delta_{g_+} + \lambda(n-\lambda)) (r^\lambda f) = r^{\lambda+2} (\Delta - \lambda \bar{\J}) (f)
   + r^{\lambda+3} (\Delta' - \lambda \bar{\J}') (f) + r^{\lambda+4} (\Delta'' - \lambda/2 \bar{\J}'') (f) + \cdots,
$$
where we use the expansion $\Delta_{h_r} = \Delta + r \Delta' + r^2 \Delta'' + \cdots$. It follows that the
{\em boundary value} $f \in C^\infty(M)$ of $u$ is free and $f_1=0$. Moreover, we find
$$
   (\Delta - \lambda \bar{\J}) (f)  - (\lambda+2)(n-\lambda-2) f_2 + \lambda(n-\lambda) f_2 = 0,
$$
i.e.,
\begin{equation}\label{sol-2}
    \mathcal{T}_2(\lambda) (f) = f_2 = \frac{1}{2(n-2\lambda-2)} (\Delta - \lambda \bar{\J})(f).
\end{equation}
Next, we get
$$
   (\Delta' - \lambda \bar{\J}')(f) - (\lambda+3)(n-\lambda-3) f_3 + \lambda(n-\lambda) f_3 = 0,
$$
i.e.,
\begin{equation}\label{sol-3}
    \mathcal{T}_3(\lambda)(f) = f_3 = \frac{1}{3(n-2\lambda-3)} (\Delta' - \lambda \bar{\J}')(f).
\end{equation}
Similarly, we find
$$
    (\Delta'' - \lambda/2 \bar{\J}'')(f) + (\Delta - (\lambda+2) \bar{\J}) f_2 = 4 (n-2\lambda-4) f_4.
$$
Hence
\begin{align}\label{sol-4}
    \mathcal{T}_4(\lambda)(f) & =
    \frac{1}{32(\frac{n}{2}-\lambda-1)(\frac{n}{2}-\lambda-2)} \notag \\
    & \times ((\Delta - (\lambda+2) \bar{\J})(\Delta - \lambda \bar{\J}) + 4 (\tfrac{n}{2}-\lambda-1) \Delta''(f)
   - 2\lambda(\tfrac{n}{2}-\lambda-1) \bar{\J}'').
\end{align}

We summarize these results in

\begin{lem}\label{sol} If $r^{-2} (dr^2 + h_r)$ has constant scalar curvature $-n(n+1)$, then $\T_1 = 0$ and
\begin{align*}
   \T_2(\lambda) = \frac{1}{2(n-2\lambda-2)} (\Delta - \lambda \bar{\J}),
   \quad  \T_3(\lambda)  =  \frac{1}{3(n-2\lambda-3)} (\Delta' - \lambda \bar{\J}')
\end{align*}
and
\begin{align*}
    \T_4(\lambda) & = \frac{1}{32(\frac{n}{2}-\lambda-1)(\frac{n}{2}-\lambda-2)} \\
    & \times ((\Delta - (\lambda+2) \bar{\J})(\Delta - \lambda \bar{\J}) + 4 (\tfrac{n}{2}-\lambda-1) \Delta''(f)
   - 2\lambda(\tfrac{n}{2}-\lambda-1) \bar{\J}'').
\end{align*}
\end{lem}

Thus \eqref{P-sol} and Lemma \ref{sol} imply the preliminary formulas
$$
   {\bf P}_2 = - 4 \res_{\lambda=\frac{n}{2}-1}(\mathcal{T}_2) = \Delta - \tfrac{n-2}{2} \hat{\bar{\J}} \quad
   \mbox{and} \quad
   {\bf P}_3 \sim \res_{\lambda=\frac{n-3}{2}}(\mathcal{T}_3) \sim \Delta'-\tfrac{n-3}{2}\hat{\bar{\J}}'.
$$
We recall that in these formulas, the right-hand sides are defined with respect to the metric $\hat{\bar{g}}$. In
order to further evaluate these results, we apply Lemma \ref{key-form}. By expanding the relation \eqref{basic}
into power series in $r$, we obtain identities for the normal derivatives of  $\bar{\J}$ in terms of the volume coefficients
$v_j$. Comparing these formulas with known formulas for the volume coefficients in terms of the metric will
imply formulas for normal  derivatives of $\bar{\J}$ in terms of the metric. First, we note that
\begin{align*}
   \frac{v'}{v} & = v_1 + (2v_2 - v_1^2) r + (3 v_3 - 3 v_1 v_2 + v_1^3) r^2
   + (4 v_4 - 4 v_1 v_3 - 2 v_2^2 + 4 v_1^2 v_2 - v_1^4) r^3 + \cdots.
\end{align*}
Next, we have

\begin{lem}\label{v-coeff-3} In general dimensions, it holds
\begin{align*}
     v_1 & = n H,\\
     2 v_2 & = - \overline{\Ric}_{00} - |\lo|^2 + n(n\!-\!1) H^2 = \overline{\Ric}_{00} + \scal - \overline{\scal},\\
     6 v_3 & = - \bar{\nabla}_0(\overline{\Ric})_{00} + 2 (\lo,\bar{\G}) - (3n\!-\!2) H \overline{\Ric}_{00}
    + 2 \tr(\lo^3) - 3 (n\!-\!2) H |\lo|^2 + n(n\!-\!1)(n\!-\!2) H^3,
\end{align*}
where $\bar{\G}_{ij} \st \bar{R}_{0ij0}$.
\end{lem}

\begin{proof} These formulas coincide with the corresponding terms in the expansion of the volume form in
\cite[Theorem 3.4]{AGV}. Note that this is obvious for $v_1$ and $v_2$ but requires applying the Gauss 
equations
\begin{equation*}
   \overline{\scal}-\scal  = 2 \overline{\Ric}_{00} + |L|^2- n^2 H^2 \quad \mbox{and} 
   \quad \overline{\Ric} - \Ric = \bar{\G} - n H L - L^2
\end{equation*}
for $v_3$. Equivalent formulas can be found in \cite[Section 2]{GG}.
\end{proof}

Note that \eqref{basic} implies $v_1 = 0$. Thus Lemma \ref{v-coeff-3} shows that $H=0$. Therefore, we get
\begin{align*}
   2v_2 = -\bar{\J}, \quad 3 v_3 = - \bar{\J}' \quad \mbox{and} \quad 4 v_4 - 2 v_2^2 = - \frac{1}{2} \bar{\J}'',
\end{align*}
or, equivalently,
\begin{align}\label{vJ}
   2v_2 = -\bar{\J}, \quad  3 v_3  = - \bar{\J}' \quad \mbox{and} \quad 4 v_4  = -\frac{1}{2} \bar{\J}'' + \frac{1}{2} \bar{\J}^2.
\end{align}

Now combining the first relation in \eqref{vJ} with Lemma \ref{v-coeff-3} gives
\begin{equation}\label{J0}
  \overline{\Ric}_{00} + |\lo|^2 = \bar{\J}.
\end{equation}
We compare this relation with the Gauss identity
$$
    2 \overline{\Ric}_{00} = \overline{\scal} - \scal - |\lo|^2
$$
(recall that $H=0$). We get
$$
   2 \bar{\J} - 2 |\lo|^2 = 2n \bar{\J} - 2(n-1) \J - |\lo|^2
$$
Equivalently, we find
\begin{equation}\label{J}
   \bar{\J} = \J - \frac{1}{2(n-1)} |\lo|^2.
\end{equation}
Now \eqref{J} gives
$$
  {\bf P}_2 = \Delta - \frac{n-2}{2} \left(\J - \frac{1}{2(n-1)}|\lo|^2 \right).
$$
This is the first part of \cite[Theorem 5.5]{CMY}. The formula is also in \cite[Proposition 8.5]{GW-LNY}
and \cite[Section 13.11]{JO1}. Hence
$$
   {\bf P}_2 = P_2 + \frac{n-2}{4(n-1)} |\lo|^2
$$
with $P_2$ being the Yamabe operator (note that in our convention, $\Delta$ is negative). This proves the first part
in Proposition \ref{P23}.

Next, we calculate  $\bar{\J}'$.  The second relation in \eqref{vJ} and Lemma \ref{v-coeff-3} imply
$$
    6 v_3 = - \bar{\nabla}_0(\overline{\Ric})_{00} + 2(\lo,\bar{\G})+ 2 \tr(\lo^3) \stackrel{!}{=} - 2 \bar{\J}'.
$$
In order to evaluate this formula, we apply the following result. Let $\bar{G}$ be the Einstein tensor
$\bar{G} = \overline{\Ric} - n \bar{g} \bar{\J} = \overline{\Ric} - \frac{1}{2} \overline{\scal} \bar{g}$ of $\bar{g}$. 
The following result calculates the normal component of the normal derivative of $\bar{G}$.

\begin{lem}\label{Nabla-1G} In general dimensions, it holds
\begin{align}\label{Einstein-d}
   \bar{\nabla}_0 (\bar{G})_{00} = - \delta (\overline{\Ric}_0) -  n H \overline{\Ric}_{00} + (L,\overline{\Ric}).
\end{align}
\end{lem}


\begin{proof} We prove a more general relation which will be important later. The metric $\bar{g}$ takes the 
form $dr^2 + h_r$ in geodesic normal coordinates. The second Bianchi identity implies $2 \delta^{\bar{g}} (\overline{\Ric}) 
= d \overline{\scal}$. Hence
\begin{align*}
   & \bar{\nabla}_0(\overline{\Ric})(\partial_0,\partial_0) \\
   & = \delta^{\bar{g}}(\overline{\Ric})(\partial_0) - \bar{g}^{ij} \bar{\nabla}_{\partial_i}(\overline{\Ric})(\partial_j,\partial_0) \\
   & =  \frac{1}{2} \langle d \overline{\scal},\partial_0 \rangle - \bar{g}^{ij} \partial_i (\overline{\Ric}(\partial_j,\partial_0))
   + \bar{g}^{ij} \overline{\Ric} (\bar{\nabla}_{\partial_i}(\partial_j),\partial_0)
   + \bar{g}^{ij} \overline{\Ric}(\partial_j,\bar{\nabla}_{\partial_i}(\partial_0)) \\
   & = \frac{1}{2} \langle d \overline{\scal},\partial_0 \rangle - h_r^{ij} \partial_i (\overline{\Ric}(\partial_j,\partial_0))
   + h_r^{ij} \overline{\Ric} (\nabla^{h_r}_{\partial_i}(\partial_j) - (L_r)_{ij} \partial_0,\partial_0)
   + h_r^{ij} \overline{\Ric}(\partial_j,\bar{\nabla}_{\partial_i}(\partial_0))  \\
   & = \frac{1}{2} \langle d \overline{\scal},\partial_0 \rangle - \delta^{h_r} (\overline{\Ric}_0)
   - n H_r \overline{\Ric}_{00} + h_r^{ij} \overline{\Ric}(\partial_j,\bar{\nabla}_{\partial_i}(\partial_0))
\end{align*}
on any level surface of $r$. Here $\delta^{h_r}$ denotes the divergence operator for the induced metric
on the level surfaces of $r$. Similarly, $L_r$ and $H_r$ are the second fundamental form and the
mean curvature of these level surfaces. Therefore, using $\bar{\nabla}_{\partial_i}(\partial_0)
= (L_r)_{ia} h_r^{ak} \partial_k$, we obtain
\begin{align*}
   \bar{\nabla}_0 (\bar{G})_{00} & = -\delta^{h_r} (\overline{\Ric}_0)
  - n H_r \overline{\Ric}_{00} + h_r^{ij} h_r^{ak} (L_r)_{ia} \overline{\Ric}_{jk},
\end{align*}
i.e., we have proved the relation
\begin{equation}\label{Bianchi}
    \bar{\nabla}_0(\bar{G})_{00}
    = - \delta^{h_r} (\overline{\Ric}_0) - n H_r \overline{\Ric}_{00} + (L_r,\overline{\Ric})_{h_r}
\end{equation}
on any level surface of $r$. The assertion is the case $r=0$.
\end{proof}

Now Lemma \ref{Nabla-1G} (using $H=0$) implies
$$
   \bar{\nabla}_0(\bar{G})_{00} = - \delta(\overline{\Ric}_0) + (\lo,\overline{\Ric}).
$$
Hence
$$
   \bar{\nabla}_0(\overline{\Ric})_{00} = n \bar{\J} '- \delta(\overline{\Ric}_0) + (\lo,\overline{\Ric}).
$$
It follows that
$$
    - 2 \bar{\J}' =  - n \bar{\J} ' + \delta(\overline{\Ric}_0) - (\lo,\overline{\Ric}) + 2(\lo,\bar{\G})+ 2 \tr(\lo^3),
$$
i.e.,
$$
   (n-2) \bar{\J}' = \delta \delta (\lo) + (\lo,\overline{\Ric}) + 2(\lo,\bar{\G} - \overline{\Ric})+ 2 \tr(\lo^3)
$$
by \eqref{CM-trace}. Thus the Gauss identity $\bar{\G} - \overline{\Ric} = - \Ric - \lo^2$ yields the desired formula
\begin{equation}\label{J-prime}
    (n-2) \bar{\J}' = \delta \delta (\lo) + (\lo,\overline{\Ric}) - 2(\lo,\Ric).
\end{equation}

We summarize these results in

\begin{prop}\label{J-der-12} If $r^{-2}(dr^2 + h_r) = r^{-2} \bar{g}$ has constant scalar curvature $-n(n+1)$,
then $H=0$ and it holds
\begin{align*}
    \bar{\J} = \J - \frac{1}{2(n-1)} |\lo|^2 \quad \mbox{and} \quad
    (n-2) \bar{\J}' = \delta \delta (\lo) + (\lo,\overline{\Ric}) - 2(\lo,\Ric).
\end{align*}
In particular, it holds $\bar{\J} = \J$ and $\bar{\J}' = 0$ if $\lo=0$.
\end{prop}

Now, in order to prove the second part of Proposition \ref{P23}, we have to calculate
$\hat{\bar{\J}}'$, i.e., $\bar{\J}'$ for the metric $\hat{\bar{g}}$. Because of \eqref{J-prime}
it only remains to calculate $\hat{\overline{\Ric}}$. The conformal transformation
law for the Ricci tensor gives
\begin{align*}
   \hat{\overline{\Ric}}_{ij} & = \overline{\Ric}_{ij} - (n-1) \overline{\Hess}_{ij}(\omega) - \bar{\Delta}(\omega) h_{ij}
   - (n-1) |d\omega|^2 h_{ij} \\
   & = \overline{\Ric}_{ij} - (n-1) L_{ij} \partial_0(\omega) - \bar{\Delta}(\omega) h_{ij}
   - (n-1) |d\omega|^2 h_{ij}
\end{align*}
using $\omega = 0$. By $\partial_0(\omega) = - H$, we get
$$
   \lo^{ij} \hat{\overline{\Ric}}_{ij} = \lo^{ij} \overline{\Ric}_{ij} + (n-1) H |\lo|^2.
$$
Now we apply the formula
$$
   \Delta' = [\delta,\mathcal{H}_1]d = \delta (\mathcal{H}_1 d) - \mathcal{H}_1 \delta d
$$
for the metric variation of the Laplacian. Here
$$
   \begin{cases}
   \mathcal{H}_1 = v_1 = 0 & \mbox{on $\Omega^0(M)$}, \\
   \mathcal{H}_1 = v_1 \id - h_{(1)} = - 2 \lo & \mbox{on $\Omega^1(M)$}.
\end{cases}
$$
We recall that $v_1 = 0$. Hence $\Delta' = - 2 \delta (\lo d)$. Thus, using \eqref{J-prime}, we get
$$
   {\bf P}_3 \sim - \delta (\lo d)
   - \frac{n-3}{4(n-2)} (\delta \delta (\lo) + (\lo,\overline{\Ric}) - 2(\lo,\Ric)+ (n-1) H |\lo|^2).
$$
This is the second part of \cite[Theorem 5.5]{CMY}.\footnote{Note that in \cite[Theorem 5.5]{CMY} the sign of
$H$ is misprinted (they use a different sign convention for $H$).} This formula can also be found in 
\cite[Proposition 8.5]{GW-LNY} and \cite[Proposition 13.10.1]{JO1}.

We continue discussing ${\bf P}_4$. Combining the formula for $\T_4(\lambda)$ in
Lemma \ref{sol} for the metric $\hat{g}_+ = \hat{r}^{-2} (d\hat{r}^2 + h_{\hat{r}}) =  \hat{r}^{-2} \hat{\bar{g}}$
with
$$
   {\bf P_4} \sim \res_{\lambda=\frac{n}{2}-2}(\mathcal{T}_4)
$$
implies the following preliminary formula for ${\bf P}_4$. Here we use the notation $\hat{\bar{\J}}$ and $\hat{\bar{\G}}$
for the quantities $\bar{\J}$ and $\bar{\G}$ for the metric $\hat{\bar{g}}$.-

\begin{theorem}\label{P4-inter} The extrinsic Paneitz operator is given by
\begin{align}\label{P4-gen}
   {\bf P}_4 (f) & =\left (\Delta - \frac{n}{2} \hat{\bar{\J}}\right) \left(\Delta - \left(\frac{n}{2}-2\right) \hat{\bar{\J}}\right) (f)
   -2(d \hat{\bar{\J}},df) - 4 \delta ((\hat{h}_{(2)} - \hat{h}_{(1)}^2) df) - (n-4) \hat{\bar{\J}}'' f \notag \\
   & = \Delta^2(f) - \delta ( (n-2)\hat{\bar{\J}} df + 4 \hat{h}_{(2)} df) + 4 \delta(\hat{h}_{(1)}^2 df)
   + \left(\frac{n}{2}-2\right) {\bf Q}_4 f,
\end{align}
where $\hat{h}_{(1)} = 2 \lo$, $\hat{h}_{(2)} = \lo^2 - \hat{\bar{\G}}$ and
\begin{equation}\label{Q4-gen}
   {\bf Q}_4 = \frac{n}{2} \hat{\bar{\J}}^2 - 2 \hat{\bar{\J}}'' - \Delta (\hat{\bar{\J}}).
\end{equation}
In particular, ${\bf P}_4$ is self-adjoint.
\end{theorem}

\begin{proof} The formula
$$
    \res_{\lambda = \frac{n}{2}-2}(\mathcal{T}_4) \sim \left(\Delta - \frac{n}{2} \bar{\J}\right)
   \left(\Delta - \left(\frac{n}{2}-2\right) \bar{\J}\right) + 4 \Delta'' - (n-4) \bar{\J}''
$$
shows that
$$
   {\bf P}_4 = \left(\Delta - \frac{n}{2} \hat{\bar{\J}}\right) \left(\Delta - \left(\frac{n}{2}-2\right) \hat{\bar{\J}}\right)
   + 4 \hat{\Delta}'' - (n-4) \hat{\bar{\J}}''.
$$
Now we apply the general formula
$$
    \Delta'' = ([\delta,\mathcal{H}_2] - \mathcal{H}_1 [\delta,\mathcal{H}_1])d
$$
for the second metric variation of the Laplacian (see the discussion at the end of the section). Here
$$
\begin{cases}
    \mathcal{H}_2 = v_2     & \mbox{on  } \Omega^0(M), \\
    \mathcal{H}_2 = v_2 \id - v_1 h_{(1)} + (h_{(1)}^2 - h_{(2)}) & \mbox{on $\Omega^1(M)$} .
\end{cases}
$$
Since $v_1 = 0$, the variation formula simplifies to
$$
    \Delta'' = [\delta,\mathcal{H}_2] d
$$
with $\mathcal{H}_2 = v_2 \id + (h_{(1)}^2 - h_{(2)})$. But $v_2 = -1/2 \bar{\J}$. This proves the assertion.
\end{proof}

We finish with a discussion of the variation formulas for the Laplacian used in the above proofs.
More precisely, let $h_r = h + r h_{(1)} + r^2 h_{(2)} + \cdots$ be a family of metrics on $M$.
Then $\Delta_{h_r} = \Delta_h + r \Delta'_h + r^2 \Delta''_h + \cdots$. We refer to $\Delta'_h$
and $\Delta_h''$ as to the first and second metric variation of the Laplacian at $h$.

The arguments rest on the following observation. Let $h_0 = h$ and let the Hodge star operators for $h_0$ and
$h_r$ be denoted by $\star_0$ and $\star_r$, respectively. Let $\H(r) \st \star_0^{-1} \star_r$ acting on $\Omega^*(M)$.
Then it holds
$$
   \H(r) = v(r)
$$
as multiplication operators acting on $C^\infty(M)$. Next, we establish an analogous relation for $\H(r)$ acting on
$\Omega^1(M)$. The relation
$$
   h_r (X,Y) = h_0(T_r X,Y) \quad \mbox{for $X,Y \in \mathfrak{X}(M)$}
$$
defines an isomorphism $T_r \in \End(TM)$, i.e., $T_r (\partial_i) = (T_r)^k_i \partial_k$ with $(h_r)_{ij} = (T_r)_i^k h_{kj}$. In other
words, $T_r$ arises by regarding $h_r$ as an endomorphism using $h_0$. Let $T_r^* \in \End(T^*M)$ be its dual.
Then
$$
   h_r(\alpha,\beta) = h_0(\alpha,(T_r^{-1})^* \beta) \quad \mbox{for $\alpha,\beta \in \Omega^1(M)$}.
$$
Let $\I(r) \st (T_r^{-1})^*$ acting on $\Omega^1(M)$.

\begin{lem}\label{H-form} For any metric $h$, it holds
\begin{equation}\label{H-simple}
    \H(r) = v(r) \I(r): \Omega^1(M) \to \Omega^1(M).
\end{equation}
\end{lem}

\begin{proof} On the one hand, we rewrite the defining relation
$$
    \omega \wedge \star_r \eta = h_r(\omega, \eta) dvol_r, \; \omega, \eta \in \Omega^1(M)
$$
for the star-operator $\star_r$ as
$$
   \omega \wedge \star_r \eta = h_0(\omega, \I(r) \eta) v(r) dvol_0.
$$
On the other hand, the defining relation for the star-operator $\star_0$ implies
$$
    \omega \wedge \eta' = h_0(\omega,\star_0^{-1}\eta') dvol_0, \; \eta' \in \Omega^{n-1}(M).
$$
Hence for $\eta'=\star_r \eta$, we obtain
$$
    h_0(\omega,\star_0^{-1} \star_r \eta) dvol_0 = \omega \wedge \star_r \eta =
    h_0 (\omega,\I(r) \eta) v(r) dvol_0.
$$
It follows that
$$
   \star_0^{-1} \star_r = v(r) \I(r).
$$
The proof is complete.
\end{proof}

We expand $\H(r) = \id + r \H_1 + r^2 \H_2 + \cdots$. Then 
$$
   \H(r)^{-1} = \id - r \H_1 + r^2(-\H_2 + \H_1^2) + \cdots.
$$
Now formula \eqref{H-simple} for $\H(r)$ may be used to expand the Laplacian $\Delta_{h_r}$. We write $\Delta_{h_r} = \delta_r d$
with $\delta_r = \star_r^{-1} d \star_r$ acting on $\Omega^1(M)$. Now
$$
   \delta_r = \H(r)^{-1} \star_0^{-1} d \star_0 \H(r) = \H(r)^{-1} \delta_0 \H(r).
$$
Hence
$$
   \delta_r = \delta_0 + r [\delta_0,\H_1] + r^2 ([\delta_0,\H_2] - \H_1 [\delta_0,\H_1]) + \dots.
$$
Therefore,
$$
   \Delta_{h_r} = \delta_r d = \Delta_h + r [\delta, \H_1] d + r^2 ([\delta,\H_2] d - \H_1 [\delta,\H_1] d) + \cdots.
$$
In other words, the first variation of $\Delta$ at $h$ under the perturbation $h_r$ is given by the operator
$$
    \Delta_h' = [\delta, \H_1] d.
$$
Now Lemma \ref{H-form} implies that $\H_1 = v_1 \id - h_{(1)} = \frac{1}{2} \tr_h(h_{(1)}) \id - h_{(1)}$
on $\Omega^1(M)$. This proves the variation formula
\begin{align*}
   (d/dt)|_0 (\Delta_{h + t h_{(1)}}) (u) & = \delta (\H_1 du) - \H_1 \delta du \\
   & = \frac{1}{2} \delta (\tr_h (h_{(1)}) du) - \delta (h_{(1)} du) - \frac{1}{2} \tr(h_{(1)}) \delta d u \\
   & = \frac{1}{2} (d \tr_h(h_{(1)}),du) - (\Hess(u),h_{(1)}) - (\delta(h_{(1)}),du)
\end{align*}
(with scalar products, $\delta$ and $\Hess$ defined by $h$) which is well-known \cite[(1.185]{Besse}. Similarly,
for the second variation, we obtain
$$
   \Delta_h'' = [\delta,\H_2] d - \H_1 [\delta,\H_1] d.
$$
Lemma \ref{H-form} implies that $\H_2 = v_2 \id - v_1 h_{(1)} + (h_{(1)}^2 - h_{(2)})$ on $\Omega^1(M)$.
Note that if $\H_1$ on functions vanishes, then the formula for the second variation reduces to
$\Delta'' = [\delta,\H_2] d$. These arguments establish the formulas used in the proofs of Proposition \ref{J-der-12}
and Theorem \ref{P4-inter}.

\section{Proof of Corollary \ref{Q4-g-int}}\label{critical}

Theorem \ref{P4-inter} shows that the further discussion of ${\bf P}_4$ and ${\bf Q}_4$ requires a good
understanding of the term $\hat{\bar{\J}}''$. Therefore, we next consider the quantity $\bar{\J}''$ if the metric
$r^{-2}(dr^2 + h_r)$ has scalar curvature $-n(n+1)$. Then we apply the results to prove Corollary
\ref{Q4-g-int}.

We first describe the volume coefficient $v_4$ in terms of the background metric and $L$.

\begin{lem}\label{v4-form}  In general dimensions, it holds
\begin{align*}
    24 v_4 & =  - \bar{\nabla}_0^2(\overline{\Ric})_{00}
    + 2 L^{ij} \bar{\nabla}_0 (\bar{R})_{0ij0} - 4 n H \bar{\nabla}_0 (\overline{\Ric})_{00} \notag \\
    & + 3 (\overline{\Ric}_{00})^2 - 2 |\bar{\G}|^2 + 8 n H (L,\bar{\G}) - 8 (L^2,\bar{\G})
    + 6 \overline{\Ric}_{00} (|L|^2 - n^2 H^2) + 24 \sigma_4 (L),
\end{align*}
where $\sigma_4(L)$ is the fourth elementary symmetric function in the eigenvalues of $L$. Equivalently, it holds
\begin{align*}
    24 v_4 & =  -\bar{\nabla}_0^2(\overline{\Ric})_{00} + 2 \lo^{ij} \bar{\nabla}_0 (\bar{R})_{0ij0}
   - (4n\!-\!2) H \bar{\nabla}_0 (\overline{\Ric})_{00} \notag \\
   & + 3 (\overline{\Ric}_{00})^2 - 2 |\bar{\G}|^2 + 8(n\!-\!2) H (\lo,\bar{\G}) - 8 (\lo^2,\bar{\G}) \notag \\
   &  + 6 |\lo|^2 \overline{\Ric}_{00} - 2(n\!-\!1)(3n\!-\!4) H^2 \overline{\Ric}_{00}+ 24 \sigma_4(L).
\end{align*}
\end{lem}

\begin{proof} This is \cite[Lemma 6.7]{JO2}. Its equivalence to \cite[Theorem 3.4]{AGV} follows by combining the
calculation on page 483 of this reference with the Gauss equation. The proofs in these references differ.
\end{proof}

Now we use Lemma \ref{v4-form} to describe the second normal derivative of $\bar{\J}$ if the metric $r^{-2}(dr^2 + h_r)$
has scalar curvature $-n(n+1)$. We recall that this condition implies $H=0$. Thus Lemma \ref{v4-form} gives
\begin{align}\label{h1}
    & 24 v_4 \notag \\ & = -\bar{\nabla}_0^2(\overline{\Ric})_{00} + 2 \lo^{ij} \bar{\nabla}_0 (\bar{R})_{0ij0}
    + 3 (\overline{\Ric}_{00})^2  - 2 |\bar{\G}|^2 - 8 (\lo^2,\bar{\G}) + 6 |\lo|^2 \overline{\Ric}_{00}
    + 24 \sigma_4(\lo).
\end{align}
The following unconditional result generalizes Lemma \ref{Nabla-1G}. It calculates the normal component of the 
second normal derivative of the Einstein tensor $\bar{G} = \overline{\Ric} - n \bar{g} \bar{\J}$ of a general metric 
$\bar{g}$.

\begin{lem}\label{Nabla-2G} In general dimensions, it holds
\begin{align}\label{Einstein-d2}
     \bar{\nabla}_0^2 (\bar{G})_{00} =  & - (n\!+\!1) H \bar{\nabla}_0 (\overline{\Ric})_{00} + 2n H \bar{\J}' \notag \\
    & + 2 (\lo, \nabla (\overline{\Ric}_0))
     - \delta (\bar{\nabla}_0(\overline{\Ric})_{0}) + (\lo, \bar{\nabla}_0(\overline{\Ric})) \notag \\
     & +  H \delta (\overline{\Ric}_0) - (n\!-\!1) (dH,\overline{\Ric}_0) + 2 (\delta(\lo),\overline{\Ric}_0)
     - \delta ((\lo \overline{\Ric})_{0}) \notag \\
     & + |L|^2  \overline{\Ric}_{00} - (L^2, \overline{\Ric}) + (\overline{\Ric}_{00})^2 - (\bar{\G}, \overline{\Ric}).
\end{align}
\end{lem}

\begin{proof} We recall that in normal geodesic coordinates the metric $g$ takes the form $dr^2 +h_r$ with
$h_r = h + 2r L + \cdots$. We also recall the formulas 
\begin{align}\label{BV}
    n H' = - |L|^2 - \overline{\Ric}_{00} \quad \mbox{and} \quad L' = L^2 - \bar{\G}
\end{align}
for the variation of $H$ and $L$ under the normal exponential map \cite{HP}. Here $'$ denotes the derivative
in the variable $r$.  Moreover, let ${\delta}' \st (d/dr)|_0(\delta^{h_r})$. Then
\begin{align}\label{vdelta}
    {\delta}' (\omega) = - 2 (L, \nabla (\omega))_h - 2 (\delta (L),\omega)_h + n (dH,\omega)_h
\end{align}
for $\omega \in \Omega^1(M^4)$ \cite[(1.185]{Besse}. Now differentiating \eqref{Bianchi} implies
\begin{align*}
    \bar{\nabla}_0^2 (\bar{G})_{00} & = - {\delta}' (\overline{\Ric}_0) - \delta (\partial_r(\overline{\Ric}_0))
   - n {H}' \overline{\Ric}_{00} - n H \bar{\nabla}_0(\overline{\Ric})_{00} \\
    & + ({L}', \overline{\Ric}) + L^{ij} \partial_r (\overline{\Ric}_{ij}) - 4 (L^2,\overline{\Ric}),
\end{align*}
where $L' = (d/dr)|_0(L_r)$.  Note that the last term is caused by the derivative of $h_r$. Hence \eqref{vdelta} implies
\begin{align*}
    \bar{\nabla}_0^2 (\bar{G})_{00} & = 2 (L, \nabla (\overline{\Ric}_0)) + 2 (\delta(L), \overline{\Ric}_0)
    - n (dH,\overline{\Ric}_0) \notag \\
    & - \delta (\bar{\nabla}_0(\overline{\Ric})_{0}) - \delta ((L \overline{\Ric})_{0})
    + |L|^2  \overline{\Ric}_{00} + (\overline{\Ric}_{00})^2 - n H \bar{\nabla}_0 (\overline{\Ric})_{00} \notag \\
    & + (L^2, \overline{\Ric}) - (\bar{\G}, \overline{\Ric}) + (L, \bar{\nabla}_0(\overline{\Ric}))
   + 2 (L^2, \overline{\Ric}) - 4 (L^2,\overline{\Ric})
\end{align*}
at $r=0$. Here we used the relations
\begin{align*}
   \bar{\nabla}_0(\overline{\Ric}_0) = \partial_r (\overline{\Ric}_{0})  - (L \overline{\Ric})_{0} \quad \mbox{and} \quad
   \bar{\nabla}_0(\overline{\Ric})_{ij} = \partial_r (\overline{\Ric}_{ij}) - (L \overline{\Ric} + \overline{\Ric} L)_{ij}.
\end{align*}
Now, separating the trace-free part of $L$, we obtain
\begin{align*}
    \bar{\nabla}_0^2 (\bar{G})_{00} & = 2 (\lo,\nabla (\overline{\Ric}_0)) + 2 H \delta (\overline{\Ric}_0)
   + 2 (\delta(\lo),\overline{\Ric}_0) - (n-2)(dH,\overline{\Ric}_0) \\
   & - \delta (\bar{\nabla}_0(\overline{\Ric})_{0}) - \delta ((\lo \overline{\Ric})_{0}) - \delta (H \overline{\Ric}_0) \\
   & + |L|^2  \overline{\Ric}_{00} + (\overline{\Ric}_{00})^2 - n H \bar{\nabla}_0 (\overline{\Ric})_{00} \\
   & - (L^2, \overline{\Ric}) - (\bar{\G}, \overline{\Ric})
   + (\lo, \bar{\nabla}_0(\overline{\Ric})) + H \overline{\scal}' - H \bar{\nabla}_0(\overline{\Ric})_{00}.
\end{align*}
Simplification leads to the result
\begin{align*}
     \bar{\nabla}_0^2 (\bar{G})_{00} =  & - (n+1)  H \bar{\nabla}_0 (\overline{\Ric})_{00} + H \overline{\scal}' \\
    & + 2 (\lo, \nabla (\overline{\Ric}_0))
    - \delta (\bar{\nabla}_0(\overline{\Ric})_{0}) + (\lo, \bar{\nabla}_0(\overline{\Ric})) \\
    & +  H \delta (\overline{\Ric}_0) - (n-1) (dH,\overline{\Ric}_0) + 2 (\delta(\lo),\overline{\Ric}_0)
    - \delta ((\lo \overline{\Ric})_{0}) \\
    & + |L|^2  \overline{\Ric}_{00} - (L^2, \overline{\Ric}) + (\overline{\Ric}_{00})^2 - (\bar{\G}, \overline{\Ric}).
\end{align*}
The proof is complete.
\end{proof}

This result is an extension of \cite[Lemma 6.12]{JO2}. Note that in the second formula, only the coefficients of
$H \bar{\nabla}_0 (\overline{\Ric})_{00}$, $(dH,\overline{\Ric}_0)$ and $H \bar{\J}'$ depend on the dimension of $M$.

Now we apply Lemma \ref{Nabla-2G} for the background metric $\bar{g} = dr^2 + h_r$ so that $g_+ = r^{-2} \bar{g}$
has scalar curvature $-n(n+1)$. We obtain
\begin{align}\label{h2}
     \bar{\nabla}_0^2 (\bar{G})_{00} & = 2 (\lo, \nabla (\overline{\Ric}_0)) - \delta (\bar{\nabla}_0(\overline{\Ric})_{0})
     + (\lo, \bar{\nabla}_0(\overline{\Ric})) \notag \\
     & + 2 (\delta(\lo),\overline{\Ric}_0) - \delta ((\lo \overline{\Ric})_{0})
    + |\lo|^2  \overline{\Ric}_{00} - (\lo^2, \overline{\Ric}) + (\overline{\Ric}_{00})^2 - (\bar{\G}, \overline{\Ric})
\end{align}
using $H=0$. Now combining \eqref{h1} and \eqref{h2} gives
\begin{align*}
   24 v_4 & = - n \bar{\J}'' - 2 (\lo, \nabla (\overline{\Ric}_0)) + \delta (\bar{\nabla}_0(\overline{\Ric})_{0})
    - (\lo, \bar{\nabla}_0(\overline{\Ric})) \notag \\
    & - 2 (\delta(\lo),\overline{\Ric}_0) + \delta ((\lo \overline{\Ric})_{0})
    - |\lo|^2  \overline{\Ric}_{00} + (\lo^2, \overline{\Ric}) - (\overline{\Ric}_{00})^2 + (\bar{\G}, \overline{\Ric}) \\
   & + 2 \lo^{ij} \bar{\nabla}_0 (\bar{R})_{0ij0}
   + 3 (\overline{\Ric}_{00})^2 - 2 |\bar{\G}|^2 - 8 (\lo^2,\bar{\G}) + 6 |\lo|^2 \overline{\Ric}_{00} + 24 \sigma_4(\lo).
\end{align*}
By $24 v_4 = -3 \bar{\J}'' + 3 \bar{\J}^2$ (see \eqref{vJ}), this result implies the formula
\begin{align}\label{J-dprime}
   (n-3) \bar{\J}'' & = - 3 \bar{\J}^2 - 2 (\lo, \nabla (\overline{\Ric}_0)) + \delta (\bar{\nabla}_0(\overline{\Ric})_{0})
    - (\lo, \bar{\nabla}_0(\overline{\Ric})) \notag \\
    & - 2 (\delta(\lo),\overline{\Ric}_0) + \delta ((\lo \overline{\Ric})_{0})
   + (\lo^2, \overline{\Ric}) - (\overline{\Ric}_{00})^2 + (\bar{\G}, \overline{\Ric}) \notag \\
   & + 2 \lo^{ij} \bar{\nabla}_0 (\bar{R})_{0ij0}
   + 3 (\overline{\Ric}_{00})^2 - 2 |\bar{\G}|^2 - 8 (\lo^2,\bar{\G}) + 5 |\lo|^2 \overline{\Ric}_{00} + 24 \sigma_4(\lo).
\end{align}

We summarize these results in

\begin{prop}\label{J-der} If $g_+ = r^{-2}(dr^2 + h_r) = r^{-2} \bar{g}$ has constant scalar curvature $-n(n+1)$,
then
\begin{align*}
   (n-3) \bar{\J}'' & = - 3 \bar{\J}^2 - 2 (\lo, \nabla (\overline{\Ric}_0)) + \delta (\bar{\nabla}_0(\overline{\Ric})_{0})
    - (\lo, \bar{\nabla}_0(\overline{\Ric})) \notag \\
    & - 2 (\delta(\lo),\overline{\Ric}_0) + \delta ((\lo \overline{\Ric})_{0})
   + (\lo^2, \overline{\Ric}) + (\bar{\G}, \overline{\Ric}) + 2 \lo^{ij} \bar{\nabla}_0 (\bar{R})_{0ij0}  \notag \\
   & + 2 (\overline{\Ric}_{00})^2 - 2 |\bar{\G}|^2 - 8 (\lo^2,\bar{\G}) + 5 |\lo|^2 \overline{\Ric}_{00} + 24 \sigma_4(\lo).
\end{align*}
In particular, it holds
\begin{align*}
     (n-3) \bar{\J}'' & = - 3 \J^2 + (\bar{\G}, \overline{\Ric}) + 2 (\overline{\Ric}_{00})^2
   - 2 |\bar{\G}|^2 + \delta (\bar{\nabla}_0(\overline{\Ric})_{0})
\end{align*}
if $\lo=0$.
\end{prop}

Now partial integration shows

\begin{cor}\label{J-2prime} Let $M$ be closed. Then
\begin{align}\label{J-2p-gen}
     & (n-3) \int_M \bar{\J}'' dvol_h \notag\\
     & =  \int_M \left(-3 \bar{\J}^2 + (\bar{\G}, \overline{\Ric}) + 2 (\overline{\Ric}_{00})^2
    - 2 |\bar{\G}|^2 + (\lo^2,\overline{\Ric}) + 5 |\lo|^2 \overline{\Ric}_{00} - 8 (\lo^2,\bar{\G})
    + 24 \sigma_4(\lo) \right) dvol_h \notag \\
    & + \int_M \left(- (\lo, \bar{\nabla}_0(\overline{\Ric})) + 2 \lo^{ij} \bar{\nabla}_0 (\bar{R})_{0ij0}\right) dvol_h.
\end{align}
In particular, it holds
\begin{equation}\label{S0}
   (n-3)  \int_{M} \bar{\J}'' dvol_h = \int_{M} \left(-3 \J^2 + (\bar{\G}, \overline{\Ric}) + 2 (\overline{\Ric}_{00})^2
   - 2 |\bar{\G}|^2 \right) dvol_h.
\end{equation}
if $\lo=0$.
\end{cor}

Corollary \ref{J-2prime} shows that $\bar{\J}''$ and its integral substantially simplify under the assumption $\lo=0$.

Now, to further evaluate the integrals in Corollary \ref{J-2prime}, we will derive formulas for some of the ingredients.

\begin{lem}\label{C-CSC} If $r^{-2}(dr^2 + h_r) = r^{-2} \bar{g}$ has constant scalar curvature $-n(n+1)$,
then it holds
\begin{align*}
    \overline{\Ric}_{00} & = \J - \frac{2n-1}{2(n-1)}|\lo|^2 \\
    \bar{\G} & = \Rho + \frac{1}{n-2} \lo^2 -\frac{2n-3}{2(n-1)(n-2)} |\lo|^2 h + \frac{n-1}{n-2} \W.
\end{align*}
In particular, it holds
$$
   \overline{\Ric}_{00} = \J \quad \mbox{and} \quad \bar{\G} = \Rho + \frac{n-1}{n-2} \W
$$
if $\lo = 0$.
\end{lem}

\begin{proof} By the Gauss equation, it holds
$$
   2 \overline{\Ric}_{00} = 2n \bar{\J} - 2(n-1) \J - |\lo|^2 - n(n-1) H^2.
$$
Using \eqref{J}, we get
$$
   2 \overline{\Ric}_{00} = 2n \J - 2(n-1) \J - \frac{n}{n-1} |\lo|^2 - |\lo|^2 - n(n-1) H^2.
$$
This proves the first relation using $H=0$. Next, we evaluate the decomposition
$\bar{\G} = \bar{\Rho} + \bar{\Rho}_{00} h + \W$ under the assumption $H=0$. We apply the Fialkow
equation
$$
   \overline{\Rho} = \Rho + \frac{1}{n-2} \left(\lo^2 - \frac{1}{2(n-1)} |\lo|^2 h + \W\right)
$$
(see \eqref{Fial} for $H=0$) and
\begin{align*}
    \bar{\Rho}_{00} & = \frac{1}{n-1} (\overline{\Ric}_{00} - \bar{\J}) = - \frac{1}{n-1} |\lo|^2
\end{align*}
(see \eqref{J0}). Hence
\begin{align*}
    \bar{\G} & = \Rho + \frac{1}{n-2} \left(\lo^2 - \frac{1}{2(n-1)} |\lo|^2 h + \W\right)
    - \frac{1}{n-1} |\lo|^2 h + \W \\
    & = \Rho +  \frac{1}{n-2} \lo^2 - \frac{2n-3}{2(n-1)(n-2)} |\lo|^2 h + \frac{n-1}{n-2} \W.
\end{align*}
The proof is complete.
\end{proof}

Now we can give the

\begin{proof}[Proof of Corollary \ref{Q4-g-int}]
We calculate the integrals in \eqref{J-2p-gen} for the metric $\hat{\bar{g}} = e^{2\omega} \bar{g}$. We start with a
discussion of the last two terms. First, we note that
$$
   (\lo,\bar{\nabla}_0(\overline{\Ric})) = (n-1) (\lo,\bar{\nabla}_0 (\bar{\Rho})) \quad \mbox{and} \quad
    \lo^{ij} \bar{\nabla}_0 (\bar{R})_{0ij0} = \lo^{ij} \bar{\nabla}_0(\bar{\Rho})_{ij} + \lo^{ij} \bar{\nabla}_0(\W)_{ij}.
$$
Hence
$$
    \int_M \left(- (\lo, \bar{\nabla}_0(\overline{\Ric})) + 2 \lo^{ij} \bar{\nabla}_0 (\bar{R})_{0ij0}\right) dvol_h
    = \int_M \left(-(n-3) (\lo, \bar{\nabla}_0 (\bar{\Rho}))
   + 2 \lo^{ij} \bar{\nabla}_0(\W)_{ij}\right) dvol_h.
$$
We will apply these terms for the metric $\hat{\bar{g}} = e^{2\omega} \bar{g}$. First, we note that
$$
   \hat{\bar{\nabla}}_0(\hat{\W})_{ij} = \bar{\nabla}_0(\W)_{ij} + 2 H \W_{ij}.
$$
Second, the results in Section \ref{C} on the conformal variation of $(\lo,\bar{\nabla}_0(\bar{\Rho}))$
(applied to the conformal factor $e^{2\omega}$) (see \eqref{CV-NRho}) show that
\begin{align*}
   \int_M (\lo,\hat{\bar{\nabla}}_0 & (\hat{\bar{\Rho}})) dvol_h = \int_M (\lo,\bar{\nabla}_0(\bar{\Rho})) dvol_h \\
   & + \int_M \left( (\lo,\Hess(H)) + 2 H(\lo,\bar{\Rho}) - H (\lo,\bar{\G})  - 2 H^2 |\lo|^2 - H \tr(\lo^3)
   - \omega'' |\lo|^2 \right) dvol_h \\
   & + \int_M 2 H^2 |\lo|^2 dvol_h.
\end{align*}
Here we utilized the properties $\omega = 0$ and $\partial_0(\omega) = - H$ on $M$. The term in the last line is
caused by the non-linear contributions of $\omega$ (see Remark \ref{NL}). Hence formula \eqref{J-2p-gen}
for $\hat{\bar{g}}$ reads
\begin{align}\label{J-double-prime}
    &  (n-3) \int_M \hat{\bar{\J}}'' dvol_h \notag \\
    & = \cdots + \int_M \left(- (\lo,\hat{\bar{\nabla}}_0 (\hat{\overline{\Ric}}))
    + 2 \lo^{ij} \hat{\bar{\nabla}}_0 (\hat{\bar{R}})_{0ij0}\right) dvol_h \notag \\
    & = \cdots + \int_M \left(-(n-3) (\lo, \hat{\bar{\nabla}}_0 (\hat{\bar{\Rho}}))
   + 2 \lo^{ij} \hat{\bar{\nabla}}_0(\hat{\W})_{ij}\right) dvol_h \notag \\
    & = \cdots + \int_M \left(- (n-3) (\lo,\bar{\nabla}_0(\bar{\Rho})) + 2 \lo^{ij} \bar{\nabla}_0(\W)_{ij}
    + 4 H (\lo,\W)\right) dvol_h
   \notag \\
    & - (n-3)  \int_M \left((\lo,\Hess(H)) + 2 H(\lo,\bar{\Rho}) - H (\lo,\bar{\G}) - H \tr(\lo^3) - \omega'' |\lo|^2\right) dvol_h.
\end{align}
Now we restrict to $n=4$ and find
\begin{align}\label{diff}
    \hat{\bar{\J}}^2 - \hat{\bar{\J}}''  & =
    (\lo, \bar{\nabla}_0 (\bar{\Rho})) - 2 \lo^{ij} \bar{\nabla}_0(\W)_{ij} - 4 H(\lo,\W) \\
    & + (\lo,\Hess(H)) + 2 H(\lo,\bar{\Rho}) - H (\lo,\bar{\G}) - H \tr(\lo^3) - \omega'' |\lo|^2 \notag \\
    & + \hat{\bar{\J}}^2 + 3 \hat{\bar{\J}}^2 - (\hat{\bar{\G}}, \hat{\overline{\Ric}}) - 2 (\hat{\overline{\Ric}}_{00})^2
    + 2 |\hat{\bar{\G}}|^2 -  (\lo^2,\hat{\overline{\Ric}}) - 5 |\lo|^2 \hat{\overline{\Ric}}_{00} + 8 (\lo^2,\hat{\bar{\G}})
   - 24 \sigma_4(\lo), \notag
\end{align}
up to a total divergence. By \cite[Lemma 5.4]{CMY}, the term $\omega''$ equals
$$
    \omega'' = \frac{n+1}{2} H^2 + \bar{\J} - \J + \frac{1}{2(n-1)}|\lo|^2
    = \frac{1}{2} H^2 + \bar{\Rho}_{00} + \frac{1}{n-1} |\lo|^2
$$
using the Gauss equation
$$
   \bar{\J} - \J = \bar{\Rho}_{00} - \frac{n}{2} H^2 + \frac{1}{2(n-1)} |\lo|^2.
$$
The terms in the last line of \eqref{diff} do not depend on $H$ and can be expressed in terms of the original metric
using the formulas derived above. Here we apply the formulas
$$
   \hat{\bar{\J}} = \J - \frac{1}{2(n-1)} |\lo|^2, \quad \hat{\overline{\Ric}}_{00} = \J - \frac{2n-1}{2(n-1)} |\lo|^2
$$
(see \eqref{J} and Lemma \ref{C-CSC}), the decomposition
\begin{equation}\label{Ric-0}
   \hat{\overline{\Ric}} = \Ric + \hat{\bar{\G}} + \lo^2 = (n-2) \Rho + \J h + \hat{\bar{\G}} + \lo^2
\end{equation}
(by the Gauss identity for the Ricci tensor), and
\begin{equation}\label{G-0}
    \hat{\bar{\G}} = \Rho + \frac{1}{n-2} \lo^2 -\frac{2n-3}{2(n-1)(n-2)} |\lo|^2 h + \frac{n-1}{n-2} \W
\end{equation}
(see Lemma \ref{C-CSC}). Recall that $H=0$ for the metric $\hat{\bar{g}}$. The terms in the second line of \eqref{diff}
do not depend on $\J$. For their simplification, we apply the decomposition
$
   \bar{\G} = \bar{\Rho} + \bar{\Rho}_{00} h + \W
$
and the Fialkow equation \eqref{Fial}. Note also that $24 \sigma_4(\lo) = 3 |\lo|^4 - 6 \tr(\lo^4)$. Then
a calculation yields
\begin{align*}
    \hat{\bar{\J}}^2 - \hat{\bar{\J}}''  & = \J^2 - |\Rho|^2 + \frac{9}{4} |\W|^2 \\
   & + (\lo, \bar{\nabla}_0 (\bar{\Rho})) - 2 \lo^{ij} \bar{\nabla}_0(\W)_{ij} + (\lo, \Hess(H)) + H(\lo,\Rho)
   - \frac{9}{2} H(\lo,\W) \\
   & + \frac{21}{2}(\lo^2,\W) + 4 (\lo^2,\Rho) - |\lo|^2 \bar{\Rho}_{00}
   -\frac{3}{2} \J |\lo|^2 - \frac{3}{2} H^2 |\lo|^2  - \frac{1}{2} H \tr (\lo^3) \\
   & - \frac{7}{3} |\lo|^4 + \left(\frac{9}{4} + 6\right) \tr (\lo^4),
\end{align*}
up to a divergence term. This proves Corollary \ref{Q4-g-int}.
\end{proof}

\section{${\bf P}_4$ and ${\bf Q}_4$ for totally umbilic hypersurfaces}\label{TUH}

In the present section, we prove Corollary \ref{PQ4-gen}. Although this result is an obvious consequence of 
Theorem \ref{main}, the following proof prepares the proof of Theorem \ref{main} in Section \ref{general}.

The following result makes the right-hand side of the last relation in Proposition \ref{J-der} explicit.

\begin{lem}\label{J-N2} Assume that $r^{-2}(dr^2 + h_r) = r^{-2} \bar{g}$ has constant scalar curvature $-n(n+1)$.
If $\lo=0$ then
\begin{align*}
    (n-3) \bar{\J}'' = (n-3) |\Rho|^2 - \frac{(n-1)^2}{(n-2)^2} |\W|^2 + \frac{(n-4)(n-1)}{(n-2)} (\Rho,\W)
   + \delta (\bar{\nabla}_0(\overline{\Ric})_{0}).
\end{align*}
In particular, it holds
$$
   \bar{\J}'' = |\Rho|^2 - \frac{9}{4} |\W|^2 + \delta (\bar{\nabla}_0(\overline{\Ric})_{0})
$$
in the critical dimension $n=4$.
\end{lem}

\begin{proof} This is a direct calculation using \eqref{Ric-0}, \eqref{G-0} and $\overline{\Ric}_{00} = \J$.
\end{proof}

Now
$$
   {\bf Q}_4 = \frac{n}{2} \hat{\bar{\J}}^2 - 2 \hat{\bar{\J}}'' - \Delta (\hat{\bar{\J}})
$$
(see \eqref{Q4-gen}) implies

\begin{cor}\label{Q4-inter} The extrinsic $Q$-curvature ${\bf Q}_4$ of a totally umbilic hypersurface is given by
\begin{align}\label{Q4-L0}
    {\bf Q}_4 & = \frac{n}{2}\J^2 - 2 |\Rho|^2 - \Delta (\J) \notag \\
    & + \frac{2}{n-3} \frac{(n-1)^2}{(n-2)^2} |\W|^2 - 2\frac{(n-4)(n-1)}{(n-3)(n-2)} (\Rho,\W)
     - \frac{2}{n-3} \delta (\hat{\bar{\nabla}}_0(\hat{\overline{\Ric}})_{0}).
\end{align}
Note that the first three terms define the intrinsic $Q$-curvature $Q_4$ of $h$.
In particular, in the critical dimension $n=4$, this formula reads
\begin{align}\label{Q4-L0-crit}
    {\bf Q}_4 = Q_4 + \frac{9}{2} |\W|^2 - 2 \delta (\hat{\bar{\nabla}}_0(\hat{\overline{\Ric}})_{0}).
\end{align}
Hence
$$
    \int_{M^4} {\bf Q}_4 dvol_h = \int_{M^4} \left(Q_4 + \frac{9}{2} |\W|^2\right) dvol_h
$$
if $M$ is closed.
\end{cor}

It remains to evaluate the divergence term $\delta (\hat{\bar{\nabla}}_0(\hat{\overline{\Ric}})_{0})$
in terms of the original metric instead of $\hat{r}^2 \sigma^{-2} \bar{g}$. This will be done below. The
general case is much more complicated due to a complicated structure of $\hat{\bar{\J}}''$. It will be
discussed in the next section.

In order to calculate the divergence term $\delta (\hat{\bar{\nabla}}_0(\hat{\overline{\Ric}})_{0})$, we
apply the properties
\begin{itemize}
\item $\omega = 0$ on $M$,
\item $\partial_0(\omega) = - H$ on $M$,
\item $\partial_0^2(\omega) = (n+1)/2 H^2 + \bar{\J} - \J$ on $M$
\end{itemize}
(if $\lo=0$) \cite[Lemma 5.4]{CMY}.

\begin{lem}\label{div-ex} Assume that $\lo=0$. Then in any dimensions
\begin{equation*}
    \delta (\hat{\bar{\nabla}}_0(\hat{\bar{\Rho}})_0) = \delta(\bar{\nabla}_0(\bar{\Rho})_0)
    - \Delta (\bar{\Rho}_{00}+ H^2).
\end{equation*}
\end{lem}

\begin{proof} In the following, we work in general dimensions. We first note that
$$
    \hat{\bar{\nabla}}_0(\hat{\overline{\Rho}})_{0a}
    = \hat{\bar{\nabla}}_0 \left(\bar{\Rho} - \overline{\Hess}(\omega) + d\omega \otimes d\omega
    - \frac{1}{2} |d\omega|^2 \bar{g}\right)_{0a}.
$$
Now
\begin{align*}
   \hat{\bar{\nabla}}_0(\bar{\Rho})_{0a} & = \partial_0 (\bar{\Rho}_{0a}) - \bar{\Rho}(\hat{\bar{\nabla}}_0 (\partial_0), \partial_a)
    - \bar{\Rho}(\partial_0,\hat{\bar{\nabla}}_0(\partial_a)) \\
    & = \partial_0 (\bar{\Rho}_{0a}) - \bar{\Rho}(\bar{\nabla}_0(\partial_0)
   + 2 \partial_0(\omega) \partial_0 - \overline{\grad}(\omega),\partial_a) - \bar{\Rho}(\partial_0,\bar{\nabla}_0(\partial_a)
   + \partial_0(\omega) \partial_a + \partial_a(\omega) \partial_0) \\
   & = \bar{\nabla}_0 (\bar{\Rho})_{0a} - 2 \partial_0(\omega) \bar{\Rho}_{0a} + \bar{\Rho}(\overline{\grad}(\omega),\partial_a)
   - \partial_0(\omega) \bar{\Rho}_{0a} - \partial_a(\omega) \bar{\Rho}_{00} \\
   & =  \bar{\nabla}_0 (\bar{\Rho})_{0a} + 2 H \bar{\Rho}_{0a}.
\end{align*}
Similarly, we find
\begin{align*}
     \hat{\bar{\nabla}}_0 (\overline{\Hess}(\omega))_{0a}
     & =  \bar{\nabla}_0 (\overline{\Hess}(\omega))_{0a} + 2 H \overline{\Hess}(\omega)_{0a}, \\
     \hat{\bar{\nabla}}_0 (d\omega \otimes d\omega)_{0a} & = \bar{\nabla}_0 (d\omega \otimes d\omega)_{0a}, \\
     \hat{\bar{\nabla}}_0 (|d\omega|^2 \bar{g})_{0a} & =  \bar{\nabla}_0 ((|d\omega|^2 \bar{g})_{0a}  = 0
\end{align*}
using $\omega = 0$ on $M$. Hence
$$
    \hat{\bar{\nabla}}_0(\hat{\overline{\Rho}})_{0a} = \bar{\nabla}_0(\overline{\Rho})_{0a}
    - \bar{\nabla}_0(\overline{\Hess}(\omega))_{0a} + \bar{\nabla}_0 (d\omega \otimes d\omega)_{0a}
    + 2 H \bar{\Rho}_{0a}  - 2 H \overline{\Hess}(\omega)_{0a}.
$$
Now
\begin{equation*}
    \bar{\nabla}_0 (d\omega \otimes d\omega)_{0a} = \partial_0 (\partial_0(\omega) \partial_a(\omega))
    - (d\omega \otimes d\omega)(\partial_0, \bar{\nabla}_0(\partial_a)) = \partial_0(\omega) \partial^2_{0a}(\omega)
     = H \partial_a(H)
\end{equation*}
using $\partial_0(\omega) = -H$, and
\begin{align*}
    \bar{\nabla}_0(\overline{\Hess}(\omega))_{0a} & = \partial_0 (\overline{\Hess}(\omega)_{0a})
     - \overline{\Hess}(\omega)(\partial_0,\bar{\nabla}_0(\partial_a)) \\
    & = \partial_0 (\overline{\Hess}(\omega)_{0a}) - L_a^b \overline{\Hess}(\omega)_{0b} \\
    & = \partial_0 (\overline{\Hess}(\omega)_{0a}) + H \partial_a(H)
\end{align*}
using
\begin{equation*}
   \overline{\Hess}(\omega)_{0a} = \partial_{0a}^2(\omega) - \bar{\Gamma}_{0a}^b \omega_b
   = \partial_{0a}^2(\omega)  = - \partial_a(H).
\end{equation*}
But
\begin{align*}
   \partial_0 (\overline{\Hess}(\omega)_{0a}) & = \partial_0 ( \partial_{0a}^2(\omega) - \bar{\Gamma}_{0a}^b \omega_b) \\
   & = \partial_a \partial_0^2(\omega) - \bar{\Gamma}_{0a}^b \partial_{0b}^2(\omega) \\
   & = \partial_a \partial_0^2(\omega) + L_a^b \partial_b(H) \\
   & = \partial_a \partial_0^2(\omega) + H \partial_a(H).
\end{align*}
We combine these results with the formula $\partial_0^2(\omega) = \frac{n+1}{2} H^2 + \bar{\J} - \J$. Then
$$
   \partial_0 (\overline{\Hess}(\omega)_{0a}) = \partial_a \left(\frac{n+1}{2} H^2 
   + \bar{\J} - \J \right) + H \partial_a(H),
$$
i.e.,
$$
   \partial_0 (\overline{\Hess}(\omega)_{0}) = d \left(\frac{n+1}{2} H^2 + \bar{\J} - \J \right) + H dH
$$
Hence
$$
    \bar{\nabla}_0(\overline{\Hess}(\omega))_{0a} = d \left(\frac{n+1}{2} H^2 + \bar{\J} - \J\right) + 2 H \partial_a(H).
$$
Combining these results yields
\begin{align*}
   \hat{\bar{\nabla}}_0(\hat{\overline{\Rho}})_{0a} = \bar{\nabla}_0(\overline{\Rho})_{0a}
   - \partial_a \left(\frac{n+1}{2} H^2 + \bar{\J} - \J\right) - 2 H \partial_a(H) + H \partial_a(H) 
   + 2 H \bar{\Rho}_{0a} + 2 H \partial_a (H).
\end{align*}
By $\bar{\Rho}_0 = - dH$ (see \eqref{CM-trace}), we find
$$
    \hat{\bar{\nabla}}_0(\hat{\overline{\Rho}})_{0a} = \bar{\nabla}_0(\overline{\Rho})_{0a}
   - \partial_a \left(\frac{n+2}{2} H^2 + \bar{\J} - \J\right).
$$
Now the Gauss equation $\bar{\J} - \J = \bar{\Rho}_{00} - \frac{n}{2} H^2$ completes the proof.
\end{proof}

As a consequence of Corollary \ref{Q4-inter} and Lemma \ref{div-ex}, we obtain the formula
\begin{align}\label{Q4-closed}
   {\bf Q}_4 & = Q_4 \notag \\
    & + \frac{2}{n-3} \frac{(n-1)^2}{(n-2)^2} |\W|^2
   - 2\frac{(n-4)(n-1)}{(n-3)(n-2)} (\Rho,\W) - \frac{2(n-1)}{n-3} \left(\delta (\bar{\nabla}_0(\bar{\Rho})_0)
   - \Delta (\bar{\Rho}_{00} + H^2) \right)
\end{align}
if $\lo=0$. The following result further simplifies this formula.

\begin{lem}\label{div-2} Assume that $\lo=0$. Then
$$
   \delta (\bar{\nabla}_0(\bar{\Rho})_0) - \Delta (\bar{\Rho}_{00} + H^2) = - \frac{1}{n-2} \delta \delta (\W).
$$
\end{lem}

\begin{proof} Let the Cotton tensor $\bar{C}$ of $\bar{g}$ be defined by
$
   \bar{C}_{ikj} = \bar{\nabla}_j(\bar{\Rho})_{ik} - \bar{\nabla}_k(\bar{\Rho})_{ij}.
$
Then
$$
   \bar{C}_{ikj} = \frac{1}{n-2} \overline{\Div}_1(\overline{W})_{ikj}
$$
and we calculate
\begin{align*}
   \bar{\nabla}_0(\bar{\Rho})_{0a} & = (\bar{\nabla}_0(\bar{\Rho})_{0a} - \bar{\nabla}_a (\bar{\Rho})_{00})
  + \bar{\nabla}_a (\bar{\Rho})_{00} \\
   & = \bar{C}_{0a0} + \bar{\nabla}_a (\bar{\Rho})_{00} \\
   & = \frac{1}{n-2} \overline{\Div}_1(\overline{W})_{0a0} + \bar{\nabla}_a (\bar{\Rho})_{00}.
\end{align*}
But
\begin{align*}
    & \overline{\Div}_1(\overline{W})_{0a0} = \sum_{i=0}^n \bar{\nabla}^i (\overline{W})_{i0a0}
    = \sum_{i=0}^n \partial^i (\overline{W}_{i0a0}) \\
    & - \overline{W}(\bar{\nabla}^i(\partial_i),\partial_0,\partial_a,\partial_0)
       - \overline{W} (\partial_i,\bar{\nabla}^i(\partial_0),\partial_a,\partial_0)
       - \overline{W}(\partial_i,\partial_0,\bar{\nabla}^i(\partial_a),\partial_0)
       - \overline{W}(\partial_i,\partial_0,\partial_a,\bar{\nabla}^i(\partial_0) \\
    & = \sum_{i=1}^n -\partial^i(\W_{ia})
       + \W(\bar{\nabla}^i(\partial_i),\partial_a)
       + \W(\partial_i,\bar{\nabla}^i(\partial_a))
    - \sum_{i,j=1}^n H h^{ij} \overline{W}_{ija0} - \sum_{i,j=1}^n H h^{ij} \overline{W}_{i0aj}
\end{align*}
using $\bar{\nabla}^i(\partial_0) = L^{ij} \partial_j = H h^{ij} \partial_j$ since $\lo=0$.\footnote{The components for $i=0$ vanish.} Hence we get
$$
     \overline{\Div}_1(\overline{W})_{0a0} = - \delta(\W)_a
$$
since $\overline{W}$ is trace-free. Therefore,
\begin{equation}\label{H1}
    \bar{\nabla}_0(\bar{\Rho})_{0a} = - \frac{1}{n-2} \delta(\W)_a +\bar{\nabla}_a (\bar{\Rho})_{00}.
\end{equation}
Moreover, we find
\begin{align}\label{H2}
    \bar{\nabla}_a (\bar{\Rho})_{00} & = \partial_a (\bar{\Rho}_{00})
    - 2 \bar{\Rho}(\bar{\nabla}_a(\partial_0),\partial_0) \notag \\
    & =  \partial_a (\bar{\Rho}_{00}) - 2 \bar{\Rho} (L_a^b \partial_b,\partial_0) \notag \\
    & = \partial_a (\bar{\Rho}_{00}) - 2 H \bar{\Rho}_{a0} \notag \\
    & = \partial_a (\bar{\Rho}_{00}) + 2 H \partial_a(H) \notag \\
    & = \partial_a (\bar{\Rho}_{00}) + \partial_a (H^2)
\end{align}
using $\lo=0$ and $\bar{\Rho}_0 = - dH$ (by \eqref{CM-trace}). Now we apply $\delta$ to the relations \eqref{H1}
and \eqref{H2} of $1$-forms. Combining the resulting identities proves the assertion.
\end{proof}

Combining \eqref{Q4-closed} with Lemma \ref{div-2} proves the formula for ${\bf Q}_4$ in Corollary \ref{PQ4-gen}.

Finally, we can make the formula for ${\bf P}_4$ in Theorem \ref{P4-inter} fully explicit in the totally umbilic case.
By $\lo=0$, we get $\hat{h}_{(1)} = 2 \lo = 0$,  $\bar{\J} = \J$ (Proposition \ref{J-der}) and
$$
   \hat{h}_{(2)} \stackrel{!}{=} - \hat{\bar{\G}} = -\Rho - \frac{n-1}{n-2} \W
$$
(see Lemma \ref{C-CSC}) so that
\begin{align*}
    {\bf P}_4  & = \Delta^2 - \delta (((n-2) \J h + 4 \Rho) d) + 4 \frac{n-1}{n-2} \delta(\W d)
    + \left(\frac{n}{2}-2\right) {\bf Q}_4.
\end{align*}
This formula completes the proof of Corollary \ref{PQ4-gen}.

\begin{rem} Assume that $\lo=0$. We combine the conformal transformation law for ${\bf P}_4$
with the formula in \eqref{P4-final}. By taking the residue at $n=3$, it follows that the scalar
$
   {\bf R} \st  (\Rho,\W) + \delta \delta (\W)
$
satisfies
$
    e^{4 \varphi} \circ \hat{{\bf R}} = {\bf R}.
$
Lemma \ref{var-W} shows that this relation actually follows from the identity
$$
    \delta(\W d\varphi) - (\delta(\W),d\varphi) - (\Hess(\varphi),\W) = 0.
$$
\end{rem}

The fundamental transformation property of the critical ${\bf Q}_4$ also can be seen as a consequence of
the conformal covariance of ${\bf P}_4$ for general $n$. Therefore, we finish this section with a direct proof of
the conformal covariance of the operator ${\bf P}_4$. In fact, we confirm the transformation law
$$
   e^{(\frac{n}{2}+2) \varphi} \hat{\bf{P}}_4 (f) = {\bf P}_4 (e^{(\frac{n}{2}-2) \varphi} f)
$$
using the explicit formula \eqref{P4-final}. By the known covariance of the intrinsic $P_4$ and the invariance 
of $|\W|^2$, it suffices to prove the conformal covariance of the operator
$$
   f \mapsto 2 (n-3)\delta(\W df)
   + \left(\frac{n}{2}-2\right) \left(- (n-4) (\Rho,\W) + \delta \delta (\W) \right) f.
$$

\begin{lem}\label{var-W} It holds
$$
    e^{(\frac{n}{2}+2) \varphi} \hat{\delta} (\hat{\W} df) - \delta(\W d)(e^{(\frac{n}{2}-2)\varphi} f) =
   - \left(\frac{n}{2}-2\right) \delta (\W d\varphi) e^{(\frac{n}{2}-2)\varphi} f
$$
and
$$
    e^{(\frac{n}{2}+2) \varphi} \hat{\delta} \hat{\delta} (\hat{\W}) - \delta \delta (\W) e^{(\frac{n}{2}-2)\varphi}
    = ((n-2) \delta (\W d\varphi) + (n-4) (\delta (\W), d\varphi)) e^{(\frac{n}{2}-2)\varphi},
$$
up to non-linear terms in $\varphi$.
\end{lem}

\begin{proof} First, we recall that the symmetric bilinear form $\W$ satisfies $\hat{\W} = \W$ and $\tr(\W)=0$.
Now the conformal transformation laws \eqref{conform-div-forms} and \eqref{conform-div-BLF} imply
\begin{align*}
    e^{(\frac{n}{2}+2) \varphi} \hat{\delta} (\hat{\W} df)  & = \delta (e^{(\frac{n}{2}-2) \varphi} \W df) +
   \left(\frac{n}{2}-2\right) (d\varphi, \W df) e^{(\frac{n}{2}-2)\varphi} \\
   & = \delta (\W d(e^{(\frac{n}{2}-2)\varphi} f)) - \left(\frac{n}{2}-2\right) \delta (\W d\varphi e^{(\frac{n}{2}-2)\varphi} f)
   + \left(\frac{n}{2}-2\right) (df, \W d\varphi) e^{(\frac{n}{2}-2)\varphi} \\
   & =  \delta(\W d)(e^{(\frac{n}{2}-2)\varphi} f) - \left(\frac{n}{2}-2\right) \delta (\W d\varphi)  e^{(\frac{n}{2}-2)\varphi} f.
\end{align*}
This proves the first assertion. Next, we calculate
\begin{align*}
    e^{(\frac{n}{2}+2) \varphi} \hat{\delta} \hat{\delta} (\hat{\W}) & = \delta (e^{\frac{n}{2} \varphi} 
    \hat{\delta}(\W))
    + \left(\frac{n}{2}-2\right) (d\varphi, e^{\frac{n}{2} \varphi} \hat{\delta}(\W))  \\
    & = \delta ( \delta ( e^{(\frac{n}{2}-2)\varphi} \W))
    + \frac{n}{2} \delta (\iota_{\grad(\varphi)}(e^{(\frac{n}{2}-2)\varphi} \W))
    + \left(\frac{n}{2}-2\right) (d\varphi, \delta (e^{(\frac{n}{2}-2)\varphi} \W)) \\
    & = \delta ( \delta ( e^{(\frac{n}{2}-2)\varphi} \W))
    + \frac{n}{2} \delta (\W d\varphi) e^{(\frac{n}{2}-2)\varphi}
    + \left(\frac{n}{2}-2\right) (\delta (\W),d\varphi) e^{(\frac{n}{2}-2)\varphi}.
\end{align*}
Now combing this relation with the identity
$$
    \delta \delta (e^{\lambda \varphi} \W) = (\delta \delta (\W) + \lambda (\delta(\W),d\varphi)
   + \lambda \delta (\W d\varphi)) e^{\lambda \varphi}
$$
for $\lambda \in \R$, we find
$$
    e^{(\frac{n}{2}+2) \varphi} \hat{\delta} \hat{\delta} (\hat{\W}) = (\delta (\delta (\W))
   + (n-2) \delta (\W d\varphi) + (n-4) \delta(\W d\varphi)) e^{(\frac{n}{2}-2)\varphi}.
$$
This proves the second relation.
\end{proof}

Therefore, it remains to prove that
$$
   - (n-4) (n-3)\delta (\W d\varphi) + \left(\frac{n}{2}-2\right) \left((n\!-\!4) (\Hess(\varphi),\W)
   +(n\!-\!2) \delta (\W d\varphi) + (n\!-\!4) (\delta (\W), d\varphi)\right) = 0.
$$
This equation is trivial for $n=4$. For $n\ne 4$, it is equivalent to
$$
   -2(n-3) \delta (\W d\varphi) + (n-4)(\Hess(\varphi),\W)
   + (n-2) \delta (\W d\varphi) + (n-4) (\delta (\W), d\varphi) = 0
$$
or in turn to
$$
   - \delta (\W d\varphi) + (\Hess(\varphi),\W) + (\delta (\W), d\varphi) = 0.
$$
However, this is obvious.

\section{The general case of Theorem \ref{main} and the singular Yamabe obstruction $\B_3$}\label{general}

In the present section, we prove Theorem \ref{main} and discuss its consequence for $n=3$.

We first determine the second-order part of ${\bf P}_4$. By \eqref{P4-gen}, it suffices to calculate the sum
$$
   - \delta ( (n-2)\hat{\bar{\J}} d + 4 \hat{h}_{(2)} d) + 4 \delta(\hat{h}_{(1)}^2 d).
$$
In order to express that operator in terms of the given data, we apply the formulas $\hat{h}_{(1)} = 2 \lo$,
\begin{equation}\label{J-bar-hat}
    \hat{\bar{\J}} = \J - \frac{1}{2(n-1)} |\lo|^2
\end{equation}
(by \eqref{J}) and
\begin{align*}
    \hat{h}_{(2)} & = \lo^2 - \hat{\bar{\G}} \\
    & = \lo^2- \left(\Rho + \frac{1}{n-2} \lo^2 -\frac{2n-3}{2(n-1)(n-2)} |\lo|^2 h + \frac{n-1}{n-2} \W\right) \\
    & = \frac{n-3}{n-2} \lo^2 - \Rho + \frac{2n-3}{2(n-1)(n-2)} |\lo|^2 h - \frac{n-1}{n-2} \W
 \end{align*}
(see Lemma \ref{C-CSC}). Now a direct calculation yields
\begin{align*}
   &  (n-2) \hat{\bar{\J}} h + 4  \hat{h}_{(2)} - 4 \hat{h}_{(1)}^2 \\
   & = (n-2) \J h - 4 \Rho - 4 \frac{3n-5}{n-2} \lo^2 - \frac{n^2-12n+16}{2(n-1)(n-2)} |\lo|^2 h - 4 \frac{n-1}{n-2} \W.
\end{align*}
This implies the displayed terms in Theorem \ref{main}. In particular, this proves Corollary \ref{P4-crit}.

\begin{rem}\label{BGW-1-compare} \cite[Corollary 1.1]{BGW-1} states that the second-order part of ${\bf P}_4$ equals
$$
   \delta (4 \Rho - (d-3) \J h) d + \delta \left( 8 \lo^2 + \frac{d^2-4d-1}{2(d-1)(d-2)} |\lo|^2 + 4(d-2) \F \right) d,
$$
where $d=n+1$ and $(d-3) \F = (\lo^2 - \frac{1}{n} |\lo|^2 + \W)$. In terms of the dimension $n$ of $M$, this formula reads
$$
   \delta (4 \Rho - (n-2) \J h) d
   + \delta \left( \left(8 + 4 \frac{n-1}{n-2} \right) \lo^2 \right) d
   + \delta \left( \frac{n^2-12n+16}{2(n-1)(n-2)} |\lo|^2 h + 4  \frac{n-1}{n-2} \W \right) d.
$$
This formula matches with the result in \eqref{P4-gen-dim}.
\end{rem}

We continue with the
\begin{proof}[Proof of the formula for ${\bf Q}_4$ in Theorem \ref{main}]
1. As in the critical dimension $n=4$, the arguments rest on \eqref{Q4-gen}. First, we extend the arguments
in the proof of Corollary \ref{Q4-g-int} in Section \ref{critical} to determine the non-divergence terms.
In fact, formula \eqref{J-double-prime} implies
\begin{align*}
    \hat{\bar{\J}}'' & = \frac{1}{n-3}
    \left(-3 \hat{\bar{\J}}^2 + (\hat{\bar{\G}}, \hat{\overline{\Ric}}) + 2 (\hat{\overline{\Ric}}_{00})^2
    - 2 |\hat{\bar{\G}}|^2 + (\lo^2,\hat{\overline{\Ric}}) + 5 |\lo|^2 \hat{\overline{\Ric}}_{00} - 8 (\lo^2,\hat{\bar{\G}})
    + 24 \sigma_4(\lo) \right) \\
     & + \left(- (\lo,\bar{\nabla}_0(\bar{\Rho})) + \frac{2}{n-3} \lo^{ij} \bar{\nabla}_0(\W)_{ij}
    + \frac{4}{n-3} H (\lo,\W)\right)  \notag \\
    & -  \left((\lo,\Hess(H)) + 2 H(\lo,\bar{\Rho}) - H (\lo,\bar{\G})  - H \tr(\lo^3) - \omega'' |\lo|^2\right),
\end{align*}
up to divergence terms. Simplification of these terms by utilizing the relations indicated in the proof of
Corollary \ref{Q4-g-int} yields the displayed non-divergence terms.

2. It remains to determine the divergence terms in ${\bf Q}_4$. These consist of $-\Delta (\hat{\bar{\J}})$ and the
divergence terms in $-2 \hat{\bar{\J}}''$. Proposition \ref{J-der} shows that the divergence terms in $\hat{\bar{\J}}''$
are given by
\begin{align*}
    & \frac{1}{n-3} \left( - 2 (\lo,\nabla(\overline{\Ric}_0)) - 2 (\delta(\lo),\overline{\Ric}_0)
    + \delta (\bar{\nabla}_0(\overline{\Ric})_0) + \delta ((\lo \overline{\Ric})_0) \right) \\
    & = \frac{1}{n-3} \left( - 2 \delta (\lo \overline{\Ric}_0) + \delta (\bar{\nabla}_0(\overline{\Ric})_0)
    + \delta ((\lo \overline{\Ric})_0) \right) \\
    &  = \frac{1}{n-3} \left(- \delta (\lo \overline{\Ric}_0)  + \delta (\bar{\nabla}_0(\overline{\Ric})_0) \right) \\
    & = - \frac{n-1}{n-3} \delta (\lo \bar{\Rho}_0) + \frac{n-1}{n-3} \delta (\bar{\nabla}_0(\bar{\Rho})_0)
\end{align*}
for the metric $\hat{\bar{g}}$. Now
$$
   \hat{\bar{\Rho}}_{0a} = \bar{\Rho}_{0a} - \overline{\Hess}_{0a}(\omega) = \bar{\Rho}_{0a} + \partial_a(H)
   = \frac{1}{n-1} \delta (\lo)_a
$$
(by \eqref{CM-trace} and the vanishing of $\omega$ on $M$) shows that
$$
    \delta (\lo \hat{\bar{\Rho}}_0) = \frac{1}{n-1} \delta (\lo \delta (\lo)).
$$
Moreover, the calculation of $\delta (\hat{\bar{\nabla}}_0(\hat{\bar{\Rho}})_0)$ rests on
the following generalizations of Lemma \ref{div-ex} and Lemma \ref{div-2}.

\begin{lem}\label{div-ex-g} In general dimensions, it holds
\begin{equation*}
    \delta (\hat{\bar{\nabla}}_0(\hat{\bar{\Rho}})_0) = \delta(\bar{\nabla}_0(\bar{\Rho})_0)
    - \Delta \left(\bar{\Rho}_{00}+ H^2 + \frac{1}{n-1} |\lo|^2\right) - 2 \delta (\lo dH) + \frac{2}{n-1} \delta (H \delta(\lo)).
\end{equation*}
\end{lem}

\begin{proof} An extension of the arguments in the proof of Lemma \ref{div-ex} shows that
$$
    \bar{\nabla}_0(\overline{\Hess}(\omega))_{0a} = \partial_a \left(\frac{n+1}{2} H^2 + \bar{\J} - \J
    + \frac{1}{2(n-1)} |\lo|^2\right) + 2 H \partial_a(H) + 2 (\lo dH)_a.
$$
Hence
\begin{align*}
    \hat{\bar{\nabla}}_0(\hat{\overline{\Rho}})_{0a} & = \bar{\nabla}_0(\overline{\Rho})_{0a}
     - \partial_a \left(\frac{n+1}{2} H^2 + \bar{\J} - \J + \frac{1}{2(n-1)} |\lo|^2\right) - 2 H \partial_a(H)
     - 2 (\lo dH)_a \\
     & + H \partial_a(H) + 2 H \bar{\Rho}_{0a} + 2 H \partial_a(H) \\
     & = \bar{\nabla}_0(\overline{\Rho})_{0a}
     - \partial_a \left(\frac{n+1}{2} H^2 + \bar{\J} - \J + \frac{1}{2(n-1)} |\lo|^2\right)
     - 2 (\lo dH)_a + H \partial_a(H) \\
     & + \frac{2}{n-1} H \delta(\lo)_a - 2 H \partial_a(H)
\end{align*}
by \eqref{CM-trace}. Now the Gauss equation
$$
   \bar{\J} - \J = \bar{\Rho}_{00} + \frac{1}{2(n-1)} |\lo|^2 - \frac{n}{2} H^2
$$
simplifies the latter result to
$$
   \hat{\bar{\nabla}}_0(\hat{\overline{\Rho}})_{0a} = \bar{\nabla}_0(\overline{\Rho})_{0a}
   - \partial_a \left( \bar{\Rho}_{00} + H^2 + \frac{1}{n-1} |\lo|^2 \right)  - 2 (\lo dH)_a
   + \frac{2}{n-1} H \delta(\lo)_a.
$$
Now an application of the operator $\delta$ proves the assertion.
\end{proof}

\begin{lem}\label{div-2-g} In general dimensions, it holds
$$
   \delta (\bar{\nabla}_0(\bar{\Rho})_0) - \Delta (\bar{\Rho}_{00} + H^2)
   = - \frac{1}{n-2} \delta \delta (\W) - \frac{1}{n-2} \delta (\lo^{ij} \overline{W}_{0ij\cdot})
   - \frac{2}{n-1} \delta (L \delta(\lo)) + 2 \delta (\lo dH).
$$
\end{lem}

\begin{proof} We extend the arguments in the proof of Lemma \ref{div-2}. The relation
$$
   \overline{\Div}_1(\overline{W})_{0a0} = - \delta(\W)_a - \lo^{ij} \overline{W}_{ija0} - \lo^{ij} \overline{W}_{i0aj}
   = - \delta(\W)_a - \lo^{ij} \overline{W}_{0ija}
$$
implies
\begin{align*}
   \bar{\nabla}_0(\bar{\Rho})_{0a} & = \frac{1}{n-2} \overline{\Div}_1(\overline{W})_{0a0}
   + \bar{\nabla}_a (\bar{\Rho})_{00} \\
   & = - \frac{1}{n-2} \delta(\W)_a - \frac{1}{n-2} \lo^{ij} \overline{W}_{0ija} + \bar{\nabla}_a(\bar{\Rho})_{00}.
\end{align*}
But
\begin{align*}
    \bar{\nabla}_a (\bar{\Rho})_{00} & = \partial_a (\bar{\Rho}_{00})
    - 2 \bar{\Rho}(\bar{\nabla}_a(\partial_0),\partial_0) \notag \\
    & =  \partial_a (\bar{\Rho}_{00}) - 2 L_a^b \bar{\Rho}_{b0} \notag \\
    & = \partial_a (\bar{\Rho}_{00}) + 2 L_a^b \partial_b(H) - \frac{2}{n-1} (L \delta(\lo))_a \\
    & = \partial_a (\bar{\Rho}_{00}) + 2 H \partial_a(H) + 2 \lo_a^b \partial_b(H) - \frac{2}{n-1} (L \delta (\lo))_a
\end{align*}
by $\bar{\Rho}_0 = - dH + \frac{1}{n-1} \delta(\lo)$ (see \eqref{CM-trace}). Therefore,
\begin{align*}
    \delta (\bar{\nabla}_0(\bar{\Rho})_0) = - \frac{1}{n-2} \delta \delta (\W) - \frac{1}{n-2} \delta (\lo^{ij}
     \overline{W}_{0ij\cdot}) + \Delta(\bar{\Rho}_{00} + H^2) + 2 \delta (\lo dH) - \frac{2}{n-1} \delta(L \delta(\lo)).
\end{align*}
The proof is complete.
\end{proof}

Using these results, we simplify the divergence terms in $\hat{\bar{\J}}''$ as
\begin{align}\label{div-inter}
     & \frac{n-1}{n-3} \Big( -\frac{1}{n-1} \delta (\lo \delta(\lo)) + \delta(\bar{\nabla}_0(\bar{\Rho})_0)
    - \Delta \left(\bar{\Rho}_{00} + H^2\right) -  \frac{1}{n-1} \Delta (|\lo|^2) \notag \\
     & - 2 \delta (\lo dH) + \frac{2}{n-1} \delta (H \delta(\lo)) \Big) \notag \\
     & = \frac{n-1}{n-3} \Big( -\frac{1}{n-1} \delta (\lo \delta(\lo)) - \frac{1}{n-2} \delta \delta (\W)
     - \frac{1}{n-2} \delta (\lo^{ij} \overline{W}_{0ij\cdot}) - \frac{2}{n-1} \delta (L \delta(\lo)) + 2 \delta (\lo dH) \notag \\
     & -  \frac{1}{n-1} \Delta (|\lo|^2) - 2 \delta (\lo dH) + \frac{2}{n-1} \delta (H \delta(\lo)) \Big) \notag \\
     & =  \frac{n-1}{n-3} \left( - \frac{1}{n-2} \delta \delta (\W) - \frac{3}{n-1} \delta (\lo \delta(\lo))
      - \frac{1}{n-2} \delta (\lo^{ij} \overline{W}_{0ij\cdot}) - \frac{1}{n-1} \Delta (|\lo|^2) \right).
\end{align}

Next, the following identity enables us to replace the contribution by the Weyl tensor in the latter sum by contributions
in terms of $\lo$.

\begin{lem}\label{van} In general dimensions, it holds
$$
    \delta (\lo^{ij} \overline{W}_{0ij\cdot}) =  \delta \delta (\lo^2) - \frac{1}{2} \Delta (|\lo|^2)
    - \frac{n-2}{n-1} \delta(\lo \delta(\lo)).
$$
\end{lem}

\begin{proof} First, the trace-free Codazzi-Mainardi equation (see \eqref{CM-TF-3}) implies
\begin{align}\label{A}
    \lo^{ij} \overline{W}_{aij0} & = \lo^{ij} \left( \nabla_i (\lo)_{aj} - \nabla_a(\lo)_{ij}
    + \frac{1}{n-1} \delta(\lo)_i h_{aj} - \frac{1}{n-1} \delta(\lo)_a h_{ij} \right) \notag \\
    & = \lo^{ij} \nabla_i(\lo)_{aj} - \lo^{ij} \nabla_a(\lo)_{ij} + \frac{1}{n-1} \lo_a^i \delta(\lo)_i.
\end{align}
Second, the relations $d(|\lo|^2)_a = 2 \lo^{ij} \nabla_a(\lo)_{ij}$ and
$$
    (\lo \delta(\lo))_a = \lo^j_a \nabla^i(\lo)_{ij} \quad \mbox{and} \quad
    \delta(\lo^2)_a = \nabla^i(\lo)_i^j \lo_{ja} + \lo^j_i \nabla^i(\lo)_{ja}
$$
imply
\begin{align}\label{B}
    \delta (\lo^2)_a - \frac{1}{2} d (|\lo|^2)_a - \frac{1}{2} (\lo \delta(\lo))_a
    & = \nabla^i(\lo)_i^j \lo_{ja} + \lo^j_i \nabla^i(\lo)_{ja}
    - \lo^{ij} \nabla_a(\lo)_{ij} - \frac{1}{2} \lo^j_a \nabla^i (\lo)_{ij} \notag \\
    & = \frac{1}{2} \delta(\lo)^j \lo_{ja} + \lo^j_i \nabla^i(\lo)_{ja} - \lo^{ij} \nabla_a(\lo)_{ij}.
\end{align}
Taking the difference of \eqref{A} and \eqref{B} shows the identity
$$
   \lo^{ij} \overline{W}_{0ija} - \delta (\lo^2)_a + \frac{1}{2} d (|\lo|^2)_a + \frac{1}{2} (\lo \delta(\lo))_a
   = \frac{1}{n-1} (\lo \delta(\lo))_a - \frac{1}{2} (\lo \delta(\lo))_a
$$
of $1$-forms. Equivalently, we find
\begin{equation}\label{WN}
   \lo^{ij} \overline{W}_{0ij\cdot} =  \delta (\lo^2) - \frac{1}{2} d (|\lo|^2) - \frac{n-2}{n-1} \lo \delta(\lo).
\end{equation}
Now an application of the operator $\delta$ to the last relation proves the assertion.
\end{proof}

Combining \eqref{div-inter} with Lemma \ref{van} yields

\begin{lem}\label{div-Q} The divergence terms of ${\bf Q}_4$ are given by the sum of
\begin{align*}
     & \frac{2(n-1)}{(n-3)(n-2)} \delta \delta (\W) + \frac{2(n-1)}{(n-3)(n-2)} \delta \delta (\lo^2)
    + \frac{4}{n-3} \delta (\lo \delta(\lo))\\
    & + \left(\frac{2}{n-3} - \frac{n-1}{(n-2)(n-3)}\right) \Delta (|\lo|^2)
\end{align*}
and
$$
    - \Delta (\J) + \frac{1}{2(n-1)} \Delta (|\lo|^2).
$$
\end{lem}

Now Lemma \ref{div-Q} implies \eqref{div-terms}. This completes the proof of Theorem \ref{main}.
\end{proof}

Finally, we discuss the residue formula \cite{JO1}
\begin{equation}\label{res-BQ}
    24 \B_3 = \res_{n=3} ({\bf Q}_4).
\end{equation}
We shall use this identity to reproduce the known formula for $\B_3$. Theorem \ref{main} shows that
$
   \res_{n=3} ({\bf Q}_4)
$
equals the sum of the divergence terms
$$
    4 \delta \delta (\W) + 4 \delta \delta (\lo^2)  + 4 \delta (\lo \delta(\lo))
$$
(see also Lemma \ref{div-Q}), the terms
$$
    - 4 \lo^{ij} \bar{\nabla}_0(\overline{W})_{0ij0} - 8 H (\lo,\W)
$$
and
$$
     8 |\W|^2 + 4 (\Rho,\W) + 32 (\lo^2,\W) + 12 (\lo^2,\Rho) - 4 \J |\lo|^2 - 8 |\lo|^4  + 24 \tr(\lo^4).
$$
Note also that $\sigma_4(\lo)=0$ for $n=3$ and that this relation is equivalent to  $2 \tr(\lo^4) = |\lo|^4$.
Hence it holds $- 8 |\lo|^4  + 24 \tr(\lo^4) = 4 |\lo|^4$.

Now we compare this result with the formula for $\B_3$ in \cite[Theorem 1]{JO2}.\footnote{Alternative formulas for
$\B_3$ were given in \cite[Proposition 1.1]{GGHW} and its arXiv-version. We refer to \cite[Section 6.5]{JO2} 
for discussing the relations between these formulas.} This formula implies
\begin{align}\label{B3-complete}
    24 \B_3 & =  8 |\W|^2 + 4 (\Rho,\W) + 32 (\lo^2,\W) + 12 (\lo^2,\Rho) - 4 \J |\lo|^2  + 4 |\lo|^4  - 8 H(\lo,\W)
    - 4  \lo^{ij} \bar{\nabla}_0(\overline{W})_{0ij0} \notag \\
    & - 8 \lo^{ij} \nabla^k \overline{W}_{kij0} + 4 |\overline{W}_0|^2 +
   12 \delta \delta (\lo^2) - 4 \Delta (|\lo|^2) + 4 \delta \delta (\W).
\end{align}

Hence
$$
    \res_{n=3} ({\bf Q}_4) - 24 \B_3
$$
equals
\begin{equation}\label{BQ-diff}
  - 8 \delta \delta (\lo^2) + 4 \delta (\lo \delta(\lo)) + 4 \Delta (|\lo|^2) + 8 \lo^{ij} \nabla^k \overline{W}_{kij0}
  - 4 |\overline{W}_0|^2.
\end{equation}

In order to simplify this sum, we apply the following result.

\begin{lem}\label{div-term} In any dimension $n \ge 2$, it holds
$$
    |\overline{W}_0|^2 - 2 \lo^{ij} \nabla^k \overline{W}_{kij0} = - 2 \delta (\lo^{ij} \overline{W}_{\cdot ij0}).
$$
\end {lem}

\begin{proof} First, we note that
$$
   \delta  (\lo^{ij} \overline{W}_{\cdot ij0}) = \nabla^k (\lo)^{ij} \overline{W}_{kij0} + \lo^{ij} \nabla^k \overline{W}_{kij0}.
$$
Hence
$$
    |\overline{W}_0|^2 - 2 \lo^{ij} \nabla^k \overline{W}_{kij0} =  |\overline{W}_0|^2 + 2 \nabla^k (\lo)^{ij} \overline{W}_{kij0}
    - 2  \delta (\lo^{ij} \overline{W}_{\cdot ij0}).
$$
Thus, it remains to prove that
\begin{equation}\label{W-vanish}
    |\overline{W}_0|^2 + 2 \nabla^k (\lo)^{ij} \overline{W}_{kij0}  = 0.
\end{equation}
The trace-free Codazzi-Mainardi equation \eqref{CM-TF-3}
implies that $|\overline{W}_0|^2 = \overline{W}_{ijk0} \overline{W}^{ijk0}$ equals the sum of
\begin{align*}
    & \frac{1}{(n-1)^2} \left(\delta(\lo)_j \delta(\lo)^j h_{ik}h^{ik}
    + \delta(\lo)_i \delta(\lo)^i h_{jk}h^{jk} - \delta(\lo)_j \delta(\lo)^i h_{ik} h^{jk}
    - \delta(\lo)_i \delta(\lo)^j h_{jk} h^{ik}\right) \\
    & =  \frac{1}{(n-1)^2} (2n (\delta(\lo),\delta(\lo)) - 2 (\delta(\lo),\delta(\lo))) = \frac{2}{n-1}(\delta(\lo),\delta(\lo)),
\end{align*}
\begin{align*}
     & \frac{2}{n-1} \left(\nabla_j(\lo)_{ik} \delta(\lo)^j h^{ik} - \nabla_j (\lo)_{ik} \delta(\lo)^i h^{jk}
     - \nabla_i(\lo)_{jk} \delta(\lo)^j h^{ik} + \nabla_i(\lo)_{jk} \delta(\lo)^i h^{jk} \right)\\
     & = - \frac{2}{n-1} \left(\nabla^k(\lo)_{ik} \delta(\lo)^i + \nabla^k (\lo)_{jk} \delta(\lo)^j\right)
     = - \frac{4}{n-1} (\delta(\lo),\delta(\lo))
\end{align*}
and
\begin{align*}
    & \nabla_j(\lo)_{ik} \nabla^i(\lo)^{ik} - \nabla_j(\lo)_{ik} \nabla^i(\lo)^{jk} - \nabla_i(\lo)_{jk} \nabla^j(\lo)^{ik}
    + \nabla_i(\lo)_{jk} \nabla^i(\lo)^{jk} \\
    & = 2 \nabla_j(\lo)_{ik} \nabla^j(\lo)^{ik} - 2 \nabla_j(\lo)_{ik} \nabla^i(\lo)^{jk}.
\end{align*}
Hence
\begin{equation}\label{W-1}
     |\overline{W}_0|^2 = - \frac{2}{n-1} |\delta(\lo)|^2 +  2 \nabla_j(\lo)_{ik} \nabla^j(\lo)^{ik} - 2 \nabla_j(\lo)_{ik} \nabla^i(\lo)^{jk}.
\end{equation}
On the other hand, \eqref{CM-TF-3} shows that
\begin{align}\label{W-2}
    \nabla^k (\lo)^{ij} \overline{W}_{kij0} & =  \nabla^k(\lo)^{ij} \left (\nabla_i(\lo)_{kj} - \nabla_k(\lo)_{ij}
    + \frac{1}{n-1} \delta(\lo)_i h_{kj} - \frac{1}{n-1} \delta(\lo)_k h_{ij} \right) \notag\\
    & = \nabla^k(\lo)^{ij} \nabla_i (\lo)_{kj} - \nabla^k(\lo)^{ij} \nabla_k(\lo)_{ij}
    + \frac{1}{n-1} \delta(\lo)_i \delta(\lo)^i.
\end{align}
Combining \eqref{W-1} and \eqref{W-2} proves \eqref{W-vanish}. The proof is complete.
\end{proof}

Lemma \ref{div-term} shows that the sum \eqref{BQ-diff} equals
\begin{equation}\label{BQ-diff-2}
    8 \delta (\lo^{ij} \overline{W}_{0ij\cdot}) - 8 \delta \delta (\lo^2) + 4 \delta (\lo \delta(\lo)) + 4 \Delta (|\lo|^2).
\end{equation}
But Lemma \ref{van} for $n=3$ shows that the sum \eqref{BQ-diff-2} vanishes.

\begin{cor}\label{QB} $\res_{n=3}({\bf Q}_4) = 24 \B_3$.
\end{cor}

\section{The invariants $\Jv_1$, $\Jv_2$ and $\mathcal{C}$. Proofs of Theorem \ref{LCI} and Theorem \ref{alex}}\label{C}


The first main result of the present section is

\begin{prop}\label{J-12} In general dimensions, the quantities
\begin{align}\label{J1-g}
   \Jv_1 & \st  \lo^{ij} \bar{\nabla}_0(\overline{W})_{0ij0} + 2 H (\lo,\W) + \tfrac{n-2}{(n-1)^2} |\delta(\lo)|^2 \notag \\
    & - \tfrac{n-2}{n-3} (\lo^2,\Rho) - \tfrac{n-2}{(n-3)(n-6)} \J |\lo|^2 + \tfrac{n-4}{(n-3)(n-6)} \Delta (|\lo|^2) -
    \tfrac{1}{n-3} \delta \delta (\lo^2)
\end{align}
and
\begin{align}\label{J2-g}
    \Jv_2 & \st (\lo,\bar{\nabla}_0(\bar{\Rho})) + (\lo,\Hess(H)) + H(\lo,\Rho) - \tfrac{n-3}{n-2} H (\lo,\W) \notag \\
    & + \tfrac{n}{(n-1)^2} |\delta(\lo)|^2 - \bar{\Rho}_{00} |\lo|^2 - \tfrac{1}{(n-3)(n-6)} \J |\lo|^2
    - \tfrac{3}{2} H^2 |\lo|^2 - \tfrac{n-3}{n-2} H \tr (\lo^3) + \tfrac{n-4}{n-3} (\lo^2,\Rho) \notag \\
    & - \tfrac{1}{n-3} \delta \delta (\lo^2) + \tfrac{n-5}{2(n-3)(n-6)} \Delta (|\lo|^2)
\end{align}
are local conformal invariants of weight $-4$ of an embedding $M^n \hookrightarrow X^{n+1}$, i.e., it holds
$$
   e^{4 \iota^*(\varphi)} \hat{\Jv}_i = \Jv_i
$$ 
for $i=1,2$ and all $\varphi \in C^\infty(X)$.
\end{prop}

Some comments are in order.

Both invariants $\Jv_i$ vanish if $\lo=0$ and have a simple formal pole at $n=3$. The formal residue of
$\Jv_1$ at $n=3$ equals
$$
    - \delta \delta (\lo^2) + \frac{1}{3} \Delta (|\lo|^2) - (\lo^2,\Rho) + \frac{1}{3} \J |\lo|^2
    = - (\delta \delta ((\lo^2)_\circ) + (\Rho, (\lo^2)_\circ)) = - \mathcal{D}((\lo^2)_\circ),
$$
where the operator $\mathcal{D}(b) = \delta \delta (b) + (\Rho,b)$ acts on trace-free symmetric bilinear forms
and $(\lo^2)_\circ$ denotes the trace-free part of $\lo^2$. It is well-known that
$\mathcal{D}: b \mapsto \delta \delta (b) + (\Rho,b)$ is a conformally covariant operator $S_0^2(M) \to C^\infty(M)$
on trace-free symmetric bilinear forms on $M^3$. We recall that the term $\mathcal{D}((\lo^2)_\circ)$ contributes to the
singular Yamabe obstruction $\B_3$ of $M^3 \hookrightarrow X^4$ (see \eqref{B3-complete}). The formal residue of
$\Jv_2$ at $n=3$ also equals $-\mathcal{D}((\lo^2)_\circ)$.

We also note that both invariants $\Jv_i$ have a simple formal pole at $n=6$ with residues being proportional to the local
invariant $P_2 (|\lo|^2)$ of weight $-4$.

We shall see that, in dimension $n=4$, a linear combination of $\Jv_1$ and $\Jv_2$ equals $\mathcal{C}$
(defined in \eqref{new invariant}). Thus, Proposition \ref{J-12} proves Theorem \ref{LCI}. Theorem \ref{alex}
then is an easy consequence.

In order to prove Proposition \ref{J-12}, it suffices to prove the vanishing of the respective conformal variations
$$
    (\Jv_i(g))^\bullet[\varphi] = (d/dt)|_0 (e^{4t \iota^*(\varphi)} \Jv_i(e^{2t\varphi} g))
$$
of $\Jv_i$ at the metric $g$. We shall use the bullet notation also for the conformal variation of other scalar curvature quantities.

\begin{lem}\label{V-1} It holds
\begin{align}\label{CV-P-g}
    & (\lo, \Hess(H) + \bar{\nabla}_0 (\bar{\Rho}))^\bullet[\varphi] \notag \\
    & = \frac{1}{2} |\lo|^2 \Delta(\varphi) - \frac{2n}{n\!-\!1} (\lo \delta(\lo),d\varphi) + H (\lo,\Hess (\varphi)) \notag \\
    & - |\lo|^2  \partial^2_{0} (\varphi) - (\lo,\Rho) \partial_0(\varphi) + \frac{n\!-\!3}{n\!-\!2} \tr (\lo^3) \partial_0(\varphi)
    + 3 H |\lo|^2 \partial_0(\varphi) + \frac{n\!-\!3}{n\!-\!2} (\lo,\W) \partial_0(\varphi) \notag \\
    & - \frac{1}{2} \delta (|\lo|^2 d\varphi) + 2 \delta (\lo^2 d\varphi)
\end{align}
and
\begin{align}\label{CV-W-g}
   & (\lo^{ij} \bar{\nabla}_0 (\overline{W})_{0ij0})^\bullet[\varphi] \notag \\
   & = -2 (\lo,\W) \partial_0(\varphi) + 2 \delta (\lo^2 d\varphi)- 2 (\lo^2,\Hess(\varphi))
   - 2 \frac{n-2}{n-1} (\lo \delta(\lo),d\varphi) \notag \\
   & +  |\lo|^2 \Delta(\varphi) - \delta(|\lo|^2 d\varphi).
\end{align}
\end{lem}

\begin{proof} In the following calculations, all non-linear terms in $\varphi$ will be omitted without mentioning.
We first calculate the conformal variation of
$$
    (\lo,\Hess(H)) + (\lo,\bar{\nabla}_0(\bar{\Rho})).
$$
We recall that
$$
   e^\varphi \hat{H} = H + \partial_0(\varphi) \quad \mbox{and} \quad \widehat{\Hess}_{ij}(u)
   = \Hess_{ij}(u) - u_i \varphi_j - u_j \varphi_i + h_{ij} (d\varphi,du).
$$
Hence
\begin{align*}
   \widehat{\Hess}_{ij} (\hat{H}) & = \Hess_{ij}(\hat{H}) - (\hat{H}_i \varphi_j + \hat{H}_j \varphi_i)
   + h_{ij}(d\varphi,d\hat{H}) \\
   & = \Hess_{ij} (e^{-\varphi} H + e^{-\varphi} \partial_0(\varphi)) -  e^{-\varphi} (H_i \varphi_j + H_j \varphi_i)
   + h_{ij} e^{-\varphi} (d\varphi, dH + \partial_0 (d\varphi)) \\
   & = e^{-\varphi} \Hess_{ij}(H) -  e^{-\varphi}  (H_i \varphi_j + H_j \varphi_i) +  H \Hess_{ij}(e^{-\varphi}) \\
   & +  e^{-\varphi} \Hess_{ij}(\partial_0(\varphi)) -  e^{-\varphi} (H_i \varphi_j + H_j \varphi_i)
   + h_{ij} e^{-\varphi} (d\varphi,dH)  \\
   & = e^{-\varphi} \Hess_{ij}(H) -  e^{-\varphi}  2 (H_i \varphi_j + H_j \varphi_i) +  e^{-\varphi} \Hess_{ij}(\partial_0(\varphi))
   - e^{-\varphi} H \Hess_{ij}(\varphi) \\
   & +  h_{ij} e^{-\varphi} (d\varphi,dH).
\end{align*}
Therefore, we get
\begin{equation*}
    e^{4 \varphi} \hat{\lo}^{ij}\widehat{\Hess}_{ij} (\hat{H}) = \lo^{ij} (\Hess_{ij}(H) - 2 (H_i \varphi_j + H_j \varphi_i)
   + \Hess_{ij}(\partial_0(\varphi)) - H \Hess_{ij}(\varphi)).
\end{equation*}
Hence
\begin{align}\label{CT-1}
   (\lo,\Hess(H))^\bullet[\varphi] = - 4 (\lo dH,d\varphi) + (\lo,\Hess (\partial_0(\varphi))) - H (\lo,\Hess(\varphi)).
\end{align}
Next, we calculate
\begin{align*}
   \widehat{\bar{\nabla}_0(\bar{\Rho})_{ij}} & = e^{-\varphi} \hat{\bar{\nabla}}_0(\bar{\Rho}
   - \overline{\Hess}(\varphi))_{ij} \\
   & =  e^{-\varphi} (\bar{\nabla}_0(\bar{\Rho} - \overline{\Hess}(\varphi))_{ij}
   - (\bar{\Rho}_{0j} \varphi_i + \bar{\Rho}_{i0} \varphi_j + 2 \bar{\Rho}_{ij} \partial_0(\varphi))) \\
   & =  e^{-\varphi} \bar{\nabla}_0(\bar{\Rho})_{ij} - e^{-\varphi} \bar{\nabla}_0 (\overline{\Hess}(\varphi))_{ij}
    -  e^{-\varphi}  (\bar{\Rho}_{0j} \varphi_i + \bar{\Rho}_{i0} \varphi_j)
   - 2 e^{-\varphi} \bar{\Rho}_{ij}  \partial_0(\varphi)
\end{align*}
using the general transformation law  
\begin{equation}\label{CTL-nabla}
   \hat{\nabla}_i(\partial_j) = \nabla_i(\partial_j) + \partial_i(\varphi) \partial_j
   + \partial_j(\varphi) \partial_i - g_{ij} \grad(\varphi).
\end{equation}
Now we contract with $\lo$. Then
\begin{equation}\label{var-1}
   (\lo,\bar{\nabla}_0 (\bar{\Rho}))^\bullet[\varphi]
   = - (\lo, \bar{\nabla}_0 (\overline{\Hess}(\varphi)))
    - 2 (\lo \bar{\Rho}_0,d\varphi) - 2 (\lo,\bar{\Rho}) \partial_0 (\varphi).
\end{equation}
We continue discussing the term $\bar{\nabla}_0(\overline{\Hess}(\varphi))$. We find
\begin{align*}
    \bar{\nabla}_0 (\overline{\Hess} (\varphi))_{ij} & = \partial_0 (\overline{\Hess}_{ij}(\varphi))
    - \overline{\Hess}(\varphi)(\bar{\nabla}_0 (\partial_i), \partial_j)
    - \overline{\Hess}(\varphi)(\partial_i,\bar{\nabla}_0 (\partial_j)) \\
    & = \partial_0 (\overline{\Hess}_{ij} (\varphi)) - L_i^l \Hess_{lj}(\varphi) - L_j^l \Hess_{il}(\varphi)
    - 2 L^2_{ij} \partial_0(\varphi)
\end{align*}
using
$$
    \overline{\Hess}_{ij}(u) = \Hess_{ij}(u) + L_{ij} \partial_0(u) \quad \mbox{and}
   \quad \bar{\nabla}_0(\partial_i) = L_i^k \partial_k.
$$
Here we use $\bar{\Gamma}_{0i}^j = L_i^j$ and $\bar{\Gamma}_{0i}^0 = 0$. Hence
\begin{align}\label{n-Hess}
    \bar{\nabla}_0(\overline{\Hess}(\varphi))_{ij}
    & = \partial_0(\overline{\Hess}_{ij}(\varphi)) - \lo_i^l \Hess_{lj}(\varphi)
   - \lo_j^l \Hess_{il}(\varphi) - 2 H \Hess_{ij}(\varphi) - 2 L^2_{ij} \partial_0(\varphi).
\end{align}
But
\begin{align*}
   \partial_0(\overline{\Hess}_{ij}(\varphi)) & = \partial_0 (\partial_{ij}^2(\varphi) - \bar{\Gamma}_{ij}^k \varphi_k) \\
   & = \partial_{ij}^2 (\partial_0(\varphi)) - \bar{\Gamma}_{ij}^k \partial^2_{0k}(\varphi)
   - \partial_0 (\bar{\Gamma}_{ij}^k) \varphi_k \\
   & = \partial_{ij}^2 (\partial_0(\varphi)) - \bar{\Gamma}_{ij}^l \partial^2_{0l}(\varphi)
   - \bar{\Gamma}_{ij}^0 \partial^2_{0} (\varphi)
   - \partial_0 (\bar{\Gamma}_{ij}^l) \varphi_l - \partial_0 (\bar{\Gamma}_{ij}^0) \partial_0(\varphi) \\
    & = \Hess_{ij} (\partial_0(\varphi)) - \partial_0 (\bar{\Gamma}_{ij}^l)  \varphi_l
   - \bar{\Gamma}_{ij}^0 \partial^2_{0} (\varphi)  - \partial_0 (\bar{\Gamma}_{ij}^0) \partial_0(\varphi).
\end{align*}
Here $k$ runs from $0$ to $n$ and the tangential index $l$ runs from $1$ to $n$ ($n=4$). In the last equality, 
we used the fact that for tangential indices, the restriction of the Christoffel symbols $\bar{\Gamma}_{ij}^l$ 
to $M$ coincide with the Christoffel symbols
$\Gamma_{ij}^l$ of the induced metric $h$ on $M$. Now the general variation formula
\begin{align*}
   \delta (\Gamma_{ij}^k) = \frac{1}{2} \left( \nabla_i(\delta(g))_j^k + \nabla_j(\delta(g))_i^k - \nabla^k (\delta(g))_{ij}\right)
\end{align*}
(with $\nabla$ for $g$) implies
\begin{align*}
   \partial_0(\bar{\Gamma}_{ij}^l)  & = \nabla_i(L)_j^l + \nabla_j(L)_i^l - \nabla^l (L)_{ij} \\
   & = \nabla_i(\lo)_j^l + \nabla_j(\lo)_i^l - \nabla^l(\lo)_{ij} + H_j h_i^l + H_i h_j^l - H^l h_{ij}
\end{align*}
using $h_r = h + 2L r + \cdots$. Moreover, the identity $\bar{\Gamma}_{ij}^0 = - \frac{1}{2} h_{ij}'$ and the expansion
$h_r = h +2L r +(L^2 - \bar{\G}) r^2 + \cdots$ imply
\begin{align}
     \bar{\Gamma}_{ij}^0 & = - L_{ij}, \label{Gamma-1} \\
     \partial_0 (\bar{\Gamma}_{ij}^0) & = - L^2_{ij} + \bar{\G}_{ij}. \label{Gamma-2}
\end{align}
Hence
\begin{align*}
   \partial_0(\overline{\Hess}_{ij}(\varphi)) & =  \Hess_{ij} (\partial_0(\varphi)) \\
   & - (\nabla_i(\lo)_j^l + \nabla_j(\lo)_i^l - \nabla^l (\lo)_{ij})\varphi_l - H_j \varphi_i - H_i \varphi_j + H^l \varphi_l h_{ij} \\
   & + L_{ij} \partial^2_{0} (\varphi)  + L^2_{ij} \partial_0(\varphi) - \bar{\G}_{ij} \partial_0(\varphi).
\end{align*}
Thus \eqref{n-Hess} implies
\begin{align*}
    \bar{\nabla}_0(\overline{\Hess}(\varphi))_{ij} & =  \Hess_{ij} (\partial_0(\varphi)) \\
    & - (\nabla_i(\lo)_j^l + \nabla_j(\lo)_i^l - \nabla^l (\lo)_{ij})\varphi_l - H_j \varphi_i - H_i \varphi_j + H^l \varphi_l h_{ij} \\
    & - ((\lo^2)_{ij} + 2 H \lo_{ij} + H^2 h_{ij}) \partial_0(\varphi)
   + L_{ij} \partial^2_{0} (\varphi) - \bar{\G}_{ij} \partial_0(\varphi) \\
    & - \lo_i^l \Hess_{lj}(\varphi) - \lo_j^l \Hess_{il}(\varphi) - 2 H \Hess_{ij}(\varphi).
\end{align*}
By contraction with $\lo$, we obtain
\begin{align*}
    (\lo, \bar{\nabla}_0 (\overline{\Hess}(\varphi))) & =  (\lo, \Hess (\partial_0(\varphi))) \\
    & - 2 \lo^{ij} \nabla_i (\lo)_j^l \varphi_l + \lo^{ij} \nabla^l (\lo)_{ij} \varphi_l - 2 \lo^{ij} H_i \varphi_j \\
    & - \tr (\lo^3) \partial_0 (\varphi) - 2 H |\lo|^2 \partial_0 (\varphi) + |\lo|^2  \partial^2_{0} (\varphi)
    - (\lo,\bar{\G}) \partial_0(\varphi) \\
    & - 2 (\lo^2, \Hess (\varphi)) - 2 H (\lo,\Hess(\varphi)).
\end{align*}
Reordering gives
\begin{align*}\label{pi-1}
   (\lo, \bar{\nabla}_0 (\overline{\Hess}(\varphi))) & =  (\lo, \Hess (\partial_0(\varphi))) \notag \\
   & + \lo^{ij} \nabla^l(\lo)_{ij} \varphi_l - 2 \lo^{ij} \nabla_i (\lo)_j^l \varphi_l - 2 \lo^{ij} H_i \varphi_j - 2 (\lo^2, \Hess (\varphi)) - 2 H (\lo,\Hess(\varphi)) \notag \\
   & + |\lo|^2  \partial^2_{0} (\varphi) - \tr (\lo^3) \partial_0(\varphi) - 2 H |\lo|^2 \partial_0(\varphi) - (\lo,\bar{\G}) \partial_0(\varphi).
\end{align*}
Now we apply the identities
\begin{equation}\label{pi-h1}
    \lo^{ij} \nabla^l(\lo)_{ij} \varphi_l = \frac{1}{2} (d (|\lo|^2),d\varphi)
    = \frac{1}{2} \delta(|\lo|^2 d\varphi) - \frac{1}{2} |\lo|^2 \Delta(\varphi)
\end{equation}
and
\begin{equation}\label{pi-h2}
    \lo^{ij} \nabla_i(\lo)_j^l \varphi_l = \delta (\lo^2 d\varphi) - (\lo^2,\Hess(\varphi)) - (\lo \delta(\lo),d\varphi).
\end{equation}
These relations show that
\begin{align*}
   (\lo,\bar{\nabla}_0(\overline{\Hess}(\varphi))) & = (\lo, \Hess (\partial_0(\varphi))) \\
   & - \frac{1}{2} |\lo|^2 \Delta (\varphi) + 2 (\lo \delta(\lo), d\varphi) - 2 H (\lo,\Hess(\varphi)) - 2 (\lo dH,d\varphi) \\
   & + |\lo|^2  \partial^2_0 (\varphi) - \tr (\lo^3) \partial_0(\varphi) - 2 H |\lo|^2 \partial_0(\varphi) - (\lo,\bar{\G}) \partial_0(\varphi) \\
   & + \frac{1}{2} \delta (|\lo|^2 d\varphi) - 2 \delta (\lo^2 d\varphi).
\end{align*}
Combining this result with \eqref{var-1} yields
\begin{align}\label{CV-NRho}
   (\lo,\bar{\nabla}_0 (\bar{\Rho}))^\bullet[\varphi]  & = - (\lo, \Hess (\partial_0 (\varphi))) \notag \\
   & + \frac{1}{2} |\lo|^2 \Delta (\varphi) - 2 (\lo \delta(\lo), d\varphi)
   - 2 (\lo \bar{\Rho}_0,d \varphi) - 2 (\lo,\bar{\Rho}) \partial_0(\varphi)  \notag \\
   & + 2 H (\lo,\Hess(\varphi)) + 2 (\lo dH,d\varphi)  \notag \\
   & - |\lo|^2  \partial^2_{0} (\varphi)  + \tr (\lo^3) \partial_0(\varphi) + 2 H |\lo|^2 \partial_0(\varphi)
  + (\lo,\bar{\G}) \partial_0(\varphi) \notag \\
  & - \frac{1}{2} \delta (|\lo|^2 d\varphi) + 2 \delta (\lo^2 d\varphi).
\end{align}
Now summarizing the conformal variations \eqref{CT-1} and \eqref{CV-NRho} gives
\begin{align*}
   & (\lo, \Hess (H) + \bar{\nabla}_0 (\bar{\Rho}))^\bullet[\varphi] \\
   & =  \frac{1}{2} |\lo|^2 \Delta(\varphi) - 2  (\lo \delta(\lo),d\varphi) - 2  (\lo dH,d\varphi) - 2 (\lo \bar{\Rho}_0,d\varphi)
   + H (\lo,\Hess (\varphi)) \notag \\
   & - |\lo|^2  \partial^2_{0} (\varphi) - 2 (\lo,\bar{\Rho}) \partial_0(\varphi) + (\lo,\bar{\G}) \partial_0(\varphi)
   + \tr (\lo^3)  \partial_0(\varphi) + 2 H |\lo|^2 \partial_0(\varphi) \\
   & - \frac{1}{2} \delta (|\lo|^2 d\varphi) + 2 \delta (\lo^2 d\varphi).
\end{align*}
By $\delta(\lo) = (n\!-\!1) dH + (n\!-\!1) \bar{\Rho}_0$ (Codazzi-Mainardi), we have
$$
    - 2 (\lo \delta(\lo),d\varphi) - 2 (\lo dH,d\varphi) - 2 (\lo \bar{\Rho}_0,d\varphi)
    = \left(-2-\frac{2}{n\!-\!1}\right)  (\lo \delta(\lo),d\varphi) = -\frac{2n}{n\!-\!1} (\lo \delta(\lo),d\varphi).
$$
Therefore, we conclude that
\begin{align}\label{Var-IRho}
   & (\lo, \Hess(H) + \bar{\nabla}_0 (\bar{\Rho}))^\bullet[\varphi] \notag \\
   & = \frac{1}{2} |\lo|^2 \Delta(\varphi)  - \frac{2n}{n-1} (\lo \delta(\lo),d\varphi)
   + H (\lo,\Hess (\varphi)) \notag \\
   & - |\lo|^2  \partial^2_{0} (\varphi) - 2 (\lo,\bar{\Rho}) \partial_0(\varphi) + (\lo,\bar{\G}) \partial_0(\varphi)
   + \tr (\lo^3) \partial_0(\varphi) + 2 H |\lo|^2 \partial_0(\varphi) \notag \\
   & - \frac{1}{2} \delta (|\lo|^2 d\varphi) + 2 \delta (\lo^2 d\varphi).
\end{align}
Now the decomposition $\bar{\G} = \bar{\Rho} + \bar{\Rho}_{00} h + \W$ implies $(\lo,\bar{\G}) = (\lo,\bar{\Rho}) + (\lo,\W)$.
Hence
\begin{align*}
   & (\lo, \Hess(H) + \bar{\nabla}_0 (\bar{\Rho}))^\bullet[\varphi] \\
   & = \frac{1}{2} |\lo|^2 \Delta(\varphi) - \frac{2n}{n\!-\!1} (\lo \delta(\lo),d\varphi) + H (\lo,\Hess (\varphi)) \\
   & - |\lo|^2  \partial^2_{0} (\varphi) - (\lo,\bar{\Rho}) \partial_0(\varphi)
   + (\lo,\W)  \partial_0(\varphi) + \tr (\lo^3) \partial_0(\varphi) + 2 H |\lo|^2 \partial_0(\varphi) \\
   & - \frac{1}{2} \delta (|\lo|^2 d\varphi) + 2 \delta (\lo^2 d\varphi)
\end{align*}
Next, we use the Fialkow equation\eqref{Fial} to write
\begin{align*}
   (\lo,\bar{\Rho}) & = \left(\lo,\Rho - H \lo + \frac{1}{n\!-\!2} \lo^2 + \frac{1}{n\!-\!2} \W \right) \\
   & = (\lo,\Rho) - H |\lo|^2 + \frac{1}{n\!-\!2} \tr(\lo^3)  + \frac{1}{n\!-\!2} (\lo,\W).
\end{align*}
Therefore, we find
\begin{align*}
    & (\lo, \Hess(H) + \bar{\nabla}_0 (\bar{\Rho}))^\bullet[\varphi] \\
    & = \frac{1}{2} |\lo|^2 \Delta(\varphi) - \frac{2n}{n\!-\!1} (\lo \delta(\lo),d\varphi) + H (\lo,\Hess (\varphi)) \\
    & - |\lo|^2  \partial^2_{0} (\varphi) - (\lo,\Rho) \partial_0(\varphi) + \frac{n\!-\!3}{2} \tr (\lo^3) \partial_0 (\varphi)
    + 3 H |\lo|^2 \partial_0(\varphi) + \frac{n\!-\!3}{n\!-\!2} (\lo,\W) \partial_0(\varphi) \\
    & - \frac{1}{2} \delta (|\lo|^2 d\varphi) + 2 \delta (\lo^2 d\varphi).
\end{align*}
This proves the first variation formula. For the proof of the second variation formula, we first observe that
\begin{align}\label{CTL-DW}
    \widehat{ (\lo^{ij} \bar{\nabla}_0(\overline{W})_{0ij0})}
    & = e^{-3\varphi} \hat{\lo}^{ij} \hat{\bar{\nabla}}_0 (\hat{\overline{W}})_{0ij0} \notag \\
    & = e^{-6 \varphi} \lo^{ij}  \hat{\bar{\nabla}}_0 (e^{2\varphi} \overline{W})_{0ij0} \notag \\
    & = e^{-6 \varphi} \lo^{ij} (e^{2\varphi} \hat{\bar{\nabla}}_0 (\overline{W})_{0ij0}
   + 2 \partial_0(\varphi) e^{2\varphi} \overline{W}_{0ij0}) \notag \\
    & = e^{-4\varphi} \lo^{ij} \hat{\bar{\nabla}}_0(\overline{W})_{0ij0} 
   + 2 e^{-4 \varphi} \partial_0(\varphi) (\lo,\W)  \notag \\
    & = e^{-4\varphi} (\lo^{ij} \hat{\bar{\nabla}}_0(\overline{W})_{0ij0} + 2  \partial_0(\varphi) (\lo,\W)).
\end{align}
Now the general transformation law \eqref{CTL-nabla} implies
\begin{align*}
   \hat{\bar{\nabla}}_0(\overline{W})_{0ij0} & = \bar{\nabla}_0(\overline{W})_{0ij0}
   - 2 \overline{W}(\partial_0(\varphi) \partial_0,\partial_i,\partial_j,\partial_0)
   + \overline{W}(\grad(\varphi),\partial_i,\partial_j,\partial_0) \\
   & - \overline{W}(\partial_0, \partial_0(\varphi) \partial_i,\partial_j,\partial_0)
   - \overline{W}(\partial_0, \partial_i, \partial_0(\varphi) \partial_j,\partial_0) \\
    & - 2\overline{W}(\partial_0,\partial_i,\partial_j,\partial_0(\varphi) \partial_0)
   + \overline{W}(\partial_0,\partial_i,\partial_j,\grad(\varphi)) \\
   & = \bar{\nabla}_0(\overline{W})_{0ij0}
    - 6 \partial_0(\varphi) \W_{ij} + \overline{W}_{\grad(\varphi)ij0} + \overline{W}_{0ij\grad(\varphi)}.
\end{align*}
Hence the right-hand side of \eqref{CTL-DW} equals
$$
    e^{-4\varphi} (\lo^{ij} \bar{\nabla}_0(\overline{W})_{0ij0} - 4 \partial_0(\varphi) (\lo,\W)
   + 2 \lo^{ij} \overline{W}_{\grad(\varphi) ij0}).
$$
Therefore, we get
\begin{equation}\label{CI-W-term}
   (\lo^{ij} \bar{\nabla}_0 (\overline{W})_{0ij0})^\bullet[\varphi] =
   -2(\lo,\W) \partial_0(\varphi) + 2 \lo^{ij} \overline{W}_{\grad^t(\varphi) ij0}.
\end{equation}
Now we further simplify the term
$$
    \lo^{ij} \overline{W}_{\grad^t(\varphi) ij0}
$$
using the trace-free Codazzi-Mainardi equation
$$
   \nabla_j(\lo)_{ik} - \nabla_i(\lo)_{jk} + \frac{1}{n\!-\!1} \delta(\lo)_j h_{ik} - \frac{1}{n\!-\!1} \delta(\lo)_{i} h_{jk}
   = \overline{W}_{ijk0}
$$
(see \eqref{CM-TF-3}). It follows that
\begin{align*}
   \lo^{ij} \overline{W}_{\grad^t(\varphi)ij0} & = \lo^{ij} \nabla_i(\lo)_{\grad^t(\varphi)j}
   - \lo^{ij} \nabla_{\grad^t(\varphi)} (\lo)_{ij} + \frac{1}{n\!-\!1} \lo^{ij} \delta(\lo)_i h_{\grad^t (\varphi) j} \\
   & = \lo^{ij} \nabla_i(\lo)_{\grad^t(\varphi)j} - \frac{1}{2} (d (|\lo|^2),d\varphi)
  + \frac{1}{n\!-\!1} \lo^{ij} \delta(\lo)_i h_{\grad^t (\varphi) j}.
\end{align*}
Now \eqref{pi-h2} shows that
\begin{equation*}
    \lo^{ij} \nabla_i(\lo)_j^l \varphi_l = \delta (\lo^2 d\varphi) - (\lo^2,\Hess(\varphi)) - (\lo \delta(\lo),d\varphi).
\end{equation*}
Hence
\begin{align*}
   & \lo^{ij} \overline{W}_{\grad^t(\varphi)ij0} \notag \\
   & = \delta (\lo^2 d\varphi) - (\lo^2,\Hess(\varphi)) - (\lo \delta(\lo), d\varphi)
    + \frac{1}{2} |\lo|^2 \Delta (\varphi) - \frac{1}{2} \delta(|\lo|^2 d\varphi) + \frac{1}{n\!-\!1} (\lo \delta(\lo),d\varphi) \notag \\
   & = \delta (\lo^2 d\varphi) - (\lo^2,\Hess(\varphi)) - \frac{n\!-\!2}{n\!-\!1} (\lo \delta(\lo),d\varphi)
   + \frac{1}{2} |\lo|^2 \Delta (\varphi) - \frac{1}{2} \delta(|\lo|^2 d\varphi)
\end{align*}
(see also \eqref{WN}). Thus, we obtain
\begin{align*}
   & (\lo^{ij} \bar{\nabla}_0 (\overline{W})_{0ij0})^\bullet[\varphi] \notag \\
   & = -2 (\lo,\W) \partial_0(\varphi) + 2 \delta (\lo^2 d\varphi)- 2 (\lo^2,\Hess(\varphi))
   - 2\frac{n\!-\!2}{n\!-\!1} (\lo \delta(\lo),d\varphi) +  |\lo|^2 \Delta(\varphi) - \delta(|\lo|^2 d\varphi).
\end{align*}
This proves the second formula.
\end{proof}

\begin{rem}\label{invariant-E} The arguments in the second part of the above proof show that
$$
    e^\varphi\widehat{\bar{\nabla}_0 (\overline{W})_{0ij0}}
    = \bar{\nabla}_0(\overline{W})_{0ij0} - 4 \partial_0(\varphi) \W_{ij} 
    + \overline{W}_{\grad(\varphi)ij0} + \overline{W}_{0ij\grad(\varphi)}.
$$
Thus, the conformal transformation law
$
   e^\varphi \widehat{\bar{C}_{ij0}} = \hat{\bar{C}}_{ij0} = \bar{C}_{ij0} + \overline{W}_{\grad(\varphi)ij0}
$
for the Cotton tensor $\bar{C}$ implies the conformal invariance $e^\varphi \hat{S}_{ij} =  S_{ij}$
of the trace-free symmetric bilinear form
\begin{equation}
   S_{ij} \st  \bar{\nabla}_0(\overline{W})_{0ij0} - \bar{C}_{ij0} - \bar{C}_{ji0} + 4 H \W_{ij}
\end{equation}
on $M$. As a consequence, the scalar curvature quantity $\Jv_5 \st (\lo,S)$ is a conformal invariant of 
weight $-4$. For more details on $\Jv_5$, we refer to Section \ref{ECI}.

For $n=3$, the conformally invariant symmetric tensor $S$ recently appeared in 
\cite[Lemma 2.1]{CG-1} in connection with the study of the metric variation of the conformally invariant 
functional
\begin{equation}\label{Weyl-gravity}
   \int_{X^4} |W|^2 dvol_g + 8 \int_{M^3} (\lo,\W) dvol_h
\end{equation}
on a four-manifold $X$ with boundary $M$ (for more details, we also refer to \cite{GZ1}). This functional 
generalizes the conformally invariant functional
$$
  \int_{X^4} |W|^2 dvol_g 
$$
of closed four manifolds $X$. Critical metrics of the latter functional are Bach-flat. The boundary term in \eqref{Weyl-gravity} 
may be regarded as an analog of the Gibbons-Hawking-York term leading to a well-defined variational problem for 
the Einstein-Hilbert functional on a manifold with boundary.\footnote{The fact that the variational problem of the 
functional \eqref{Weyl-gravity} is well-defined also suggests to expect that the same combination of a bulk and a 
boundary term contributes to the integrated conformal anomaly of CFT's on a four-manifold with boundary \cite[(14)]{Sol}.} 
It was noted in \cite{GZ1} that critical points of the functional \eqref{Weyl-gravity} are Bach-flat and satisfy the equation 
$S=0$.

For $n=3$, the conformal invariant $\Jv_5 = (\lo,S)$ of weight $-4$ also contributes to the singular Yamabe obstruction 
$\B_3$ of $M^3 \hookrightarrow X^4$.  Indeed, we calculate
\begin{align*}
   (\lo,S) & = - \lo^{ij} \bar{\nabla}_0(\overline{W})_{0ij0} - 2 \lo^{ij} \bar{\nabla}^k (\overline{W})_{kij0} 
   + 4 H (\lo,\W) \\
   & = - \lo^{ij} \bar{\nabla}_0(\overline{W})_{0ij0} - 2 ( \lo^{ij} \nabla^k \overline{W}_{kij0} - (\lo^2,\W) 
   + 3H (\lo,\W) - \lo^{ij} \lo^{kl} \overline{W}_{kijl})  + 4 H (\lo,\W) \\
   & = - \lo^{ij} \bar{\nabla}_0(\overline{W})_{0ij0} - 2\lo^{ij} \bar{\nabla}^k \overline{W}_{kij0} - 2 H (\lo,\W) 
   + 2(\lo^2,\W) + 2 \lo^{ij} \lo^{kl} \overline{W}_{kijl}
\end{align*}
(for the second equality see \eqref{n-reduce} in the proof of Proposition \ref{B-JJD}). Comparing this formula with 
\begin{align*}
   12 \B_3 & = 6 \LOP ((\lo^2)_\circ) + 2 \LOP (\W)  \notag \\
   & - 2 \lo^{ij} \bar{\nabla}^0(\overline{W})_{0ij0} - 4 \lo^{ij} \nabla^k \overline{W}_{kij0} - 4 H(\lo,\W)
  + 2 |\lo|^4 + 16 (\lo^2,\W) + 4 |\W|^2 + 2 |\overline{W}_{0}|^2
\end{align*}
(see \cite[Theorem 1]{JO2}) yields
\begin{equation}\label{B3-S}
   12 \B_3 = 2 \Jv_5 + 6 \LOP ((\lo^2)_\circ) + 2 |\lo|^4  + 2 \LOP (\W) + 12 (\lo^2,\W) + 4 |\W|^2 
   + 2 |\overline{W}_{0}|^2.
\end{equation}
\end{rem}

\begin{rem}\label{NL} In the proof of Corollary \ref{Q4-g-int} in Section \ref{critical}, we need to know
how $(\lo,\bar{\nabla}_0(\bar{\Rho}))$ transforms under the conformal change from $\bar{g}$ to $\hat{\bar{g}}
= e^{2\omega} \bar{g}$. In addition to the terms which are linear in $\omega$, this also requires determining
the non-linear contributions by $\omega$. By the conformal transformation law
$$
    \hat{\bar{\Rho}} = \bar{\Rho} - \overline{\Hess}(\omega) + d\omega \otimes d \omega
   - \frac{1}{2} |d\omega|^2 \bar{g},
$$
all non-linear contributions by $\omega$ are caused by $(\lo,\hat{\bar{\nabla}}_0 (\overline{\Hess}(\omega)))$
and
$$
    (\lo,\hat{\bar{\nabla}}_0(d\omega \otimes \omega)) - \frac{1}{2} (\lo,\hat{\bar{\nabla}}_0 (|d\omega|^2 \bar{g})).
$$
But $\omega = 0$ on $M$ implies $\hat{\bar{\nabla}}_0(d\omega \otimes \omega)_{ij} = 0$, and one easily sees
that $\hat{\bar{\nabla}}_0 (|d\omega|^2 \bar{g})_{ij}$ is a multiple of $h_{ij}$. The latter term vanishes by
contraction with $\lo$. It remains to determine the non-linear contributions which are caused by $\hat{\bar{\nabla}}_0 (\overline{\Hess}(\omega))_{ij}$.
But
\begin{align*}
   \hat{\bar{\nabla}}_0(\overline{\Hess}(\omega))_{ij} & = \bar{\nabla}_0(\overline{\Hess}(\omega))_{ij}
    - 2 \overline{\Hess}(\omega)_{ij} \partial_0(\omega) \\
    & = \bar{\nabla}_0(\overline{\Hess}(\omega))_{ij} - 2 L_{ij} (\partial_0(\omega))^2 \\
    & = \bar{\nabla}_0(\overline{\Hess}(\omega))_{ij} - 2 L_{ij} H^2
\end{align*}
using $\omega = 0$, $\partial_0(\omega) = - H$ and $\overline{\Hess}_{ij}(\omega) = \Hess_{ij}(\omega) + L_{ij} \partial_0(\omega) = L_{ij} \partial_0(\omega)$ on $M$. Therefore, the additional contribution
to $(\lo,\bar{\nabla}_0(\bar{\Rho}))$ is the term $2 H^2 |\lo|^2$.
\end{rem}

With the above preparations, we are able to give a

\begin{proof}[Proof of Proposition \ref{J-12}]
For the proofs of the conformal invariance of $\Jv_i$, we combine \eqref{CV-P-g} and \eqref{CV-W-g} with the
variation formulas
$$
    (H (\lo,\W))^\bullet[\varphi] = (\lo,\W) \partial_0(\varphi),
    \quad  (\lo^2,\Rho)^\bullet[\varphi] =  -(\lo^2,\Hess(\varphi)), \quad
    (\J |\lo|^2)^\bullet[\varphi] = - |\lo|^2 \Delta(\varphi),
$$
\eqref{var-2}, \eqref{var-L} and
\begin{equation}\label{CV-square}
    (|\delta(\lo)|^2)^\bullet[\varphi] = 2(n-1) (\lo \delta(\lo),d\varphi),
\end{equation}
The latter relation follows from
$$
   e^{4 \varphi} \widehat{(\delta(\lo),\delta(\lo))}
   = e^{\varphi} \widehat{\delta(\lo)}_i e^\varphi \widehat{\delta(\lo)}_j h^{ij}
   = (\delta (\lo) + (n-1) \lo d\varphi)_i (\delta(\lo) + (n-1) \lo d\varphi)_j h^{ij}
$$
using \eqref{conform-div-BLF} (for $\lambda = -1$). These variation formulas easily imply that the
conformal variations of $\Jv_1$ and $\Jv_2$ vanish. We omit the details.
\end{proof}

As an application of the local invariants $\Jv_1$, $\Jv_2$ in general dimensions, we obtain the following
decomposition of ${\bf Q}_4$.

\begin{theorem}\label{Q-deco-I} In general dimensions, the extrinsic $Q$-curvature ${\bf Q}_4$ admits the decomposition
\begin{align*}
    {\bf Q}_4 & = Q_4 - \tfrac{15n^4-49n^3+36n^2+24n-32}{8(n-3)(n-2)^2(n-1)^2} \Iv_1
   + \tfrac{2(5n^2-14n+9)}{(n-3)(n-2)^2}  \Iv_2 + \tfrac{2 (n-1)^2}{(n-3)(n-2)^2} \Iv_4
   + \tfrac{4 (3n-5)(n-1)}{(n-3)(n-2)^2} \Iv_6  \\
   & - 4 \Jv_1 + 2 \Jv_2  + (n-4) {\bf E}_4 + \mbox{total divergence},
\end{align*}
where
\begin{align*}
    {\bf E}_4 & \st \tfrac{4}{n-3} \lo^{ij} \bar{\nabla}_0(\overline{W})_{0ij0} - \tfrac{2(n-1)}{(n-3)(n-2)} (\Rho,\W)
   + \tfrac{8}{n-3} H (\lo,\W) \\
    & - \tfrac{2(4n-7)}{(n-3)(n-2)} (\lo^2,\Rho) - \tfrac{n^3+n^2+8n-20}{2(n-1)(n-2)(n-3)(n-6)} \J |\lo|^2
    + \tfrac{2}{(n-1)^2} |\delta(\lo)|^2).
\end{align*}
\end{theorem}

\begin{rem}\label{E-var} In the critical dimension $n=4$, we have
$$
    {\bf E}_4  = 4 \lo^{ij} \bar{\nabla}_0(\overline{W})_{0ij0} - 3 (\Rho,\W)  + 8 H (\lo,\W) - 9  (\lo^2,\Rho)
    - \frac{5}{2}  \J |\lo|^2 + \frac{2}{9} |\delta(\lo)|^2.
$$
A calculation shows that
$$
    \left(\int_M {\bf E}_4 dvol_h \right)^\bullet[\varphi] = \int_M \varphi
    \left[ \delta \delta (\lo^2) + \frac{1}{6} \Delta (|\lo|^2) + 4 \delta (\lo \delta(\lo)) + 3 \delta \delta (\W) \right] dvol_h.
$$
The integrand on the right-hand side is given by the divergence part of ${\bf Q}_4$ (see \eqref{Q4-ex1}). This result
also follows using general principles. Let $n > 4$. Combining the variation formula
$$
    \left(\int_M {\bf Q}_4 dvol_h \right)^\bullet[\varphi] = (n-4) \int_M \varphi {\bf Q}_4 dvol_h
$$
(and similarly for $Q_4$) with the decomposition ${\bf Q}_4 = Q_4 + \Iv + (n-4) {\bf E}_4 + \delta$ with a local conformal
invariant $\Iv$ of weight $-4$ and a total divergence $\delta$ yields
$$
   (n-4) \int_M \varphi {\bf Q}_4 dvol_h = (n-4) \int_M \varphi Q_4 dvol_h + (n-4) \int_M \varphi \Iv dvol_h
  + (n-4)  \left(\int_M {\bf E}_4 dvol_h \right)^\bullet[\varphi].
$$
We divide by $n-4$ and conclude that
$$
    (n-4) \int_M \varphi {\bf E}_4 dvol_h + \int_M \varphi \delta dvol_h
   =  \left(\int_M {\bf E}_4 dvol_h \right)^\bullet[\varphi].
$$
Now continuation in dimension implies
$$
    \left(\int_M {\bf E}_4 dvol_h \right)^\bullet[\varphi] =  \int_M \varphi \delta dvol_h
$$
for $n=4$.
\end{rem}

The following result is the special case of Proposition \ref{J-12} in the critical dimension $n=4$.

\begin{lem}\label{new-invariants} In dimension $n=4$, the quantities
\begin{equation}\label{EI-1}
   \Jv_1 \st \lo^{ij} \bar{\nabla}_0(\overline{W})_{0ij0} + 2 H (\lo,\W) + \frac{2}{9} |\delta(\lo)|^2
   - 2 (\lo^2,\Rho) + \J |\lo|^2 - \delta \delta (\lo^2)
\end{equation}
and
\begin{align}\label{EI-2}
   \Jv_2 \st & (\lo, \bar{\nabla}_0(\bar{\Rho})) + H (\lo,\Rho) - \frac{1}{2} H (\lo,\W) + (\lo,\Hess(H)) \notag \\
   & + \frac{4}{9} |\delta(\lo)|^2
   - \bar{\Rho}_{00} |\lo|^2 + \frac{1}{2} \J |\lo|^2 - \frac{3}{2} H^2 |\lo|^2 - \frac{1}{2} H \tr(\lo^3)
   - \delta \delta (\lo^2) + \frac{1}{4} \Delta (|\lo|^2)
\end{align}
are local conformal invariants  of weight $-4$ of the embedding $M^4 \hookrightarrow X^5$, i.e., it holds
$e^{4 \iota^*(\varphi)} \hat{\Jv}_i =\Jv_i$ for all $\varphi \in C^\infty(X)$, $i=1,2$.
\end{lem}

The integrated invariant $\Jv_1$ in the critical dimension $n=4$ was discovered in \cite{AS} (for more details, we refer to 
Section \ref{deco-general}).

We continue to consider the invariants $\Jv_1i$ and $\Jv_2$  for $n=4$. Then the relation $\mathcal{C} = - 4 \Jv_1 + 2 \Jv_2$ 
immediately implies

\begin{cor}
$\mathcal{C}^\bullet[\varphi] = 0$.
\end{cor}

This also completes the proof of Theorem \ref{LCI}.

We continue with the

\begin{proof}[Proof of Theorem \ref{alex}]
Corollary \ref{Q4-g-int} and Lemma \ref{div-Q} show  that, in the critical dimension $n=4$, the extrinsic
$Q$-curvature ${\bf Q}_4$ is given by the sum  of
\begin{equation}\label{Q-GB-LI}
    2 \J^2 - 2|\Rho|^2 + \mathcal{C} + \frac{9}{2} |\W|^2 + 21 (\lo^2,\W) + \frac{33}{2} \tr (\lo^4) 
   - \frac{14}{3} |\lo|^4
\end{equation}
and the divergence terms
\begin{align*}
      - \Delta (\J) + 2 \Delta (|\lo|^2) + \frac{1}{6} \Delta (|\lo|^2) + 6 \delta (\lo \delta(\lo))
     + 3 \delta \delta (\W) + 3  \delta (\lo^{ij} \overline{W}_{0ij\cdot}).
\end{align*}
By Theorem \ref{LCI},  the individual terms in \eqref{Q-GB-LI}, except $2\J^2-2|\Rho|^2$ are
local conformal invariants.  Since $\Pf_4 = \J^2 - |\Rho|^2 + \frac{1}{8} |W|^2$, this proves 
Theorem \ref{alex}.
\end{proof}

We finish this section with the formulation of a conjectural decomposition of the critical extrinsic $Q$-curvature in
higher dimensions. Its role in more general contexts will be discussed in Section \ref{deco-general}.

\begin{conj}\label{extrinsic-DS} Let $n$ be even. Then the critical extrinsic $Q$-curvature ${\bf Q}_n(g)$ is a linear
combination of $Q_n(h)$, local conformal invariants of the embedding $M \hookrightarrow X$ and a total divergence.
By the Deser-Schwimmer decomposition of $Q_n(h)$ (see \cite{alex}), this is equivalent to the existence of a
decomposition of ${\bf Q}_n(g)$ as a linear combination of the Pfaffian of $(M,h)$, local conformal invariants of the
embedding and a total divergence.

Let $n$ be odd. Then ${\bf Q}_n(g)$ is a linear combination of local conformal invariants of the embedding
$M \hookrightarrow X$ and a total divergence.
\end{conj}

\section{The Graham-Reichert functional}\label{GR-functional}

In \cite{GR}, Graham and Reichert studied the asymptotic expansion of the volume of
minimal hypersurfaces $M$ (of arbitrary codimension) in a Poincar\'e-Einstein background $X$.
The coefficient of $\log \varepsilon$ ($\varepsilon$ being a cut-off parameter) in these expansions is
a global conformal invariant. We shall refer to it as the Graham-Reichert functional.

In the codimension-one special case, the following result describes the structure of the Graham-Reichert functional
$\mathcal{E}_{GR}$ of $M^4 \hookrightarrow X^5$ from the perspective of an analog of Conjecture 
\ref{extrinsic-DS} (see also Corollary \ref{EGR-DS}).

\begin{lem}\label{GR-I} It holds
\begin{align}\label{E-GR2}
   8 \mathcal{E}_{GR} & = \int_M (\J^2-|\Rho|^2) dvol_h \notag \\
   & + \int_M \left(-(\lo^2,\Rho) - (\Rho,\W) + \frac{1}{2} \J |\lo|^2 + \frac{1}{9} |\delta(\lo)|^2
   + \bar{B}_{00}\right) dvol_h \notag \\
  & + \int_M \left(\frac{1}{12} |\lo|^4 - \frac{1}{4} \tr(\lo^4)\right) dvol_h
   - \int_M \left(\frac{1}{4} |\W|^2 + \frac{1}{2} (\lo^2,\W)\right) dvol_h
\end{align}
or, equivalently,
\begin{align}\label{E-GR}
   8 \mathcal{E}_{GR} & = \int_M (\J^2-|\Rho|^2) dvol_h \notag \\
   & + \frac{1}{2} \int_M \Jv_1 dvol_h
   - \int_M \left(\frac{1}{2} \lo^{ij} \bar{\nabla}_0(\overline{W})_{0ij0} + H(\lo,\W) + (\Rho,\W)
    - \bar{B}_{00} \right) dvol_h \notag \\
    & + \int_M \left(\frac{1}{12} |\lo|^4 - \frac{1}{4} \tr(\lo^4)\right) dvol_h
   - \int_M \left(\frac{1}{4} |\W|^2 + \frac{1}{2} (\lo^2,\W)\right) dvol_h.
\end{align}
\end{lem}

We recall that $\J^2- |\Rho|^2 = \Pf_4- \frac{1}{8} |W|^2$. 

All integrals in Lemma \ref{GR-I} are conformally invariant (see Remark \ref{CI-Bach}).

\begin{proof} In the codimension-one case, \cite[Proposition 5.1]{GR} is equivalent to
\begin{align}\label{GR-form}
    8 \mathcal{E}_{GR} & = \int_M \left(|dH|^2 - H^2 |\lo|^2 + 3H^4\right) dvol_h \notag \\
    & + \int_M \left(2 H h^{ij} \bar{\nabla}_0(\bar{\Rho})_{ij} + 4 (\bar{\Rho}_0, dH) + 5 H^2 h^{ij} \bar{\Rho}_{ij}
    - 8 \bar{\Rho}_{00} H^2\right) dvol_h \notag \\
   & + \int_M \left (-\bar{\Rho}^{ij} \bar{\Rho}_{ij} + (\bar{\Rho}_0,\bar{\Rho}_0) + (h^{ij} \bar{\Rho}_{ij})^2
   - h^{ij} \bar{B}_{ij}\right) dvol_h,
\end{align}
where $\bar{B}$ is the Bach tensor of the background metric. Let $\bar{G} = \overline{\Ric} - 4 \bar{\J} \bar{g}$
be the Einstein tensor of $\bar{g}$ on $X^5$. Then
\begin{align*}
    \bar{\nabla}_0(\bar{\Rho})_{00} & = \frac{1}{3} \bar{\nabla}_0(\overline{\Ric} - \bar{\J} \bar{g})_{00} \\
    & = \frac{1}{3} \bar{\nabla}_0 (\bar{G})_{00} + \bar{\nabla}_0(\bar{\J} \bar{g})_{00} \\
    & = \frac{1}{3} (-\delta (\overline{\Ric}_0) - 4 H \overline{\Ric}_{00} + (L,\overline{\Ric})) + \bar{\J}'
\end{align*}
using Lemma \ref{Nabla-1G}. Hence
\begin{align*}
   h^{ij} \bar{\nabla}_0(\bar{\Rho})_{ij} & = \bar{\nabla}_0(\bar{\J}) - \bar{\nabla}_0(\bar{\Rho})_{00} \\
   & = \frac{1}{3} (\delta (\overline{\Ric}_0) + 4 H \overline{\Ric}_{00} - (L,\overline{\Ric})) \\
   & = \delta (\bar{\Rho}_0) + 4 H \overline{\Rho}_{00} -  (L,\bar{\Rho}).
\end{align*}
This yields the relation
$$
   \int_M H h^{ij} \bar{\nabla}_0(\bar{\Rho})_{ij} dvol_h = \int_M (4 H^2 \bar{\Rho}_{00} - H (L, \bar{\Rho})
   + H \delta(\bar{\Rho}_0)) dvol_h.
$$
Now, abbreviating the first integral in \eqref{GR-form} by $(\cdot)_{\text{Guven}}$, we find
\begin{align*}
   8 \mathcal{E}_{GR} =  (\cdot)_{\text{Guven}}
   & + \int_M (-2H(\lo,\bar{\Rho}) + 3 H^2 h^{ij} \bar{\Rho}_{ij} + 2 (dH,\bar{\Rho}_0)) dvol_h \\
   & + \int_M (-\bar{\Rho}^{ij} \bar{\Rho}_{ij} + (\bar{\Rho}_0,\bar{\Rho}_0)
   + (h^{ij} \bar{\Rho}_{ij})^2 - h^{ij} \bar{B}_{ij}) dvol_h
\end{align*}
using partial integration. Next, we substitute the Fialkow equation
$$
   \bar{\Rho} = \Rho - H \lo - \frac{1}{2} H^2 h + \frac{1}{2} \left(\lo^2 - \frac{|\lo|^2}{6} h + \W \right)
$$
(see \eqref{Fial}) and $h^{ij} \bar{\Rho}_{ij} = \bar{\J} - \bar{\Rho}_{00} = 
\J + \frac{1}{6} |\lo|^2 - 2 H^2$ (by the Gauss equation) into that formula and simplify the result. This gives
\begin{align*}
    8 \mathcal{E}_{GR} & = \int_M (\J^2 - |\Rho|^2) dvol_h
   + \int \left(|dH|^2 + 2(dH,\bar{\Rho}_0) + |\bar{\Rho}_0|^2\right) dvol_h \\
   & + \int_M \left( \frac{1}{2} \J |\lo|^2 - (\lo^2,\Rho) + \frac{1}{12} |\lo|^4 - \frac{1}{4} \tr(\lo^4) \right) dvol_h\\
   & - \int_M \left(( \Rho,\W) + \frac{1}{2} (\lo^2,\W) + \frac{1}{4} |\W|^2 + h^{ij} \bar{B}_{ij}\right) dvol_h.
\end{align*}
But the relation $3 \bar{\Rho}_0 + 3 dH = \delta(\lo)$ (Codazzi-Mainardi) shows that
$$
    |dH|^2 + 2(dH,\bar{\Rho}_0) + |\bar{\Rho}_0|^2 = \frac{1}{9} |\delta(\lo)|^2.
$$
Thus, we finally arrive at
\begin{align*}
    8 \mathcal{E}_{GR} &  =  \int_M (\J^2 - |\Rho|^2) dvol_h
   + \int_M \left(\frac{1}{9} |\delta(\lo)|^2 + \frac{1}{2} \J |\lo|^2
   - (\lo^2,\Rho) + \frac{1}{12} |\lo|^4 - \frac{1}{4} \tr(\lo^4) \right) dvol_h\\
    & - \int_M \left((\Rho,\W) + \frac{1}{2} (\lo^2,\W) + \frac{1}{4} |\W|^2 + h^{ij} \bar{B}_{ij}\right) dvol_h.
\end{align*}
This implies the assertion using $h^{ij} \bar{B}_{ij} = - \bar{B}_{00}$ since $\bar{B}$ is trace-free.
\end{proof}

\begin{rem}\label{flat-check} As a cross-check of the coefficients in the first integral in the last line of \eqref{E-GR},
one may verify the relation
\begin{align*}
    & \int_M (|dH|^2 - H^2 |\lo|^2 + 3 H^4) dvol_h \\
    & = \int_M (\J^2 - |\Rho|^2) dvol_h
    + \int_M \frac{1}{2} \Jv_1 dvol_h + \int_M \left(\frac{1}{12} |\lo|^4 - \frac{1}{4} \tr(\lo^4)\right) dvol_h
\end{align*}
for the flat background $\R^5$ using \eqref{J-Rho-flat} and $\delta(\lo) = 3dH$.
\end{rem}

\begin{rem}\label{CI-Bach} The first integral in \eqref{E-GR} is a global conformal invariant.
$\Jv_1$ and the integrands in the third line of \eqref{E-GR} are local conformal invariants. Thus, the
conformal invariance of the functional $\mathcal{E}_{GR}$ implies that of the integral
\begin{equation}\label{inv-2}
    \int_M \left(\frac{1}{2} \lo^{ij} \bar{\nabla}_0(\overline{W})_{0ij0} + H(\lo,\W) + (\Rho,\W)
   - \bar{B}_{00} \right) dvol_h.
\end{equation}
\end{rem}

We give an independent proof of this fact.

\begin{proof}[Second proof Remark \ref{CI-Bach}] First, we note that
\begin{align*}
   \left( \int_M \left(\frac{1}{2} \lo^{ij}\bar{\nabla}_0 (\overline{W})_{0ij0} + H (\lo,\W)\right) dvol_h \right)^\bullet[\varphi]
   & = \int_M \lo^{ij} \overline{W}_{kij0} \varphi^k dvol_h = - \int_M \delta (\lo^{ij} \overline{W}_{0ij\cdot}) \varphi dvol_h
\end{align*}
(using \eqref{CI-W-term} and partial integration). Second, we have
$$
    \left( \int_M (\Rho,\W) dvol_h \right)^\bullet[\varphi] = - \int_M (\Hess(\varphi),\W) dvol_h
    = -\int_M \delta \delta (\W) \varphi dvol_h.
$$
Finally, the conformal  transformation law 
$$
    e^{2\varphi} \hat{B}_{ij} = B_{ij} + (n-4) (C_{ijk} + C_{jik}) \varphi^k + (n-4) W_{kijl} \varphi^k \varphi^l
$$
with $C_{ijk} = \nabla_k(\Rho)_{ij} - \nabla_j(\Rho)_{ik}$ implies (for $n=5$)
$$
    \left( \int_M\bar{B}_{00} dvol_h \right)^\bullet[\varphi] = 2 \int_M \bar{C}_{00k} \varphi^k dvol_h
    = - 2 \int_M \delta (\bar{C}_{00\cdot}) \varphi dvol_h.
$$
But, the general formula
$$
    (n-3) C_{ijk} = \Div_1 (W)_{ijk} = \nabla^a (W)_{aijk}
$$
implies (for $n=5$)
$$
    2 \bar{C}_{00k} = \overline{\Div}_1(\overline{W})_{00k} = \bar{\nabla}^a (\overline{W})_{a00k}.
$$
Now
\begin{align*}
     \bar{\nabla}^a (\overline{W})_{a00k} & = \partial^a (\overline{W}_{a00k}) -
      \overline{W}(\bar{\nabla}^a(\partial_a),\partial_0,\partial_0,\partial_k)
      - \overline{W}(\partial_a, \partial_0,\partial_0,\bar{\nabla}^a(\partial_k))\\
      & - \overline{W}(\partial_a,\bar{\nabla}^a(\partial_0),\partial_0,\partial_k)
      - \overline{W}(\partial_a,\partial_0,\bar{\nabla}^a(\partial_0),\partial_k).
\end{align*}
Note that the terms for $a=0$ vanish. Thus, it suffices to let the summation run only over the tangential index
$a=1,\dots,4$. But for such $a$ it holds $\bar{\nabla}^a (\partial_k) = \nabla^a(\partial_k)$ and
$\bar{\nabla}^a(\partial_0) = L^{aj} \partial_j$. Hence
$$
    \bar{\nabla}^a (\overline{W})_{a00k} = \delta (\W)_k
    - \overline{W}_{aj0k} L^{aj} - \overline{W}_{a0jk} L^{aj} = \delta (\W)_k - \overline{W}_{a0jk} L^{aj}
     =  \delta (\W)_k - \overline{W}_{a0jk} \lo^{aj}
$$
and we find
\begin{equation}\label{CWL}
   2 \bar{C}_{00k} = \delta(\W)_k + \overline{W}_{0ajk} \lo^{aj}.
\end{equation}
Therefore, we get the variation formula
$$
     - \left( \int_M\bar{B}_{00} dvol_h \right)^\bullet[\varphi] = \int_M \left(\delta \delta (\W)
     + \delta (\overline{W}_{0ij\cdot} \lo^{ij})\right) \varphi dvol_h.
$$
These conformal variation formulas imply that the conformal variation of \eqref{inv-2} vanishes.
\end{proof}

It is also worth emphasizing the conformal invariance of the second integral in the decomposition \eqref{E-GR2} 
of the functional $\mathcal{E}_{GR}$.

\begin{cor}\label{inv-new} The integral
\begin{equation}\label{INV-3}
    \int_M \left((\lo^2,\Rho) - \frac{1}{2} \J |\lo|^2 - \frac{1}{9} |\delta(\lo)|^2 
    + (\Rho,\W) - \bar{B}_{00}\right) dvol_h
\end{equation}
is conformally invariant.
\end{cor}

This result also follows from the following local fact.

\begin{lem}\label{Invariant-4} The curvature quantity
\begin{equation}\label{INV-2}
    \Jv_4 \st  \frac{2}{9} |\delta(\lo)|^2 - 2 (\lo^2,\Rho) + \J |\lo|^2 - 2 (\Rho,\W) + 2 \bar{B}_{00}
    - \delta \delta (\lo^2) - \delta \delta (\W)
\end{equation}
of the embedding $M^4 \hookrightarrow X^5$ is conformally invariant of weight $-4$, i.e., it holds
$$
    e^{4 \iota^*(\varphi)} \hat{\Jv}_4 = \Jv_4
$$ 
for all $\varphi \in C^\infty(X)$.
\end{lem}

\begin{proof} We prove that the conformal variation of $\Jv_4$ vanishes. We recall that
$$
    (\lo^2,\Rho)^\bullet[\varphi] = - (\lo^2,\Hess(\varphi)), \quad (\J |\lo|^2)^\bullet[\varphi]
    = - |\lo|^2 \Delta(\varphi)
$$
and
$$
    (|\delta(\lo)|^2)^\bullet[\varphi] = 6 (\delta(\lo), \lo d\varphi)
$$
(see \eqref{CV-square}). Hence the conformal variation of the first three terms in \eqref{INV-2} equals
\begin{align*}
    & 2(\lo^2,\Hess(\varphi)) - |\lo|^2 \Delta(\varphi) + \frac{4}{3} (\delta (\lo), \lo d\varphi) \\
    & = - 2 (\delta(\lo^2),d\varphi) + (d (|\lo|^2),d\varphi) + \frac{4}{3} (\lo \delta(\lo),d\varphi)
    + 2 \delta(\lo^2 d\varphi) - \delta (|\lo|^2 d\varphi).
\end{align*}
By \eqref{WN} (for $n=4$), this sum equals
$$
   -2  (\lo^{ij} \overline{W}_{0ij\cdot},d\varphi) + 2 \delta(\lo^2 d\varphi) - \delta (|\lo|^2 d\varphi).
$$
We also recall that
$
   (\delta \delta (\lo^2))^\bullet[\varphi] = 2 \delta (\lo^2 d\varphi) - \delta (|\lo|^2 d\varphi)
$
(see \eqref{var-2}). Thus,
\begin{align}\label{SI-1}
    \left( \frac{2}{9} |\delta(\lo)|^2 - 2 (\lo^2,\Rho) + \J |\lo|^2 - \delta \delta (\lo^2)\right)^\bullet[\varphi]
   & = - 2 (\lo^{ij} \overline{W}_{0ij\cdot},d\varphi) \notag \\
    & = 2 \delta (\lo^{ij} \overline{W}_{0ij\cdot}) \varphi - 2 \delta(\lo^{ij} \overline{W}_{0ij\cdot} \varphi).
\end{align}
Next, we calculate
\begin{align}\label{SI-2}
   (\bar{B}_{00})^\bullet[\varphi] & = 2 \bar{C}_{00k} \varphi^k \notag \\
   & = - 2 \delta (\bar{C}_{00\cdot}) \varphi + 2 \delta (\bar{C}_{00\cdot} \varphi) \notag \\
   & = - \delta \delta (\W) \varphi - \delta (\overline{W}_{0ij\cdot} \lo^{ij}) \varphi
  + 2 \delta(\bar{C}_{00\cdot} \varphi)
\end{align}
using arguments in the proof of Remark \ref{CI-Bach} and \eqref{CWL}. Finally, we find
\begin{align}\label{SI-3}
   (\Rho,\W)^\bullet[\varphi] & = - (\Hess(\varphi),\W) \notag \\
   & = (d\varphi, \delta(\W)) - \delta (\W d\varphi) \notag \\
   & = - \delta \delta (\W) \varphi + \delta (\delta(\W) \varphi) - \delta (\W d\varphi).
\end{align}
Summarizing the results \eqref{SI-1}--\eqref{SI-3} shows that the conformal variation of
\eqref{INV-2} (up to the last term) equals
\begin{align*}
   & 4 \delta(\bar{C}_{00\cdot} \varphi) - 2 \delta(\lo^{ij} \overline{W}_{0ij\cdot} \varphi)
   - 2 \delta (\delta(\W) \varphi) + 2 \delta (\W d\varphi) \\
   & = 2 \delta (\delta (\W) \varphi) + 2 \delta (\lo^{ij} \overline{W}_{0ij\cdot} \varphi)
   - 2 \delta(\lo^{ij} \overline{W}_{0ij\cdot} \varphi) - 2 \delta (\delta(\W) \varphi) + 2 \delta (\W d\varphi)
   \qquad \mbox{(by \eqref{CWL})} \\
   & = 2 \delta (\W d\varphi) \\
   & = (\delta \delta (\W))^\bullet[\varphi].
\end{align*}
This proves that $(\Jv_4)^\bullet[\varphi] = 0$. The proof is complete.
\end{proof}

As a corollary of Lemma \ref{Invariant-4} and the conformal invariance of $\Jv_1$, we obtain the 
following improvement of Remark \ref{CI-Bach}.

\begin{cor}\label{J3} The curvature quantity
\begin{equation}\label{J3-def}
   \Jv_3 \st \Jv_4 - \Jv_1 =
   \frac{1}{2} \lo^{ij} \bar{\nabla}_0(\overline{W})_{0ij0} + H (\lo,\W) + (\Rho,\W) - \bar{B}_{00}
   + \frac{1}{2} \delta \delta (\W)
\end{equation}
is conformally invariant of weight $-4$.
\end{cor}

Note that $\Jv_3$ vanishes if $\overline{W} = 0$.

Finally, we rewrite the second relation in Lemma \ref{GR-I} in terms of the Pfaffian $\Pf_4$ and the invariants 
$\Iv_j$, $\Jv_1$. In particular, this requires expressing $|W|^2$ and $\Jv_3$ in terms of the invariants $\Iv_j$. 

\begin{lem}\label{Weyl-ext} In dimension $n=4$, it holds
\begin{equation}\label{W2-exterior}
    |W|^2 = |\overline{W}|^2 + \frac{7}{3} \Iv_1 - 4 \Iv_2 - 2 \Iv_4 + 4 \bar{\Iv}_5 - 4 \Iv_6.
\end{equation}
\end{lem}

\begin{proof} In general dimensions, the Gauss equation \eqref{Gauss-Weyl} implies
\begin{align*}
   |W|^2 & = |\overline{W} - \frac{1}{2} (\lo \owedge \lo) - \JF \owedge h|^2 \\
   & = |\overline{W}|^2 + \frac{1}{4} |\lo \owedge \lo|^2 + |\JF \owedge h|^2
   - (\lo \owedge \lo,\overline{W}) - 2 (\JF \owedge h,\overline{W}) + (\lo \owedge \lo,\JF \owedge h).
\end{align*}
Now using
\begin{align*}
     |\lo \owedge \lo|^2 & = 8 |\lo|^4 - 8 \tr(\lo^4), \\
     |\JF \owedge h|^2 & = 4(n-2) |\JF|^2 + 4 \tr (\JF)^2, \\
     (\lo \owedge \lo,\overline{W}) & = - 4 \lo^{il} \lo^{jk} \overline{W}_{ijkl}, \\
     (\JF \owedge h,\overline{W}) & = 4 (\JF,\W), \\
     (\lo \owedge \lo,\JF \owedge h) & = - 8 (\JF,\lo^2).
\end{align*}
We apply these results in dimension $n=4$. Then
$$
   |W|^2  = |\overline{W}|^2 + 2 |\lo|^4 - 2 \tr(\lo^4) 
   + 8 |\JF|^2 + 4 \tr(\JF)^2 + 4 \bar{I}_5 - 8 (\JF,\W) - 8 (\JF,\lo^2).
$$
In order to make that relation explicit, we recall that $2 \JF =  \lo^2 - \frac{1}{6} |\lo|^2 h + \W$. 
Hence 
\begin{align*}
    \tr(\JF) & = \frac{1}{6} |\lo|^2, \\
    |\JF|^2 & = - \frac{1}{18} \Iv_1 + \frac{1}{4} \Iv_2 + \frac{1}{4} \Iv_4 + \frac{1}{2} \Iv_6, \\
    (\JF,\W) & = \frac{1}{2} \Iv_4 + \frac{1}{2} \Iv_6, \\
    (\JF,\lo^2)  & = - \frac{1}{12} \Iv_1 + \frac{1}{2} \Iv_2 + \frac{1}{2} \Iv_6.
\end{align*}
Summarizing these results completes the proof.
\end{proof}

Moreover, \cite{AS2} provides the additional relation
\begin{equation}\label{J3-AS}
   \Jv_3 = \frac{1}{2} \bar{\Iv}_5 + \frac{1}{2} \Iv_6 - \frac{1}{4} \Iv_7,
\end{equation}
up to a divergence. Note that this formula implies that the integral of $\Jv_3$ vanishes if $\lo=0$ using Lemma 
\ref{div-term}.

Now combining \eqref{W2-exterior} and \eqref{J3-AS} with the formula
$$
    8 \mathcal{E}_{GR} = \int_M  \left( \Pf_4 - \frac{1}{8} |W|^2 + \frac{1}{2} \Jv_1 - \Jv_3 
    + \frac{1}{12} \Iv_1 - \frac{1}{4} \Iv_2 - \frac{1}{4} \Iv_4 - \frac{1}{2} \Iv_6 \right) dvol_h
$$
(see \eqref{E-GR}) gives the following result.

\begin{cor}\label{EGR-DS}
\begin{equation}\label{EGR-final}
    8 \mathcal{E}_{GR}  = \int_M  \left( \Pf_4 + \frac{1}{2} \Jv_1 - \frac{5}{24} \Iv_1 + \frac{1}{4} \Iv_2 
    - \frac{1}{8} \Iv_3 - \bar{\Iv}_5 - \frac{1}{2} \Iv_6 + \frac{1}{4} \Iv_7 \right) dvol_h.
\end{equation}
\end{cor}

\begin{rem}\label{GR-2} Corollary \ref{EGR-DS} represents the integrand of the functional $\mathcal{E}_{GR}$
(for a hypersurfaces $M^4 \hookrightarrow X^5$)  as a linear combination of the Pfaffian of $M$, local
conformal invariants of the embedding and total divergences. Similarly, the Graham-Reichert functional for a closed 
surface $M^2 \hookrightarrow X^{n+1}$  ($n \ge 2$) is a constant multiple of $\int_M (|H|^2 + \tr_h (\bar{\Rho})) dvol_h$
\cite[Corollary 5.3]{GR}. For $n=2$, this integral equals
$$
   \int_M (H^2 + \bar{\J} - \bar{\Rho}_{00}) dvol_h = \int_M (\J + \frac{1}{2} |\lo|^2) dvol_h
$$
(by the Gauss identity). For $n=3$, it equals
$\int_M (|H|^2 + \bar{\J} - \bar{\Rho}_{00} - \bar{\Rho}_{11}) dvol_h$, where $\{\partial_0,\partial_1\}$ is
an orthonormal basis of the normal space of $M$. This conformal invariant appears in the logarithmic term
of the entanglement entropy in \cite[(1.1),(A6)]{Sol-entropy}. For general $n$, the Gauss equation shows that
$$
   n(|H|^2 + \tr_h (\bar{\Rho}))
   = \scal + (\lo,\lo) + \left(- \bar{R}_{ab}{}^{ba} + \frac{n-2}{n-1} \overline{\Ric}_a^a + (n-2) (H,H) \right)
$$
(if $n \ge 2$). Here $\{\partial_a\}$ is an orthonormal basis of the normal space of $M$. The integral of $\J$ is a
constant multiple of the Euler characteristic of $M$ (Gauss-Bonnet), $(\lo,\lo)$ and the last term in brackets are
local conformal invariants of $M^2 \hookrightarrow X^{n+1}$.
\end{rem}

Another representation of $\mathcal{E}_{GR}$ is given in \cite[Section 4]{BGW-1}. It is a consequence of the theory 
developed in \cite{AGW}.\footnote{Even for $\lo=0$, this formulas differs from \eqref{E-GR}. In view 
of ${\text Wm} = 0$ and $\Iv_7 = 0$ (Lemma \ref{div-term}), the formula in \cite{BGW-1} is equivalent to
$8 \mathcal{E}_{GR} = \int_M \J^2 - |\Rho|^2 + \frac{1}{4} \Iv_4$. On the other hand, \eqref{E-GR} simplifies to
$8 \mathcal{E}_{GR} = \int_M \J^2- |\Rho|^2 - \frac{1}{4} \Iv_4$.}

\section{Energy functionals for a flat background}\label{flat-bg}

Here we take a closer look at the energy functionals
$$
   \mathcal{C}_M \st \int_M \mathcal{C} dvol_h \quad \mbox{and} \quad (\Jv_i)_M \st \int_M \Jv_i dvol_h
$$
for a closed hypersurface $M^4 \hookrightarrow R^5$ in a flat background. By the Gauss equation, all curvature 
data can be expressed in terms of $L$. As a consequence of the conformal invariance of $\mathcal{C}$,  we 
find

\begin{lem}\label{Mob} Assume that $M^4 \hookrightarrow \R^5$. Then
$$
    \mathcal{I}_M \st \int_{M} \left( |dH|^2 - H \tr(\lo^3)  + \frac{1}{2} H^2 |\lo|^2\right) dvol_h
$$
is Möbius invariant, i.e., it holds
$
   \mathcal{I}_{\gamma (M)} = \mathcal{I}_M
$
for all Möbius transformations $\gamma$ of $\R^5$.
\end{lem}


\begin{proof} Let $X = \R^5$. We evaluate the integral
\begin{align*}
   \int_M \mathcal{C} dvol_h = \int_M \left(2 (\lo,\Hess(H)) + 2 H(\lo,\Rho) + 8 (\lo^2,\Rho)
   - 3 \J |\lo|^2 - 3  H^2 |\lo|^2 - H \tr(\lo^3)\right) dvol_h.
\end{align*}
Partial integration and $\delta(\lo) = 3 dH$ (Codazzi-Mainardi) show that
$$
    \int_M (\lo,\Hess(H)) dvol_h = - \int_M (\delta(\lo),dH) dvol_h = - \int_M 3 (dH,dH) dvol_h.
$$
The Gauss equation and the Fialkow equation imply
\begin{align}\label{J-Rho-flat}
    \J = 2 H^2 - \frac{1}{6} |\lo|^2 \quad \mbox{and} \quad
    \Rho = -\frac{1}{2} \lo^2 + \frac{1}{12} |\lo|^2 h + H \lo + \frac{1}{2} H^2 h
\end{align}
and $2H (\lo,\Rho) = - H \tr(\lo^3) + 2H^2 |\lo|^2$. Hence
\begin{align*}
   3 |\lo|^2 \J = 6 H^2 |\lo|^2 - \frac{1}{2} |\lo|^4 \quad \mbox{and} \quad
   8 (\lo^2,\Rho) =  - 4 \tr (\lo^4) + 8 H \tr(\lo^3) + 4 H^2 |\lo|^2 + \frac{2}{3} |\lo|^4.
\end{align*}
These results imply
\begin{align}\label{C-flat}
    \mathcal{C}_M = \int_M \mathcal{C} dvol_h
    & = \int_M \left(- 6 |dH|^2 + 6 H \tr(\lo^3) - 3 H^2 |\lo|^2 + \frac{7}{6} |\lo|^4 -  4\tr(\lo^4)\right) dvol_h
\end{align}
Now the Möbius invariance of $\mathcal{C}_M$ implies the assertion.
\end{proof}

Similar arguments yield

\begin{lem}\label{J1-int} For an embedding $M^4 \hookrightarrow \R^5$, it holds
\begin{align*}
     \int_M \Jv_1 dvol_h & = \int_M \left( 2 |dH|^2 - 2 H \tr(\lo^3) + H^2 |\lo|^2 - \frac{1}{3} |\lo|^4 + \tr(\lo^4) \right) dvol_h
\end{align*}
and
\begin{align*}
     \int_M \Jv_2 dvol_h & = \int_M  \left(|dH|^2 - H \tr(\lo^3) + \frac{1}{2} H^2 |\lo|^2 - \frac{1}{12} |\lo|^4 \right) dvol_h.
\end{align*}
In particular, we find
\begin{equation}\label{diff-flat}
    2 \int_M \Jv_2  dvol_h - \int_M \Jv_1 dvol_h = \int_M  \left(\frac{1}{6} |\lo|^4 - \tr(\lo^4) \right) dvol_h.
\end{equation}
\end{lem}

The relation \eqref{diff-flat} reflects properties of the invariants $\Jv_j$ for general backgrounds. In fact, this relation 
follows by combining Proposition \ref{B-JJD} with Remark \ref{action-local} and $\Iv_5 = \frac{7}{6} \Iv_1 - 2 \Iv_2$
(see \eqref{II5}).

Another calculation using \eqref{J-Rho-flat} shows that
\begin{equation}\label{sigma2}
   \J^2 - |\Rho|^2 = 3 H^4 + H \tr(\lo^3) - \frac{3}{2} H^2 |\lo|^2 + \frac{1}{12} |\lo|^4 - \frac{1}{4} \tr(\lo^4).
\end{equation}
Hence the conformal invariance of the integral $\int_M (\J^2 - |\Rho|^2) dvol_h$  implies that the integral
$$
   \int_M  \left(3 H^4 + H \tr(\lo^3) - \frac{3}{2} H^2 |\lo|^2\right) dvol_h
$$
is Möbius invariant. Combining this with Lemma \ref{Mob} shows that the energy functional
\begin{equation}\label{Guv-corrected}
   \mathcal{E}_{G} \st \frac{1}{4} \int_M \left( |dH|^2 + 3 H^4 - H^2 |\lo|^2 \right) dvol_h
\end{equation}
is Möbius invariant. This result is related to \cite[(61)]{Guven}. Graham and
Reichert \cite{GR} proved that the functional $\mathcal{E}_G$ is a special case of a
global conformal invariant
$$
    \mathcal{E}_{GR} = \frac{1}{4} \int_M (|dH|^2 + \cdots) dvol_h
$$
of $M \hookrightarrow X^5$ which appears in the asymptotic expansion of the renormalized volume of a
minimal hypersurface with boundary $M \hookrightarrow X$ in a Poincar\'e-Einstein background with
conformal infinity $(X,[g])$ (see Section \ref{GR-functional}). Graham and Reichert noticed that the
original calculation in \cite{Guven} dropped a factor of $-2$. The corrected result \cite[(1.1)]{GR}
for $M \hookrightarrow \R^n$ (for $n=5$) states the Möbius invariance of \eqref{Guv-corrected}.\footnote{Note
that \cite{GR} uses a different normalization of $H$.} For more results on the functional
$\mathcal{E}$, we refer to \cite{GR}. The results in \cite{GR} suggest regarding the
functional $\mathcal{E}_G$ as a natural analog of the Willmore functional.

Combining \eqref{C-flat}, \eqref{sigma2} and \eqref{Guv-corrected}, we find the relation
\begin{equation}\label{inv-rel}
   \mathcal{C}_M = - 24 \mathcal{E}_G + 6 \int_M (\J^2 - |\Rho|^2) dvol_h
  + \int_M \left(\frac{2}{3} |\lo|^4 - \frac{5}{2} \tr(\lo^4)\right) dvol_h
\end{equation}
for $M^4 \hookrightarrow \R^5$. In particular, this again shows the Möbius invariance of $\mathcal{E}_G$. 
Together with the Chern-Gauss-Bonnet formula \eqref{CGB-4}, we obtain the relation
$$
    \mathcal{C}_M = 24 \pi^2 \chi(M) - 24 \mathcal{E}_G - \frac{3}{4} \int_M |W|^2 dvol_h
   + \int_M \left(\frac{2}{3} |\lo|^4 - \frac{5}{2} \tr(\lo^4)\right) dvol_h.
$$
The difference $24 \pi^2 \chi(M) - 24 \mathcal{E}_G$ can be expressed in terms of the universal invariant 
${\text Wm}$ (introduced in \cite{BGW-1}); the proof of the following formula is given in Section \ref{invariant}.

\begin{lem}\label{Guven-Wm}  $24 \pi^2 \chi(M) - 24 \mathcal{E}_G 
= \int_M ({\text Wm} + \frac{11}{6} |\lo|^4 - \frac{5}{2} \tr(\lo^4)) dvol_h$.
\end{lem}

Finally, we observe that the relation \eqref{inv-rel} generalizes to embeddings $M^4 \hookrightarrow S^5$ in the
round sphere $S^5$ if
\begin{equation}\label{energy-sphere}
    \mathcal{E}_G \st  \frac{1}{4} \int_M \left( |dH|^2 + 3 H^4 - H^2 |\lo|^2 + 6 H^2 + 3 \right) dvol_h.
\end{equation}
By \cite[(1.2)]{GR}, this energy functional again is a special case of the
Graham-Reichert energy functional $\mathcal{E}_{GR}$. In the present case, it holds
$\bar{\Rho} = \frac{1}{2} \bar{g}$ (with $\bar{g}$ being the round metric on $S^5$)
and $\bar{\J} = \frac{5}{2}$. Then $\bar{\Rho}_{00} = \frac{1}{2}$. The Gauss
identity shows that $\J = \J_{\text flat} + 2$, where $\J_{\text flat} = 2 H^2 -
\frac{1}{6} |\lo|^2$ is defined by the same formula as in the case of a flat
background. Similarly, the Fialkow equation shows that $\Rho = \Rho_{\text flat} +
\frac{1}{2} h$, where $\Rho_{\text flat}$ is defined by the same formula
\eqref{J-Rho-flat} as in the case of a flat background. Now we calculate
\begin{align*}
    \int_M \mathcal{C} dvol_h & = \int_M \left( - 6 |dH|^2 + 2 H(\lo,\Rho)
    + 8 (\lo^2,\Rho) - |\lo|^2 - 3 \J |\lo|^2 - 3  H^2 |\lo|^2 - H \tr(\lo^3) \right) dvol_h \\
    & = \int_M \mathcal{C}_{\text flat} + 4 |\lo|^2 - |\lo|^2 - 6 |\lo|^2 = \int_M \mathcal{C}_{\text flat} -3 |\lo|^2,
\end{align*}
where $\mathcal{C}_{\text flat}$ is given by the integrand in \eqref{C-flat}, and
\begin{align*}
   \J^2  - |\Rho|^2 & = (\J_{\text flat} + 2)^2 - |\Rho_{\text flat} + 1/2 h|^2 \\
   & = \J_{\text flat}^2 - |\Rho_{\text flat}|^2 + 4 \J_{\text flat} + 4 - \J_{\text flat} - 1 \\
   & =  \J_{\text flat}^2 - |\Rho_{\text flat}|^2 + 6 H^2 - \frac{1}{2}  |\lo|^2 + 3.
\end{align*}
These results imply \eqref{inv-rel} for $M^4 \hookrightarrow S^5$. In particular,
this again shows the Möbius invariance of the energy \eqref{energy-sphere}.

\section{The singular Yamabe energy of a four-dimensional hypersurface $M$}\label{SYE}


The main result of this section is a formula for the singular Yamabe energy
$$
   \mathcal{E}_M \st \int_M V_4 dvol_h
$$
of $M^4 \hookrightarrow X^5$ for a general background metric $g$ (Theorem \ref{V4-final}). Here
$V_4$ denotes the fourth-order singular Yamabe renormalized volume coefficient. The final result will
confirm the relation
\begin{equation}\label{E-Q}
    16 \mathcal{E}_M = \int_{M^4} {\bf Q}_4 dvol_h
\end{equation}
of conformal invariants proved in \cite{GW-RV, JO1}.

We first recall the definition of the coefficient $V_n$ for a hypersurface $M^n \hookrightarrow X^{n+1}$. In normal
geodesic coordinates, the metric $g$ takes the form $g = dr^2 + h_r$ with a one-parameter family $h_r$. Then the
expansion
$$
    v(r) = dvol_{h_r}/dvol_h = \sum_{k \ge 0} r^k v_k
$$
defines the volume coeffcients $v_k \in C^\infty(M)$  (see \eqref{v-coeff}).

We recall from Section \ref{Laplace-scattering} that $\sigma$ solves the singular Yamabe problem for the
hypersurface $M^n \hookrightarrow X^{n+1}$ with the background metric $g$ if the scalar curvature of $\sigma^{-2} g$
equals $-n(n+1)$.

\begin{defn}\label{def-Vn} Let
$$
    \sigma (r) = \sum_{k \ge 1} \sigma_{(k)} r^k = r +  r^2 \sigma_{(2)} +  r^3  \sigma_{(3)}+ \cdots
$$
be the expansion of the solution $\sigma(r)$ of the singular Yamabe problem of $M^n \hookrightarrow X^{n+1}$ in normal
geodesic coordinates. Let $V_n$ be the coefficient of $r^n$ in the expansion of the function
$$
    (1 + r \sigma_{(2)} +  r^2  \sigma_{(3)}+ \cdots )^{-(n+1)} v(r).
$$
\end{defn}

Note that Definition \ref{def-Vn} implies that the expansion of
\begin{equation*}
    dvol_{\sigma^{-2} g} = r^{-(n+1)} (1+ \sigma_{(2)} r + \sigma_{(3)} r^2 + \cdots )^{-(n+1)} v(r) dr dvol_h
\end{equation*}
involves a term $V_n r^{-1} dr dvol_h$. By integration, this shows that the total integral of $V_n$ defines the coefficient of
$\log(\varepsilon)$ in the expansion of the volume. The integral $\int_M V_n dvol_h$ is the singular Yamabe energy \cite{G-SY}.

The following result describes the coefficients $V_k$ for $k \le 4$ in terms of $\sigma_{(k)}$ and $v_k$ for $k \le 4$
in the respective critical dimensions. It directly follows from the definition.

\begin{lem}\label{Vn-crit} In the respective critical dimensions, it holds 
\begin{align*}
   V_2 & = 6 \sigma_{(2)}^2 - 3 \sigma_{(3)} - 3\sigma_{(2)} v_1 + v_2, \\
   V_3 & = -20 \sigma_{(2)}^3 + 20 \sigma_{(2)} \sigma_{(3)} - 4 \sigma_{(4)}
   + 10 \sigma_{(2)}^2 v_1 - 4 \sigma_{(3)} v_1 - 4 \sigma_{(2)} v_2 + v_3
\end{align*}
and
\begin{align*}
   V_4 & = 70 \sigma_{(2)}^4 - 105 \sigma_{(2)}^2 \sigma_{(3)}+ 15 \sigma_{(3)}^2 + 30 \sigma_{(2)}\sigma_{(4)}
   - 5 \sigma_{(5)} \\
   & - 35 \sigma_{(2)}^3 v_1 + 30 \sigma_{(2)}\sigma_{(3)} v_1 - 5 \sigma_{(4)} v_1 + 15 \sigma_{(2)}^2 v_2
   - 5 \sigma_{(3)} v_2 - 5 \sigma_{(2)} v_3 + v_4.
\end{align*}
\end{lem}

In order to determine the coefficients $V_k$ for $k \le 4$, we need explicit formulas for the coefficients
$\sigma_{(k)}$ for $k \le 5$ and $v_k$ for $k \le 4$. We first display such formulas for the coefficients
$\sigma_{(k)}$.

\begin{lem}[\cite{JO1}]\label{sigma-234} In general dimensions, it holds
\begin{align*}
   \sigma_{(2)} & = \frac{1}{2n} v_1, \\
   \sigma_{(3)} & = \frac{2}{3(n-1)} v_2 - \frac{1}{3n} v_1^2 + \frac{1}{3(n-1)} \bar{\J}.
\end{align*}
and
\begin{align*}
   \sigma_{(4)} & = \frac{3}{4(n-2)} v_3 - \frac{9n^2-20n+7}{12n(n-1)(n-2)} v_1 v_2
   + \frac{6n^2-11n+1}{24n^2(n-2)} v_1^3 \notag \\
   & + \frac{2n-1}{6n(n-1)(n-2)} v_1 \bar{\J} + \frac{1}{4(n-2)} \bar{\J}' + \frac{1}{8n (n-2)}  \Delta (v_1).
\end{align*}
\end{lem}

Note that $3 \sigma_{(3)} =  2 v_2 + \cdots$ for $n=2$ and $4 \sigma_{(4)} = 3 v_3 + \cdots$ for $n=3$.

Note also that $\sigma_{(4)}$ has a simple pole at $n=2$ with $\res_{n=2}(\sigma_{(4)}) \sim \B_2$.

In connection with the discussion of $V_4$, we shall apply the following consequences for $n=4$.

\begin{cor}\label{sigma-234-4} In the critical dimension $n=4$, it holds
\begin{align*}
    \sigma_{(2)} & = \frac{1}{8} v_1, \\
    \sigma_{(3)} & = -\frac{1}{12} v_1^2 + \frac{2}{9} v_2 + \frac{1}{9} \bar{\J}
\end{align*}
and
\begin{align*}
     \sigma_{(4)} & = \frac{53}{768} v_1^3 - \frac{71}{288} v_1 v_2 + \frac{3}{8} v_3
     + \frac{7}{144} \bar{\J} v_1 + \frac{1}{8} \bar{\J}' + \frac{1}{8} \Delta (\sigma_{(2)}).
\end{align*}
\end{cor}

Finally, we need the following formula for $\sigma_{(5)}$.

\begin{lem}[\cite{JO1}] \label{sigma-5} In general dimensions, it holds
\begin{align*}
    \sigma_{(5)} & = - \frac{n+1}{10(n-3)} |d\sigma_{(2)}|^2 \\
   & + \frac{1}{5(n-3)} \Delta'(\sigma_{(2)})
   + \frac{1}{5(n-3)} \Delta(\sigma_{(3)}) + \frac{3n-1}{20(n-3)(n-2)n} \Delta(\sigma_{(2)}) v_1 \\
   & + \frac{1}{10(n-3)} \bar{\J}''  + \frac{n-1}{4(n-3)(n-2)n} \bar{\J}' v_1
   + \frac{2(3n-5)}{15(n-3)(n-1)^2} \bar{\J} v_2 \\
   & - \frac{4n-3}{20 (n-2)(n-1)n} \bar{\J} v_1^2 + \frac{1}{30(n-1)^2} \bar{\J}^2 \\
   & + \frac{48n^4-247n^3+387n^2-179n+3}{60(n-3)(n-2)(n-1)n^2} v_1^2 v_2
   - \frac{2(3n^2-11n+10)}{15(n-3)(n-1)^2} v_2^2 \\
   & - \frac{24n^4-110n^3+133n^2-24n-3}{120 (n-3)(n-2)n^3} v_1^4
   - \frac{16n^2-53n+27}{20(n-3)(n-2)n} v_1 v_3 + \frac{4}{5(n-3)} v_4.
\end{align*}
\end{lem}

Note that $\sigma_{(5)}$ has a simple pole at $n=3$ with $\res_{n=3}(\sigma_{(5)}) \sim \B_3$. 


In particular, we obtain

\begin{cor}\label{sigma-5-4} In the critical dimension $n=4$, it holds
\begin{align*}
    24 \sigma_{(5)}  & = - \frac{1133}{640} v_1^4 - \frac{213}{20} v_1 v_3 + \frac{653}{80} v_1^2 v_2
   - \frac{224}{45} v_2^2 + \frac{96}{5} v_4 \\
    & - \frac{13}{20} \bar{\J} v_1^2 + \frac{112}{45} \bar{\J} v_2  + \frac{9}{4} \bar{\J}' v_1
   + \frac{4}{45} \bar{\J}^2 + \frac{12}{5} \bar{\J}'' \\
   & + \frac{24}{5} \Delta(\sigma_{(3)})
   + \frac{33}{20} v_1 \Delta(\sigma_{(2)}) + \frac{24}{5} \Delta'(\sigma_{(2)}) - 12 |d\sigma_{(2)}|^2.
\end{align*}
\end{cor}

Note that $5 \sigma_{(5)} = 4 v_4 + \cdots$ for $n=4$.

Next, formulas for the volume coefficients $v_k$ for $k \le 4$ in terms of the curvature of the background metric
$g$ were displayed in Lemma \ref{v-coeff-3} and Lemma \ref{v4-form}.

These results imply

\begin{cor}\label{sigma-34} In general dimensions, it holds $2 \sigma_{(2)} = H$ and
$$
   3(n-1) \sigma_{(3)} = - \overline{\Ric}_{00} - |\lo|^2 + \bar{\J} = -(n-1) \bar{\Rho}_{00} - |\lo|^2.
$$
Hence
$$
   3 \sigma_{(3)} = \J - \bar{\J} - \frac{1}{2(n-1)} |\lo|^2 - \frac{n}{2} H^2.
$$
Moreover, we have
$$
    24 \sigma_{(4)} = -3 \bar{\nabla}_0{\overline{\Ric}}_{00} + \bar{\J}' + 6 (\lo,\bar{\G}) + 3 \Delta (H) + 6 \tr (\lo^3)
    +13 H |\lo|^2 - 7H \overline{\Ric}_{00} + 10 H \bar{\J}
$$
for $n=3$.
\end{cor}

\begin{proof} The formula for $\sigma_{(2)}$ is obvious from Lemma \ref{v-coeff-3} and Lemma \ref{sigma-234}.
Similarly, these results imply
$$
   3(n-1) \sigma_{(3)} = (-\overline{\Ric}_{00} - |\lo|^2 + n(n-1) H^2) - n(n-1) H^2 + \bar{\J}.
$$
This identity simplifies the second claim. Finally, a direct calculation yields the third formula.
\end{proof}

The formulas in Corollary \ref{sigma-34} are equivalent to \cite[(2.16), (2.17), (2.19)]{GG}.

Next, we use the above results to find explicit formulas for the coefficients $V_k$ ($k \le 4$).

First of all, combining Lemma \ref{Vn-crit} with Lemma \ref{v-coeff-3} and Corollary \ref{sigma-34} gives

\begin{lem}\label{V2-0} Let $n=2$. Then
\begin{align*}
   2 V_2 & = 12 \sigma_{(2)}^2 - 6 \sigma_{(3)} - 6 \sigma_{(2)} v_1 + 2 v_2 \\
   & = - 2 v_2 + \frac{1}{4} v_1^2 - 2 \bar{\J} \\
   & = \overline{\Ric}_{00} + |\lo|^2 - H^2 -2 \bar{\J}.
\end{align*}
\end{lem}

\begin{cor}\label{V2-1} Let $n=2$. Then
\begin{equation*}
   V_2 = \frac{1}{4} |\lo|^2 - \frac{1}{2} \J.
\end{equation*}
Hence
\begin{equation*}\label{V2-final}
   \int_M V_2 dvol_h = \frac{1}{4} \int_M |\lo|^2 dvol_h - \pi \chi(M).
\end{equation*}
\end{cor}

\begin{proof} By the Gauss equation for $\bar{\J}$, we have
$$
    \frac{1}{4} |\lo|^2 - \frac{1}{2} \J
   = \frac{1}{4} |\lo|^2 - \frac{1}{2} (\bar{\J} - \bar{\Rho}_{00} - \frac{1}{2} |\lo|^2 + H^2).
$$
But $\bar{\Rho}_{00} = \overline{\Ric}_{00} - \bar{\J}$. This implies the first assertion. The second claim follows
from the Gauss-Bonnet formula.
\end{proof}

In particular, the total integral of $V_2$ is conformally invariant.

Similarly, Lemma \ref{Vn-crit} implies

\begin{lem}\label{V3-1} Let $n=3$. Then
$$
    V_3 = -20 \sigma_{(2)}^3 + 20 \sigma_{(2)} \sigma_{(3)} - 4 \sigma_{(4)} + 10 \sigma_{(2)}^2 v_1
   - 4 \sigma_{(3)} v_1 - 4 \sigma_{(2)} v_2 + v_3.
$$
\end{lem}

Evaluation of this formula using the above results yields

\begin{lem}\label{V3} Let $n=3$. Then
\begin{align}\label{V3-final}
    6 V_3 & = - \frac{8}{9} v_1^3 + 4  v_1 v_2 - 12 v_3 - 4  v_1 \bar{\J} - 6 \bar{\J}' -  \Delta (v_1) \notag \\
    & = 2 \bar{\nabla}_0(\overline{\Ric})_{00} - 6 \bar{\J}'  - 4 (\lo,\bar{\G}) + 8 H \overline{\Ric}_{00}
    - \Delta (v_1)  - 12 H \bar{\J} - 4 \tr(\lo^3).
\end{align}
\end{lem}

Next, we simplify \eqref{V3-final}. The second Bianchi identity enables us to remove the normal derivatives in
\eqref{V3-final}. Let $\bar{G} = \overline{\Ric} - \frac{1}{2} \overline{\scal} \bar{g}$ be the Einstein tensor of the background
metric $\bar{g} $ on $X$. The argument rests on the identity
\begin{equation}\label{Bianchi-2}
     \bar{\nabla}_0(\overline{\Ric})_{00} - 3 \bar{\J}'  =
     \bar{\nabla}_0(\bar{G})_{00} \stackrel{!}{=} - \delta (\overline{\Ric}_0) - 3 H \overline{\Ric}_{00} + (L,\overline{\Ric})
\end{equation}
(see Lemma \ref{Nabla-1G}).

\begin{cor}\label{V3-ex-total} Let $n=3$. Then it holds
\begin{equation*}\label{V3-final-2}
   6 V_3 = - 2 \delta \delta (\lo) - 4 (\lo,\JF) + \Delta (H).
\end{equation*}
Thus, $V_3$ is a sum of a linear combination of the local conformal invariant $(\lo,\JF)$ and some divergence terms.
In particular,
\begin{equation*}
    \int_M V_3 dvol_h = -\frac{2}{3} \int_M (\lo,\JF) dvol_h.
\end{equation*}
\end{cor}

Corollary \ref{V3-ex-total} implies the conformal invariance of the total integral of $V_3$ using that of the
Fialkow tensor $\JF$.

\begin{proof} By \eqref{Bianchi-2}, the formula \eqref{V3-final} simplifies to
$$
    6 V_3 = - 2 \delta (\overline{\Ric}_0) + 2 H \overline{\Ric}_{00} + 2 (L,\overline{\Ric}) - 4 (\lo,\bar{\G})
   - 12 H \bar{\J} - 4 \tr(\lo^3) - 3 \Delta (H).
$$
By $\overline{\Ric} =  2 \bar{\Rho} + \bar{\J} h$ and $\bar{\G} = \bar{\Rho} + \bar{\Rho}_{00} h + \W$, we get
$$
    2 H \overline{\Ric}_{00} + 2 (L,\overline{\Ric}) - 4 (\lo,\bar{\G}) - 12 H \bar{\J} = - 4(\lo,\W).
$$
Therefore, we obtain
$$
   6 V_3 = -2 \delta  (\overline{\Ric}_0) - 4(\lo,\W)  - 4 \tr(\lo^3) - 3 \Delta (H).
$$
On the other hand, \eqref{CM-trace} gives
$
   \delta \delta (\lo) = \delta (\overline{\Ric}_0) +  2 \Delta(H).
$
Hence
$$
  - 2 \delta \delta (\lo) - 4 (\lo,\JF) + \Delta (H) = - 2\delta (\overline{\Ric}_0) - 4 \Delta(H) - 4 \tr(\lo^3) - 4 (\lo,\W) + \Delta(H).
$$
This completes the proof.
\end{proof}

Now we turn to the discussion of $V_4$ in dimension $n=4$. Lemma \ref{Vn-crit} gives

\begin{lem}\label{V4-sigma} Let $n=4$. Then
\begin{align*}
   V_4 & = 70 \sigma_{(2)}^4 - 105 \sigma_{(2)}^2 \sigma_{(3)}+ 15 \sigma_{(3)}^2 + 30 \sigma_{(2)}\sigma_{(4)}
   - 5 \sigma_{(5)} \\
   & - 35 \sigma_{(2)}^3 v_1 + 30 \sigma_{(2)}\sigma_{(3)} v_1 - 5 \sigma_{(4)} v_1 + 15 \sigma_{(2)}^2 v_2
   - 5 \sigma_{(3)} v_2 - 5 \sigma_{(2)} v_3 + v_4.
\end{align*}
\end{lem}

The evaluation of this formula using Corollary \ref{sigma-234-4} and Corollary \ref{sigma-5-4} yields

\begin{lem}\label{V4-v} Let $n=4$. Then
\begin{align*}
   24 V_4 & = \frac{981}{256} v_1^4 - \frac{159}{8} v_1^2 v_2 + 16 v_2^2 + 27 v_1 v_3 - 72 v_4 \\
   & + 4 \bar{\J}^2 + \frac{3}{4} \bar{\J} v_1^2 - 8 \bar{\J} v_2 - 15 \bar{\J}' v_1 - 12 \bar{\J}'' \\
   & - 12 v_1 \Delta (\sigma_{(2)}) - 24 \Delta(\sigma_{(3)}) - 24 \Delta'(\sigma_{(2)}) + 60 |d\sigma_{(2)}|^2.
\end{align*}
\end{lem}

Finally, the evaluation of the latter result using the formulas for the volume coefficients $v_k$ in Lemma \ref{v-coeff-3}
and Lemma \ref{v4-form} gives

\begin{prop}\label{E1} Let $n=4$. Then
\begin{align*}
   24 V_4 & = 3 \bar{\nabla}^2_0(\overline{\Ric})_{00} - 12 \bar{\J}'' + 24 H \bar{\nabla}_0 (\overline{\Ric})_{00}
   - 60 H \bar{\J}' - 6 \lo^{ij} \bar{\nabla}_0(\bar{R})_{0ij0} \\
   &  - 48 H \Delta (\sigma_{(2)}) - 24 \Delta (\sigma_{(3)}) - 24 \Delta' (\sigma_{(2)}) + 60 |d \sigma_{(2)}|^2\\
   & + 6 |\bar{\G}|^2  - 5 (\overline{\Ric}_{00})^2  + 4 \overline{\Ric}_{00} \bar{\J}  + 4 \bar{\J}^2  \\
   & + 24 (\lo^2,\bar{\G}) - 10 |\lo|^2 \overline{\Ric}_{00} + 4 |\lo|^2 \bar{\J}
   - 12 H (\lo,\bar{\G}) + 27 H^2 \overline{\Ric}_{00}
   - 36 H^2 \bar{\J}\\
   & + 18 \tr(\lo^4) + 12 H \tr(\lo^3) - 9 H^2 |\lo|^2 - 5 |\lo|^4 + 9 H^4.
\end{align*}
\end{prop}

We continue with a simplification of the formula in Proposition \ref{E1}.

First, we again use the second Bianchi identity to remove the second-order normal derivatives of the Ricci tensor.
Lemma \ref{Nabla-2G} implies that
$$
   3 \bar{\nabla}^2_0(\overline{\Ric})_{00} - 12 \bar{\J}'' + 24 H \bar{\nabla}_0 (\overline{\Ric})_{00} - 60 H \bar{\J}'
   = 3  \bar{\nabla}_0^2 (\bar{G})_{00} + 24 H \bar{\nabla}_0 (\overline{\Ric})_{00} - 60 H \bar{\J}'
$$
equals the sum of
\begin{align}\label{S1}
    & - 15 H \bar{\nabla}_0 (\overline{\Ric})_{00} + 24 H \bar{\J}' + 24 H \bar{\nabla}_0 (\overline{\Ric})_{00} - 60 H \bar{\J}' \notag \\
    & = 9  H \bar{\nabla}_0 (\overline{\Ric})_{00} - 36  H \bar{\J}' \notag \\
    & = 9 H  \bar{\nabla}_0 (\bar{G})_{00} \notag \\
    & = - 9 H \delta (\overline{\Ric}_0) -  36H^2 \overline{\Ric}_{00} + 9 H (L,\overline{\Ric})
\end{align}
and
\begin{align}\label{S2}
     & 6 (\lo, \nabla (\overline{\Ric}_0)) - 3 \delta (\bar{\nabla}_0(\overline{\Ric})_{0})
    + 3 (\lo, \bar{\nabla}_0(\overline{\Ric}))  \notag \\
     & + 3H \delta (\overline{\Ric}_0) - 9 (dH,\overline{\Ric}_0) + 6 (\delta(\lo),\overline{\Ric}_0)
     - 3 \delta ((\lo \overline{\Ric})_{0})  \notag \\
     & + 3 |L|^2  \overline{\Ric}_{00} - 3 (L^2, \overline{\Ric}) + 3 (\overline{\Ric}_{00})^2 - 3 (\bar{\G}, \overline{\Ric}).
\end{align}

Thus, we have proved

\begin{prop}\label{V4-int} $24 V_4$ equals the sum of \eqref{S1}, \eqref{S2},
$$
  - 6 \lo^{ij} \bar{\nabla}_0(\bar{R})_{0ij0},
$$
\begin{align}\label{S3}
   & 60 |d \sigma_{(2)}|^2 - 48 H \Delta (\sigma_{(2)}) - 24 \Delta' (\sigma_{(2)})  - 24 \Delta (\sigma_{(3)})
\end{align}
and
\begin{align*}\label{S4}
   & 6 |\bar{\G}|^2  - 5 (\overline{\Ric}_{00})^2  + 4 \overline{\Ric}_{00} \bar{\J}  + 4 \bar{\J}^2   \notag  \\
   & + 24 (\lo^2,\bar{\G}) - 10 |\lo|^2 \overline{\Ric}_{00} + 4 |\lo|^2 \bar{\J}
   - 12 H (\lo,\bar{\G}) + 27 H^2 \overline{\Ric}_{00} - 36 H^2 \bar{\J} \notag \\
   & + 18 \tr(\lo^4) + 12 H \tr(\lo^3) - 9 H^2 |\lo|^2 - 5 |\lo|^4 + 9 H^4.
\end{align*}
\end{prop}

Since we are only interested in the total integral of $V_4$, we may ignore the total divergences in Proposition \ref{V4-int}. These are the terms $\delta (\bar{\nabla}_0(\overline{\Ric})_{0})$, $\delta ((\lo \overline{\Ric})_{0})$
in \eqref{S2} and $\Delta(\sigma_{(3)})$ in \eqref{S3}.  Furthermore, partial integration shows that
$$
   \int_M- 9 H \delta (\overline{\Ric}_0) - 9 (dH,\overline{\Ric}_0)  = 0
$$
and
$$
   \int_M 6 (\lo, \nabla (\overline{\Ric}_0)) +  6 (\delta(\lo),\overline{\Ric}_0) = 0.
$$
Therefore, we may omit these four terms in \eqref{S1} and \eqref{S2}.

Next, we evaluate the terms in \eqref{S3}. We recall that we omit the term $\Delta (\sigma_{(3)})$. The
variation formula
$$
   \Delta'(u) = - 2(L,\Hess(u)) - 2 (\delta(L),du) + (d \tr(L),du)
$$
implies
$$
   \Delta'(u) = - 2(L,\Hess(u)) - 4 (dH,du) - 2 (\overline{\Ric}_0,du)
$$
using $\delta(L) = 4 dH + 3 \bar{\Rho}_0$ (Codazzi-Mainardi). Hence
\begin{align*}
   & 60  |d \sigma_{(2)}|^2 - 48 H \Delta (\sigma_{(2)}) - 24 \Delta' (\sigma_{(2)}) \\
   & = 15 |dH|^2 - 24 H \Delta (H) + 24 ((L,\Hess(H)) + 2 |dH|^2 + (\overline{\Ric}_0,dH))
\end{align*}
By partial integration, the integral of this sum equals
$$
   \int_M(87 |dH|^2 - 24 (\delta(L),dH) + 24 (\overline{\Ric}_0,dH)) dvol_h = - 9 \int_M |dH|^2 dvol_h
$$
again using $\delta(L) = 4 dH + 3 \bar{\Rho}_0$. On the other hand, the integrals of \eqref{S1} and \eqref{S2} contribute
$$
   -3 \int_M (dH, \overline{\Ric}_0) dvol_h.
$$
Together with the above terms, this gives
$$
   \int_M -3 (dH, \overline{\Ric}_0) - 9 |dH|^2  dvol_h = -3 \int_M (dH, \delta(\lo)) dvol_h = 3 \int_M (\lo,\Hess(H)) dvol_h.
$$

Now simplification of the remaining terms proves

\begin{theorem}\label{V4-final} Let $n=4$. Then
\begin{align*}
   & 8 \int_M V_4 dvol_h = \int_M \left( \J^2 - |\Rho|^2 + \frac{9}{4} |\W|^2 \right) dvol_h \\
   & + \int_M \left( (\lo, \bar{\nabla}_0 (\bar{\Rho})) - 2 \lo^{ij} \bar{\nabla}_0 (\overline{W})_{0ij0}
   +(\lo,\Hess(H)) + H(\lo,\Rho) - \frac{9}{2} H (\lo,\W) \right)  dvol_h \\
   & + \int_M \left( 4 (\lo^2,\Rho) - |\lo|^2 \bar{\Rho}_{00} - \frac{3}{2} \J |\lo|^2 
   - \frac{3}{2} H^2 |\lo|^2 - \frac{1}{2} H \tr(\lo^3) \right) dvol_h \\
   &  + \int_M \left( \frac{21}{2} (\lo^2,\W) + \frac{33}{4} \tr(\lo^4)- \frac{7}{3} |\lo|^4 \right) dvol_h.
\end{align*}
\end{theorem}

\begin{proof} We first verify the first two terms in the second line. The terms in \eqref{S2} and \eqref{S0} yield
the contributions
\begin{align*}
   & - 2 \lo^{ij} \bar{\nabla}_0(\bar{R})_{0ij0} + (\lo,\bar{\nabla}_0(\overline{\Ric})) \\
   & = - 2 \lo^{ij} \bar{\nabla}_0(\overline{W})_{0ij0} - \frac{2}{3} (\lo,\bar{\nabla}_0(\overline{\Ric}))
    + (\lo,\bar{\nabla}_0(\overline{\Ric})) \\
   & =  - 2 \lo^{ij} \bar{\nabla}_0(\overline{W})_{0ij0} + (\lo,\bar{\nabla}_0(\bar{\Rho}))
\end{align*}
using the relation
$$
    3 \lo^{ij} \bar{\nabla}_0(\bar{R})_{0ij0} = 3  \lo^{ij} \bar{\nabla}_0(\overline{W})_{0ij0}
   + (\lo,\bar{\nabla}_0(\overline{\Ric})).
$$
The remaining terms follow by direct calculation using the following identities. First, we note that $8 \bar{\J}
= \overline{\scal}$ and $\overline{\Ric} = 3 \bar{\Rho} + \bar{\J} h$. Now it holds
$$
    \bar{\J} = \J + \bar{\Rho}_{00} + \frac{1}{6} |\lo|^2 - 2 H^2
$$
(Gauss equation) and
$$
    \bar{\Rho} = \Rho - H \lo - \frac{1}{2} H^2 h + \frac{1}{2} \left(\lo^2 - \frac{1}{6} |\lo|^2 h + \W \right)
$$
(Fialkow equation). Finally, $\overline{\Ric}_{00} = 3\bar{\Rho}_{00} + \bar{\J}$ and
$\bar{\G} = \bar{\Rho} + \bar{\Rho}_{00} h + \W$. Using these identities, we express all terms in terms of
$\J, \Rho, \bar{\Rho}_{00}, \W$ and $H, \lo$. We omit the details.
\end{proof}

Comparing Corollary \ref{Q4-g-int} with Theorem \ref{V4-final} confirms the relation \eqref{E-Q}.
Alternatively, Theorem \ref{V4-final} and the relation \eqref{E-Q} confirm Corollary \ref{Q4-g-int}.

\section{Appendix}\label{App}

This appendix  contains the following additional issues.
\begin{itemize}
\item A direct check of the conformal covariance of the operator ${\bf P}_4$ in Theorem \ref{main}.
\item A brief discussion of Deser-Schwimmer type decompositions of conformal invariants of 
hypersurfaces $M^4 \hookrightarrow X^5$ in geometric analysis and physics.
\item A proof of a decomposition of the local conformal invariant ${\text Wm}$ (introduced in \cite{BGW-1}) 
in terms of basic local conformal invariants (Proposition \ref{WM-deco}). 
\item A proof of the relation \eqref{diff-rel}.
\item A proof of the equivalence of the formula for ${\bf Q}_4$ in Theorem \ref{main} to a formula 
in \cite{BGW-1} (at least up to terms which are quartic in $L$).
\end{itemize}

\subsection{The conformal covariance of ${\bf P}_4$}\label{covariance-direct}

The formulation of Theorem \ref{main} contains the claim that the displayed formula defines a conformally covariant
operator. Here we verify this fact by a direct calculation extending the arguments at the end of Section \ref{TUH}.

First, we recall the conformal transformation laws
\begin{equation}\label{conform-div-forms}
    e^{(\lambda+2)\varphi} \hat{\delta} (e^{-\lambda \varphi} \omega) = \delta (\omega)
   + (n\!-\!2\!-\!\lambda) (d\varphi, \omega)
\end{equation}
for $\omega \in \Omega^1(M)$ and
\begin{equation}\label{CTL-Laplace}
    e^{(\lambda+2)\varphi} \hat{\Delta} (e^{-\lambda \varphi} f) = \Delta (f) - \lambda \delta (f d \varphi)
   + (n\!-\!2\!-\!\lambda)(df,d\varphi) - \lambda(n\!-\!2\!-\!\lambda) |d\varphi|^2 f
\end{equation}
for $f \in C^\infty(M)$. Moreover, it holds
\begin{equation}\label{conform-div-BLF}
    e^{(\lambda+2)\varphi} \hat{\delta} (e^{-\lambda \varphi} b) = \delta (\omega)
   + (n\!-\!2\!-\!\lambda) \iota_{\grad(\varphi)}(b) - \tr(b) d\varphi = \delta (\omega)
   + (n\!-\!2\!-\!\lambda) b d\varphi - \tr(b) d\varphi
\end{equation}
for symmetric bilinear forms $b$ and $\lambda \in \R$. Here we use the same symbol for a
bilinear form and the corresponding endomorphism on $\Omega^1(M)$.

Now, by the discussion at the end of Section \ref{TUH}, it suffices to prove that the operator
\begin{align*}
   f & \mapsto \delta \left( 4 \tfrac{3n-5}{n-2} \lo^2 + \tfrac{n^2-12n+16}{2(n-1)(n-2)} |\lo|^2 h \right) df \\
   & + \left(\tfrac{n}{2}-2\right)
   \Big( \tfrac{2(n-1)}{(n-3)(n-2)} \delta \delta (\lo^2) + \tfrac{4}{n-3} \delta (\lo \delta(\lo)) + \tfrac{3n-4}{2(n-1)(n-2)} \Delta (|\lo|^2) \\
   & + 2 (\lo,\Hess(H))  + 2 \lo^{ij} \bar{\nabla}_0(\bar{\Rho})_{ij} - \tfrac{4}{n-3} \lo^{ij} \bar{\nabla}_0 (\overline{W})_{0ij0}
   - \tfrac{2(n-1)^2}{(n-3)(n-2)} H (\lo,\W) \\
   & - \tfrac{2(n^2-9n+12)}{(n-3)(n-2)} (\lo^2,\Rho) - \tfrac{n^3-5n^2+18n-20}{2(n-3)(n-2)(n-1)} \J |\lo|^2 + 2 H (\lo,\Rho)
   - 2 |\lo|^2 \bar{\Rho}_{00} \\
   & - 3 H^2 |\lo|^2 - \tfrac{2(n-3)}{n-2} H \tr(\lo^3) \Big) f
\end{align*}
is conformally covariant. We denote this operator by ${\bf R}_4 = {\bf r}_4 + {\bf c}_4$ (with ${\bf c}_4$ denoting its zeroth-order term)
and prove that
$$
  e^{(\frac{n}{2}+2) \varphi} {\bf R}_4 (e^{2 \varphi} g) (f) = {\bf R}_4 (g) (e^{(\frac{n}{2}-2) \varphi} f)
$$
for all $f \in C^\infty(M)$, $\varphi \in C^\infty(X)$ and all $g$; we recall that in these formulas we suppress the pull-back 
operator
$\iota^*$. It suffices to prove that the conformal variation operator
$$
   f \mapsto  (d/dt)|_0 \left(e^{(\frac{n}{2}+2) t \varphi} {\bf R}_4 (e^{2 t \varphi} g) (e^{-(\frac{n}{2}-2) t\varphi} f) \right)
$$
vanishes for all $\varphi \in C^\infty(X)$ and all $g$. The latter operator is the sum of the conformal variation operator
of the second-order operator ${\bf r}_4$ and the conformal variation
$$
    ({\bf c}_4(g))^\bullet[\varphi] = (d/dt)|_0 (e^{4 t \varphi} {\bf c}_4(e^{2t\varphi} g))
$$
of the zeroth-order term. Only the conformal variation of ${\bf c}_4$ contains normal derivatives of $\varphi$. Lemma \ref{V-1} implies
that these terms are given by
\begin{align*}
    & - 2 |\lo|^2  \partial^2_{0} (\varphi) - 2 (\lo,\Rho) \partial_0(\varphi) + 2 \tfrac{n-3}{n-2} \tr (\lo^3) \partial_0(\varphi)
    + 6 H |\lo|^2 \partial_0 (\varphi) + 2 \tfrac{n-3}{n-2} (\lo,\W) \partial_0(\varphi) \\
    & + \tfrac{8}{n-3} (\lo,\W) \partial_0(\varphi) - \tfrac{2(n-1)^2}{(n-3)(n-2)} (\lo,\W) \partial_0(\varphi) \\
    & + 2 (\lo,\Rho) \partial_0(\varphi) + 2 |\lo|^2 \partial^2_0(\varphi)
    - 6 H |\lo|^2 \partial_0(\varphi) - \tfrac{2(n-3)}{n-2} \tr(\lo^3)  \partial_0(\varphi).
\end{align*}
However, this sum obviously vanishes. Next, the conformal variation operator of ${\bf r}_4$ equals
$$
    -\left(\tfrac{n}{2}-2\right) \delta \left( 4 \tfrac{3n-5}{n-2} \lo^2 d \varphi
    + \tfrac{n^2-12n+16}{2(n-1)(n-2)} |\lo|^2 d\varphi \right).
$$
In order to determine the tangential terms in $({\bf c}_4)^\bullet[\varphi]$, we again apply Lemma \ref{V-1} and the
variation formulas in the following result.

\begin{lem}\label{V-form} In general dimensions, it holds
\begin{equation}\label{var-2}
    (\delta \delta (\lo^2))^\bullet [\varphi] = (n\!-\!2) \delta (\lo^2 d\varphi) - \delta (|\lo|^2 d\varphi)
   + (n\!-\!4) (\delta(\lo^2),d\varphi),
\end{equation}
\begin{equation}
   (\delta(\lo \delta(\lo)))^\bullet[\varphi] = (n\!-\!1) \delta(\lo^2 d\varphi) + (n\!-\!4) (\lo \delta(\lo),d\varphi)
\end{equation}
and
\begin{equation}\label{var-L}
   (\Delta (|\lo|^2))^\bullet[\varphi] = (n\!-\!6) \delta (|\lo|^2 d\varphi) - (n\!-\!4) |\lo|^2 \Delta (\varphi).
\end{equation}
\end{lem}

\begin{proof} We recall that the trace-free part of $L$ satisfies $\hat{\lo} = e^\varphi \lo$. The transformation
laws \eqref{conform-div-forms} (for $\lambda=2$) and \eqref{conform-div-BLF} (for $\lambda=0$) imply
\begin{align*}
    e^{4 \varphi} \hat{\delta} \hat{\delta} (\hat{\lo}^2) & = \delta (e^{2\varphi} \hat{\delta} (\hat{\lo}^2))
   + (n\!-\!4) e^{2\varphi} (d\varphi, \hat{\delta}(\hat{\lo}^2)) \\
   & = \delta \delta (\lo^2) + (n\!-\!2) \delta (\lo^2 d\varphi) - \delta (|\lo|^2 d\varphi) + (n\!-\!4) (d\varphi,\delta (\lo^2))
\end{align*}
(up to non-linear terms) using $\hat{\lo}^2 = \lo^2$. This proves the first relation. Similarly,  \eqref{conform-div-forms}
(for $\lambda=2$) and \eqref{conform-div-BLF} (for $\lambda=-1$) imply
\begin{align*}
    e^{4 \varphi} \hat{\delta} (\hat{\lo} \hat{\delta} (\hat{\lo})) & = \delta (e^{2\varphi} \hat{\lo} \hat{\delta} (\hat{\lo}))
    + (n\!-\!4) e^{2\varphi} (d\varphi, \hat{\lo} \hat{\delta} (\hat{\lo})) \\
    & = \delta (\lo e^\varphi \hat{\delta} (\hat{\lo})) + (n\!-\!4) (d\varphi, \lo e^\varphi \hat{\delta} (\hat{\lo})) \\
    & = \delta (\lo \delta (\lo)) + (n\!-\!1) \delta (\lo^2 d\varphi) + (n\!-\!4) (d\varphi, \lo \delta(\lo))
\end{align*}
(up to non-linear terms). This proves the second relation. Finally, the transformation law \eqref{CTL-Laplace}
(for $\lambda = 2$) implies
\begin{align*}
    e^{4\varphi} \hat{\Delta}(|\hat{\lo}|^2) & = \Delta (|\lo|^2) - 2 \delta (|\lo|^2 d\varphi) + (n\!-\!4) (d(|\lo|^2),d\varphi) \\
    & = \Delta (|\lo|^2) + (n\!-\!6) \delta (|\lo|^2 d\varphi) - (n\!-\!4) |\lo|^2 \Delta(\varphi)
\end{align*}
(up to non-linear terms) using $|\hat{\lo}|^2 = e^{-2\varphi} |\lo|^2$. The proof is complete.
\end{proof}

Lemma \ref{V-form} shows that the tangential terms in $({\bf c}_4)^\bullet[\varphi]$ are given by the product of
$\frac{n}{2}-2$ with
\begin{align*}
    & \tfrac{2(n-1)}{(n-3)(n-2)} \left[ (n\!-\!2) \delta (\lo^2 d\varphi) - \delta (|\lo|^2 d\varphi)
   + (n\!-\!4) (\delta(\lo^2 d\varphi) - (\lo^2,\Hess(\varphi)))\right] \\
   & + \tfrac{4}{n-3} \left[ (n\!-\!1) \delta(\lo^2d\varphi) + (n\!-\!4) (\lo \delta(\lo),d\varphi)\right] \\
   & + \tfrac{3n-4}{2(n-1)(n-2)} \left[(n\!-\!6) \delta (|\lo|^2 d\varphi) - (n\!-\!4) |\lo|^2 \Delta (\varphi) \right] \\
   & + |\lo|^2 \Delta (\varphi) - \frac{4n}{n\!-\!1} (\lo \delta(\lo),d\varphi) + 2 H(\lo,\Hess(\varphi)) - \delta(|\lo|^2 d\varphi)
   + 4 \delta(\lo^2 d\varphi) \\
   & - \tfrac{4}{n-3} \left[ 2 \delta(\lo^2 d\varphi) - 2 (\lo^2,\Hess(\varphi)) - 2 \tfrac{n-2}{n-1} (\lo \delta(\lo),d\varphi)
   + |\lo|^2 \Delta(\varphi) - \delta (|\lo|^2 d\varphi)\right] \\
   & + \tfrac{2(n^2-9n+12)}{(n-3)(n-2)} (\lo^2,\Hess(\varphi)) + \tfrac{n^3-5n^2+18n-20}{2(n-3)(n-2)(n-1)} |\lo|^2 \Delta(\varphi)
  - 2 H(\lo,\Hess(\varphi)).
\end{align*}
This sum vanishes. Summarizing the above results shows that the conformal variation of ${\bf R}_4$ vanishes.


\subsection{Extrinsic conformal invariants of hypersurfaces}\label{ECI}

The scalar invariants
\begin{itemize}
\item $\Iv_1 = |\lo|^4$, $\Iv_2 = \tr(\lo^4)$,
\item $\Iv_3 = |\overline{W}|^2$, $\Iv_4 = |\W|^2$,
\item $\bar{\Iv}_5 = \lo^{ij} \lo^{kl} \overline{W}_{iklj}$, $\Iv_6 = (\lo^2,\W)$, $\Iv_7 = |\overline{W}_0|^2$
\end{itemize}
of an embedding $M^4 \hookrightarrow X^{5}$ are obvious local conformal invariants of weight $-4$. Here we set 
$$
   |\W|^2 \st \W_{ij} \W^{ij} \quad \mbox{with} \quad \W_{ij} \st \overline{W}_{0ij0}
$$
and 
$$
   |\overline{W}|^2 \st \overline{W}_{ijkl} \overline{W}^{ijkl} \quad \mbox{and} 
   \quad  |\overline{W}_0|^2 \st \overline{W}_{ijk0} \overline{W}^{ijk0};
$$ 
in these definitions, all indices $i,j,k,l$ are tangential. The above invariants are defined in terms of the trace-free part $\lo$ 
of $L$ and the Weyl tensor $\overline{W}$ of the background metric. 

Note that the local conformal invariant $|W|^2 = W_{ijkl} W^{ijkl}$ of weight $-4$ is a linear combination of 
$\Iv_3$ and the other invariants (Lemma \ref{Weyl-ext}). Likewise, the invariant $\Iv_5 = \lo^{ij} \lo^{kl} W_{iklj}$ 
is a linear combination of the other invariants (see \eqref{II5}). The above conformal invariants are invariant 
under a change of the orientation of the normal vector. 

In addition, we have the {\em non-trivial} local conformal invariants $\mathcal{J}_1$ and $\mathcal{J}_2$
(defined in \eqref{CI-1} and \eqref{CI-2}). Note that the definitions of $\Jv_1$ and $\Jv_2$ contain the 
respective normal derivative terms $\lo^{ij} \bar{\nabla}_0 (\overline{W})_{0ij0}$ and 
$(\lo,\bar{\nabla}_0(\bar{\Rho}))$. Corollary \ref{J12-LC} shows that $\Jv_1 - 2 \Jv_2$ again is a linear 
combination of the above invariants. Note also that $\Jv_1$ and $\Jv_2$ both contain non-trivial divergence terms - 
they are conformally invariant only with these divergence terms.

Next, we have the local conformal invariant $\Jv_3$ (see \eqref{CI-3}). In a forthcoming paper, Astaneh and Solodukhin 
will prove that the integral of $\Jv_3$ is a linear combination of the integrals of $\bar{\Iv}_5$, $\Iv_6$ 
and $\Iv_7$ (see \eqref{J3-AS}).\footnote{We are grateful to S. Solodukhin for informing us about this
result \cite{AS2}.} 

Finally, we note that $\lo^{ij} \nabla^k \overline{W}_{kij0}$ is a local conformal invariant in dimension $n=4$ 
(see the proof of \cite[Lemma 6.27]{JO2}). Lemma \ref{div-term} shows that its integral 
reduces to the conformal invariant $\int_M |\overline{W}_0|^2 dvol_h$. However, the divergence term 
$\delta (\lo^{ij} \overline{W}_{\cdot ij0})$ itself is a conformal invariant of weight $-4$ (see the comment 
after Remark \ref{action-local}).

The latter observation is related to the local conformal invariant $\Jv_5 = (\lo,S)$ (defined in 
Remark \ref{invariant-E}).  In fact, we find
\begin{align*}
   (\lo,S) & = \lo^{ij} (\bar{\nabla}_0 (\overline{W})_{0ij0} - \bar{C}_{ij0} - \bar{C}_{ji0}) + 4 H(\lo,W) \\
   & = - \lo^{ij} \bar{\nabla}^k (\overline{W})_{kij0} + 4 H (\lo,\W) 
\end{align*}
using $\bar{C}_{ijk} = \frac{1}{2} \bar{\nabla}^l(\overline{W})_{lijk}$. But
$$
   \lo^{ij} \bar{\nabla}^k (\overline{W})_{kij0} = \lo^{ij} \nabla^k \overline{W}_{kij0} - (\lo^2,\W) 
   + 4 H (\lo,\W) -   \lo^{ij} \lo^{kl} \overline{W}_{kijl} 
$$
in dimension $n=4$ (see \eqref{n-reduce} in the proof of Proposition \ref{B-JJD}). Hence
\begin{equation}\label{J5-expand}
   \Jv_5 = (\lo^2,\W) + \lo^{ij} \lo^{kl} \overline{W}_{kijl} - \lo^{ij} \nabla^k \overline{W}_{kij0}. 
\end{equation}
Note that this relation again implies the conformal invariance of $\lo^{ij} \nabla^k \overline{W}_{kij0}$. 

It is an open problem to classify all local conformal invariants of weight $-4$ of a hypersurface 
$M^4 \hookrightarrow X^5$. An easier problem is classifying all global conformal integrals attached to 
an embedding $M^4 \hookrightarrow X^5$.

\subsection{Decompositions of conformal anomalies}\label{deco-general}

The local conformal invariants in Section \ref{ECI} are also of interest in other parts of geometric analysis 
and theoretical physics.

Let $(X^{n+1},g)$ be a compact odd-dimensional manifold with smooth even-dimensional boundary $M^n$.
We consider the boundary value problem for the Yamabe operator $P_2(g)$ on $X$ with Dirichlet or Robin boundary
conditions. The constant term $a_{n+1}$ in the small-time asymptotic expansion of the trace of the heat kernel of
this boundary value problem is a global conformal invariant.\footnote{The coefficient $a_{n+1}$ is also called the critical
heat kernel coefficient.} It is given by an integral of curvature invariants of the embedding $M \hookrightarrow X$. This 
result is a consequence of the conformal index property of the critical heat kernel coefficient. The heat kernel coefficient
$a_{n+1}$ may be regarded as the integrated conformal anomaly of the functional determinant of the boundary 
value problem \cite[Section 2]{BG-BVP}. 

It is expected that $a_{n+1}$ is a linear combination of the Euler characteristic $\chi(M)$ and integrals of 
local conformal invariants of the embedding $M \hookrightarrow X$. However, in general, the structure of the local conformal 
invariants in that decomposition is unknown.

For $n=2$, the heat kernel coefficient $a_3$ is a linear combination of the Euler-characteristic $\chi(M)$ 
and the integral of $|\lo|^2$ (see \cite[Section 2]{CQ}). For $n=4$, the heat kernel coefficient $a_5$ has been 
determined in \cite{BGKV}, and its decomposition has been studied recently in \cite{AS}. It involves the invariants
$\Iv_j$ listed in Section \ref{ECI} and $\Jv_1$. 

The conformal invariance of the functional $\int_{M^4} \Jv_1 dvol_h$ is one of the main results 
of \cite[(29)]{AS}.\footnote{$\Jv_1$ is the functional $I_8$ in the notation of \cite{AS}.} Supported by the decomposition 
of $a_5$, Astaneh and Solodukhin state that the integral of the Pfaffian and the integrals of the conformally invariant 
curvature quantities $\Iv_j$, together with the integral of $\Jv_1$, form a basis of all conformally invariant integrals
associated to $M^4 \hookrightarrow X^5$. We recall that the existence of the local invariants $\Jv_2$ (or $\mathcal{C}$), 
$\Jv_3$, $\Jv_4$, and $\Jv_5$ does not contradict that completeness statement since their integrals are linear 
combinations of the other invariants. 

For odd $n$, the situation is different. Then $a_{n+1}$ is a conformally invariant sum of an integral on $X$ 
and a boundary integral on $M$. Its decomposition is expected to have the form
\begin{equation}\label{DS-odd}
    a \chi(X) + \sum_j c_j \int_X \bar{\Iv}_j dvol_g + \sum_j  b_j \int_M \Jv_j dvol_h
\end{equation}
with  local conformal invariants $\bar{\Iv}_j$ of $X$ (only depending on $\overline{W}$) and 
local conformal invariants $\Jv_j$ of the embedding $M \hookrightarrow X$ (only depending on 
$\overline{W}$ and $\lo$). 

For $n=3$, the heat kernel coefficient $a_4$ is a linear combination of the Euler characteristic $\chi(X)$, 
$\int_X |\overline{W}|^2 dvol_g$ and the boundary integrals 
$$
   \int_M(\lo,\W) dvol_h \quad \mbox{and} \quad \int_M \tr(\lo^3) dvol_h
$$
(\cite[Theorem 3.7]{BG-BVP}). In higher dimensions, the situation is much less understood. 

Of course, there are similar problems for more general conformally invariant boundary value problems.

The above decompositions may be regarded as analogs of the Deser-Schwimmer classification of global conformal 
anomalies of CFTs on closed manifolds established by Alexakis (see \cite{alex} and its references). 

In the framework of CFTs on manifolds $X^{n+1}$ with boundary $M^n$ (BCFT), it is a key problem to 
classify the (integrated) conformal anomalies. These quantities are expected to have analogous decompositions. 
More precisely, one expects that, for even $n$, they have the form 
\begin{equation}\label{DS-even}
   a \chi(M) + \sum_j b_j \int_M \Jv_j dvol_h
\end{equation}
with local conformal invariants $\Jv_j$ of the embedding $M \hookrightarrow X$ (consisting of extrinsic 
and intrinsic invariants). Similarly, for odd $n$, they should decompose as in \eqref{DS-odd}.

Moreover, these decompositions should follow from corresponding decompositions of the anomalies themselves.
In these decompositions, the Euler form $E_n$ of $M$ (for even $n$) is responsible for the Euler characteristic of $M$. 
Likewise, for odd $n$, the Euler form $E_{n+1}$ of $X$, together with a boundary term $E_{n+1}^\partial$ 
on $M$ - according to the Chern-Gauss-Bonnet theorem on manifolds with boundary (see \cite[Chapter 4]{Gilkey-book}) - 
are responsible for the Euler characteristic $\chi(X)$ contributing to \eqref{DS-odd}. 

For $n=2$ and $n=3$, the boundary terms of the respective anomalies decompose as $a \J + b |\lo|^2$ and 
$a E_4^\partial + b_1 \tr(\lo^3) + b_2 (\lo,\W)$, up to divergence terms. Note that ${\bf Q}_2$ and
${\bf Q}_3$ (see Proposition \ref{P23}) have the same structure as the boundary terms in these decompositions. 
Theorem \ref{alex} implies a decomposition of $\int_M {\bf Q}_4 dvol_h$ which is of the form \eqref{DS-even}.
The relevant local invariants are listed in Section \ref{ECI} and include $\Jv_1$. 

Generalizing \cite{AS}, the authors of \cite{CHBRS} determined the most general form of the boundary terms 
in the conformal anomaly of a CFT on a manifold $X$ of dimension $d \ge 5$ with a boundary (or defect) 
$M$ of dimension $4$. The identity \cite[(3.1)]{CHBRS} gives the general form of the anomaly of such a BCFT.
Apart from the Pfaffian of $M$, it contains {\em two non-trivial} invariants $\Jv_{(1)}$, $\Jv_{(2)}$, and a series 
of functionals of $\overline{W}$ and $\lo$. In the case $d=5$, this confirms the decomposition \eqref{DS-even} 
with the local invariants $\Iv_j$ and $\Jv_1$ listed in Section \ref{ECI}. In particular, in the codimension-one case, 
\cite[(6.2)]{CHBRS} states that the decomposition of the conformally invariant Graham-Reichert functional 
reduces to\footnote{The Weyl tensor in \cite{CHBRS} has the opposite sign.}
$$
    8 \mathcal{E}_{GR} = \int_M \left( \Pf_4 + \frac{1}{2} \Jv_{(1)} - \frac{3}{8} \Iv_1 + \frac{3}{4} \Iv_2 
    - \frac{1}{8} \Iv_3 + \frac{1}{2} \Iv_4 - \bar{\Iv}_5 + \frac{1}{2} \Iv_6 + \frac{1}{4} \Iv_7 \right) dvol_h.
$$
The Gauss equations imply the relation $\Jv_{(1)} = \Jv_1 + \frac{1}{3} \Iv_1 - \Iv_2 - \Iv_6$ (up to a divergence)
and we get 
\begin{equation}\label{CH-E}
    8 \mathcal{E}_{GR} = \int_M  \left( \Pf_4 + \frac{1}{2} \Jv_1 - \frac{5}{24} \Iv_1 + \frac{1}{4} \Iv_2 
    - \frac{1}{8} \Iv_3 + \frac{1}{2} \Iv_4 - \bar{\Iv}_5 + \frac{1}{4} \Iv_7 \right) dvol_h.
\end{equation}
Note that this formula (slightly) differs from \eqref{EGR-final}.

The anomalies listed in \cite[Section 3.3]{CHBRS} contain three more local invariants, one of
which is an extrinsic analog of the $4$-form $\tr (W \wedge W)$ (first Pontrjagin form).\footnote{Here $W$ 
is regarded as an $\End(TM)$-valued $2$-form.} In contrast to the invariants discussed above, these change 
signs under a simultaneous change of the orientations on $X$ and $M$.

From the perspective of the AdS/CFT duality, the conformal (quantum) anomalies of determinants (conformal index) of 
conformally covariant operators on a manifold $(X,g)$ (of even dimension) appear as duals of the anomaly of the 
renormalized volume of an associated Poincar\'e-Einstein metric with conformal infinity $[g]$ \cite{GZ}. The latter anomaly 
is proportional to the total integral $\int_X Q_n(g) dvol_g$ of Branson's critical $Q$-curvature \cite{GZ}. This result 
extends to a relation between the total integral $\int_X{\bf Q}_n(h) dvol_h$ and the anomaly of the renormalized volume 
of the singular Yamabe metric $\sigma^{-2} g$ on $X$ \cite{G-SY, GW-RV, JO1}.

There is an analog of the AdS/CFT duality for CFTs, which relates anomalies of BCFTs to geometric anomalies
of dual theories. This duality naturally involves the study of conformal invariants of submanifolds (see \cite{RT} and its references).
In particular, on the geometric side, this leads to the study of the conformal anomaly of the renormalized volume of a 
minimal hypersurface with boundary $M$ in a Poincar\'e-Einstein background with conformal infinity $X$ as initiated in \cite{GWi}.
Graham and Reichert \cite{GR} analyzed these conformal anomalies. In particular, they derived an explicit formula for this
global conformal invariant of an embedding of $M^4 \hookrightarrow X^n$ ($n \ge 5$). Parallel work
\cite{Zhang} led to equivalent formulas. Finally, from the perspective of BCFTs, these results were discussed in \cite{CHBRS}. 

For more details on anomalies, we refer to \cite{FV-book, Fu, Sol, AS, HH}. For further generalizations and a unified 
discussion of geometric anomalies, we refer to \cite{AGW}.

\subsection{The invariants $\text{Wm}$ and $(D (\lo),\lo)$}\label{invariant}

In \cite{BGW-1}, the authors derived formulas for ${\bf P}_4$ and ${\bf Q}_4$ in general dimensions using
conformal tractor calculus. A central  role plays the local conformal invariant ${\text Wm}$. In the present
section, we discuss this invariant. A full comparison of our results with the corresponding results in \cite{BGW-1}
in the critical dimension $n=4$ will be given in the Section \ref{BGW-1}.

The local conformal invariant ${\text Wm}$ is the sum of the term $(\Delta (\lo),\lo)$ and some curvature terms.
This is an interesting result on its own. In order to better understand the conformal invariant ${\text Wm}$, we next
describe this invariant in terms of obviously conformally invariant terms. The original proof of the conformal
invariance of ${\text Wm}$ depends on heavy tractor calculus machinery.

From \cite{BGW-1}, we recall the definition
\begin{equation}\label{Wm-def}
   \text{Wm} \st \frac{1}{2} (\lo,\Delta (\lo)) + \frac{4}{3} \delta (\lo \delta(\lo)) + \frac{3}{2} \Delta (|\lo|^2)
   - 6 (\lo,\bar{C}_0) + 4 (\lo^2,\Rho) - \frac{7}{2} \J |\lo|^2 + 6 H (\lo,\W),
\end{equation}
where $\bar{C}$ is the Cotton tensor of the background metric and $(\bar{C}_0)_{ij} = \bar{C}_{ij0}$. This 
step uses the fact that the formula $\mathring{\JF} = \frac{1}{2} \lo^2 + \frac{1}{2} \W - \frac{1}{8} |\lo|^2 h$ 
(see \eqref{TF-F}) implies the relation $- 6 H \tr(\lo^3) + 12 H (\lo,\mathring{\JF}) = 6 H (\lo,\W)$.

\begin{prop}\label{WM-deco} Let $n=4$. Then
\begin{equation}\label{Wm-expansion}
    \text{Wm} = \frac{1}{2} (D (\lo),\lo) - 3 \Jv_1 + 3 \bar{\Iv}_5 + 3 \Iv_6 - \frac{3}{2} \Iv_7
   - 6 \delta (\lo,\overline{W}_0),
\end{equation}
where
$$
     D (t)_{ij} \st \Delta (t)_{ij} - 2 (\Rho \circ t + t \circ \Rho)_{ij} - \J t_{ij}
     - \frac{2}{3} (\nabla_i \delta(t)_j + \nabla_j \delta(t)_i) + h_{ij} (\Rho,t) + \frac{1}{3} h_{ij} \delta \delta (t)
$$
is a conformally covariant operator $S_0^2(M) \to S_0^2(M)$ on trace-free symmetric $2$-tensors:
$e^\varphi \hat{D} (t) = D (e^{-\varphi} t)$.
\end{prop}

Some comments are in order.

We recall that the scalar product $(\lo,\overline{W}_0) \in \Omega^1(M)$ is defined as $\lo^{ij} \overline{W}_{\cdot ij0}$.

All terms on the right-hand side of  \eqref{Wm-expansion} are local conformal invariants of weight $-4$. In 
particular, we recall that the local invariant
\begin{equation*}
   \Jv_1 = \lo^{ij} \bar{\nabla}_0(\overline{W})_{0ij0} + 2 H (\lo,\W) + \frac{2}{9} |\delta(\lo)|^2
   - 2 (\lo^2,\Rho) + \J |\lo|^2 - \delta \delta (\lo^2)
\end{equation*}
(see \eqref{CI-1}) satisfies $e^{4 \iota^*(\varphi)} \hat{\Jv}_1 = \Jv_1$ (Lemma  \ref{new-invariants}). Moreover,
the divergence $\delta (\lo,\overline{W}_{0})$ is a conformal invariant of weight $-4$ (see the comment after Remark \ref{action-local}). Therefore, Proposition 
\ref{WM-deco} shows that ${\text Wm}$ is a local conformal invariant of weight $-4$. This reproves a part of
\cite[Theorem 1.2]{BGW-1} (see also \cite[Theorem 1.5]{BGW-2}).

In general dimensions, the operator
\begin{align}\label{D-g}
   D(t)_{ij} & \st \Delta(t)_{ij} - 2 (\Rho \circ t + t \circ \Rho)_{ij} -\left (\frac{n}{2}-1\right) \J t_{ij}
   - \frac{4}{n+2} (\nabla_i \delta(t)_j + \nabla_j \delta(t)_i) \notag \\
   & + \frac{4}{n} h_{ij} (\Rho,t) + \frac{8}{n(n+2)}
    h_{ij} \delta \delta (t)
\end{align}
maps trace-free symmetric $2$-tensors to trace-free symmetric $2$-tensors. $D$ is conformally covariant in the
sense that
$$
   e^{(\frac{n}{2}-1)\varphi} \hat{D} (t) = D (e^{(\frac{n}{2}-3)\varphi} t).
$$
The operator $D$ in Proposition \ref{WM-deco} is its special case in dimension $n=4$.

A conformally covariant generalization of the operator $D$ to an operator on trace-free symmetric $p$-tensors was
discovered in \cite{Wunsch}. It satisfies
$$
   e^{(\frac{n}{2}+1-p)\varphi} \hat{D} (t) = D (e^{(\frac{n}{2}-1-p) \varphi} t).
$$
We refer to \cite[Chapter 2]{Jenne} for an ambient metric derivation of it. Matsumoto \cite{Matsumoto}
used an ambient metric approach to define analogs of the GJMS-operators acting on trace-free symmetric $2$-tensors.
On divergence-free and trace-free symmetric $2$-tensors, the second-order operator $P_2$ in this sequence at an Einstein
metric acts like a linear combination of $D$ in \eqref{D-g} and the obviously conformally invariant operator $t_{ij} \mapsto
W_{iklj} t^{kl}$.

The formula \eqref{D-g} and the related action functional
$$
    \int_M (D (t),t) dvol_h
$$
also have been derived in \cite{AHR} (correcting \cite{EO}). Since $D$ acts on trace-free symmetric tensors, it is not possible
to read off the full operator $D$ from the associated action.

In \cite[Section 5]{br-NL}, Branson classified all conformally covariant second-order operators acting on irreducible
tensor bundles. Moreover, explicit formulas are given for operators on $p$-forms and trace-free symmetric $p$-tensors.
The operator in \eqref{D-g} should be a particular  case of \cite[(5.7)]{br-NL}.\footnote{The curvature terms
in this formula differ from the above result.} Branson also showed that the operator $D$ is unique, up to a constant
multiple of the above action of the Weyl tensor.

\begin{proof} The definitions of $D$ and $\Jv_1$ yield
\begin{equation}\label{D-action}
   \frac{1}{2} (D (\lo),\lo) = \frac{1}{2} (\Delta (\lo),\lo) - 2 (\Rho,\lo^2) - \frac{1}{2} \J |\lo|^2
   - \frac{2}{3} \lo^{ij} \nabla_i \delta(\lo)_j
\end{equation}
and
$$
    3 \Jv_1 = 3 \lo^{ij} \bar{\nabla}_0(\overline{W})_{0ij0} + 6 H (\lo,\W) + \frac{2}{3} |\delta(\lo)|^2
   - 6 (\lo^2,\Rho) + 3 \J |\lo|^2 - 3 \delta \delta (\lo^2).
$$
Hence
\begin{align*}
   \frac{1}{2} (D (\lo),\lo) -  3 \Jv_1 & =  \frac{1}{2} (\Delta (\lo),\lo) + 4 (\lo^2,\Rho) - \frac{7}{2} \J |\lo|^2
   - \frac{2}{3} \lo^{ij} \nabla_i \delta(\lo)_j - \frac{2}{3} |\delta(\lo)|^2 + 3 \delta \delta (\lo^2) \\
   & - 6 H (\lo,\W) - 3 \lo^{ij} \bar{\nabla}_0(\overline{W})_{0ij0}.
\end{align*}
It follows that the difference ${\text Wm} -  \frac{1}{2} (D (\lo),\lo) + 3 \Jv_1$ equals
\begin{align*}
     & \frac{3}{2} \Delta (|\lo|^2) + \frac{4}{3} \delta (\lo \delta(\lo)) - 3 \delta \delta (\lo^2)
     + \frac{2}{3} \lo^{ij} \nabla_i \delta(\lo)_j + \frac{2}{3} |\delta(\lo)|^2 \\
     & - 6 (\lo,\bar{C}_0) + 12 H (\lo,\W) + 3 \lo^{ij} \bar{\nabla}_0(\overline{W})_{0ij0}.
\end{align*}
Now the trace-free Codazzi-Mainardi equation implies the relation
$$
    3 \Delta (|\lo|^2) - 6 \delta \delta (\lo^2) + 4 \delta (\lo \delta(\lo))
    = - 6 \delta (\lo^{ij} \overline{W}_{\cdot ij0}) = - 6 \delta (\lo,\overline{W}_0)
$$
(see \eqref{WN}). Thus, the latter sum equals
\begin{align*}
    & - \frac{2}{3} \delta (\lo \delta(\lo)) + \frac{2}{3} \lo^{ij} \nabla_i \delta(\lo)_j
    + \frac{2}{3} |\delta(\lo)|^2 \\
    & + 12 H (\lo,\W) - 3 \delta (\lo,\overline{W}_0) - 6 (\lo,\bar{C}_0) + 3 \lo^{ij} \bar{\nabla}_0(\overline{W})_{0ij0} \\
    & = 12 H (\lo,\W) - 3 \delta (\lo,\overline{W}_0) - 6 (\lo,\bar{C}_0) + 3 \lo^{ij} \bar{\nabla}_0(\overline{W})_{0ij0}.
\end{align*}
Now we note that
\begin{align*}
    -2 (\lo,\bar{C}_0) + \lo^{ij} \bar{\nabla}_0(\overline{W})_{0ij0} 
    & = - \lo^{ij} \bar{\nabla}^k (\overline{W})_{kij0} \\
    & = - \lo^{ij} \nabla^k \overline{W}_{kij0} + (\lo^2,\W) - 4 H(\lo,\W) + \lo^{ij} \lo^{kl} \overline{W}_{kijl}
\end{align*}
using $2 (\bar{C}_0)_{ij} = 2 \bar{C}_{ij0} = \bar{\nabla}^a (\overline{W})_{aij0}$ and \eqref{n-reduce}. 
Hence
\begin{align*}
   & {\text Wm} -  \frac{1}{2} (D (\lo),\lo) + 3 \Jv_1 = - 3 \delta (\lo,\overline{W}_0) 
  - 3 \lo^{ij} \nabla^k \overline{W}_{kij0} + 3 (\lo^2,\W)
  + 3 \lo^{ij} \lo^{kl} \overline{W}_{kijl}.
\end{align*}
Now we apply Lemma \ref{div-term}. This completes the proof.
\end{proof}

The decomposition \eqref{Wm-expansion} contains the local conformal invariant $(D(\lo),\lo)$. We next
express the resulting conformally invariant action functional
$$
   \int_{M^4} (D(\lo),\lo) dvol_h
$$
in terms of the basic conformal invariants listed in Section \ref{ECI}. It seems remarkable that the result does
not involve derivatives of $\lo$.\footnote{The first relation in Proposition \ref{action} corrects a typo in formula (20) in \cite{AS}.}

\begin{prop}\label{action} For a closed four-manifold $M$, it holds
$$
   \int_M (D(\lo),\lo) dvol_h = \int_M \left( -\frac{7}{6} \Iv_1 + 2 \Iv_2 -\bar{\Iv}_5 + \Iv_6 
    - \frac{1}{2} \Iv_7 \right) dvol_h
   = \int_M \left(-I_5 - \frac{1}{2} I_7 \right) dvol_h.
$$
\end{prop}

\begin{proof} An identity of Simons \cite{HP}
states that for any hypersurface $M^n \hookrightarrow X^{n+1}$ with the second fundamental form $L$, it holds
\begin{align*}
   \Delta (L)_{ij} & = n \Hess_{ij}(H) + n H L^2_{ij} - L_{ij} |L|^2\notag \\
   & +  L^s_j \bar{R}_{ikks} + L_i^s \bar{R}_{jkks} - 2 L^{rs}  \bar{R}_{rijs}
   + n H \bar{R}_{0ij0} - L_{ij} \overline{\Ric}_{00} + \bar{\nabla}^k (\bar{R})_{ikj0} + \bar{\nabla}_i (\bar{R})_{jkk0}.
\end{align*}
In the following, it will be convenient to restate that identity in terms of covariant derivatives of the hypersurface only. In fact,
the Gauss identity for the curvature tensor shows that the above identity is equivalent to
\begin{align}\label{Simons-1}
    \Delta (L)_{ij} & = n \Hess_{ij}(H) + n H L^2_{ij} - L_{ij} |L|^2\notag \\
    & + L_j^s \bar{R}_{ik}{}^{k}{}_{s} -  L^{rs}  \bar{R}_{rijs}
   + \nabla^k \bar{R}_{ikj0} + \nabla_i \bar{R}_{jk}{}^{k}{}_{0}.
\end{align}
Now let $n=4$. The identity $(\Delta(L),L) = (\Delta(\lo),\lo) + 4 H \Delta(H)$ and the relation \eqref{Simons-1} imply
\begin{align*}
    (\Delta(\lo),\lo) & = 4 (\lo,\Hess(H)) + 4 H \tr(L^3) - |L|^4 \\
    & + (L^2)^{ij} \bar{R}_{ik}{}^{k}{}_{j} - L^{ij} L^{rs} \bar{R}_{rijs} + L^{ij} \nabla^k \bar{R}_{ikj0}
    + L^{ij} \nabla_i(\overline{\Ric}_0)_j.
\end{align*}
We integrate and apply partial integration. Hence
\begin{align}\label{DLL-ex}
   & \int_M (D(\lo),\lo) dvol_h \notag \\
   & = \int_M \left( (\Delta(\lo),\lo) - 4 (\lo^2,\Rho) - \frac{4}{3} (\lo,\nabla \delta(\lo)) - \J |\lo|^2 \right) dvol_h
   \qquad \qquad \mbox{(by definition)} \notag \\
   & = \int_M \left( - 4 (\delta(\lo),dH) + \frac{4}{3} (\delta(\lo),\delta(\lo)) - 4 (\lo^2,\Rho)  - \J |\lo|^2
   + 4 H \tr(L^3) - |L|^4 \right) dvol_h  \notag \\
   & + \int_M \left( (L^2)^{ij} \bar{R}_{ik}{}^{k}{}_{j} - L^{ij} L^{rs} \bar{R}_{rijs} + L^{ij} \nabla^k \bar{R}_{ikj0}
   - (\delta(L),\overline{\Ric}_0) \right) dvol_h.
\end{align}
The Codazzi-Mainardi equation $\delta(\lo) = 3 dH + 3 \bar{\Rho}_0$ shows that the sum of the latter two integrals 
equals
\begin{align*}
   \int_M 4 (\delta(\lo), \bar{\Rho}_0) - (\delta(L), \overline{\Ric}_0) - 4 (\lo^2,\Rho)  - \J |\lo|^2 + 4 H \tr(L^4) - |L|^4
   + \mbox{\text{curvature tensor terms}}.
\end{align*}
Thus, using
$$
   4 H \tr(L^3) - |L|^4 = 4 H \tr(\lo^3) + 4 H^2 |\lo|^2 - |\lo|^4,
$$
this sum can be written as the sum of the integral
\begin{equation}\label{HE-1}
   \int_M \left( - 4 (\lo^2,\Rho)  - \J |\lo|^2 + 4 H \tr(\lo^3) + 4 H^2 |\lo|^2 - |\lo|^4 \right)  dvol_h
\end{equation}
and the integral
\begin{equation}\label{HE-2}
    \int_M  \left((\delta(\lo), \bar{\Rho}_0) - (dH, \overline{\Ric}_0)  +
    (L^2)^{ij} \bar{R}_{ik}{}^{k}{}_{j} - L^{ij} L^{rs} \bar{R}_{rijs} + L^{ij} \nabla^k \bar{R}_{ikj0} \right) dvol_h
\end{equation}
of curvature terms. Now the Fialkow equation
$$
   \Rho = \bar{\Rho} + H \lo + \frac{1}{2} H^2 h - \frac{1}{2} \lo^2 + \frac{1}{12} |\lo|^2 h - \frac{1}{2} \W
$$
(see \eqref{Fial}) and the Gauss equation $\J = \bar{\J} - \bar{\Rho}_{00} - \frac{1}{6} |\lo|^2 + 2 H^2$
(see \eqref{Gauss-scalar}) imply that the integral \eqref{HE-1} equals
\begin{equation}\label{HE-3}
    \int_M -4 (\lo^2,\bar{\Rho}) -\bar{\J}|\lo|^2 + \bar{\Rho}_{00} |\lo|^2 + 2 (\lo^2,\W) +2 \tr(\lo^4) - \frac{7}{6} |\lo|^4.
\end{equation}
It remains to calculate the integral \eqref{HE-2}. First, the decomposition $\bar{R} = \overline{W} - \bar{\Rho} \owedge g$ implies     
$$
     L^{ij} L^{rs} \bar{R}_{rijs}
    =  L^{ij} L^{rs} \overline{W}_{rijs} - 2 (\lo^2,\bar{\Rho}) - 2 H (\lo,\W) + 4 H(\lo,\bar{\Rho}) + 6 H^2 (\bar{\J} - \bar{\Rho}_{00}).
$$
Second, we have
\begin{align*}
     (L^2)^{ij} \bar{R}_{ik}{}^{k}{}_{j} & = (L^2, \overline{\Ric} - \bar{\G}) \\
     & = (L^2, 3 \bar{\Rho}+ \bar{\J} h) - (L^2,\bar{\Rho}+\bar{\Rho}_{00} h + \W).
\end{align*}
We combine these results and simplify. Then
$$
    (L^2)^{ij} \bar{R}_{ik}{}^{k}{}_{j} - L^{ij} L^{rs} \bar{R}_{rijs}
   = - L^{ij} L^{rs} \overline{W}_{rijs} + 4 (\lo^2,\bar{\Rho}) + \bar{\J} |\lo|^2 - \bar{\Rho}_{00} |\lo|^2 - (\lo^2,\W).
$$
Finally, we calculate
\begin{align*}
    L^{ij} \nabla^k \bar{R}_{kij0} & = \lo^{ij} \nabla^k \bar{R}_{kij0} + H \nabla^k \overline{\Ric}_{k0} \\
   & = \lo^{ij} \nabla^k \overline{W}_{kij0} - \lo^{ij} \nabla^k (\bar{\Rho}_{0})_i h_{kj} 
   + H \nabla^k \overline{\Ric}_{k0}.
\end{align*}
Hence
\begin{equation}\label{curv-term}
    \int_M L^{ij} \nabla^k \bar{R}_{kij0} dvol_h
   = \int_M \left( \lo^{ij} \nabla^k \overline{W}_{kij0} + (\delta(\lo), \bar{\Rho}_0)
    - (dH,\overline{\Ric}_0)\right) dvol_h.
\end{equation}
Now summarizing these results shows that
$$
   \int_M (D(\lo),\lo) dvol_h = \int_M  \left( (\lo^2,\W) +2 \tr(\lo^4) - \frac{7}{6} |\lo|^4 
   - L^{ij} L^{rs} \overline{W}_{rijs} - \lo^{ij} \nabla^k \overline{W}_{kij0}\right) dvol_h.
$$
But
$$
    \int_M \lo^{ij} \nabla^k \overline{W}_{kij0} dvol_h = \frac{1}{2} \int_M |\overline{W}_0|^2 dvol_h
$$
by Lemma \ref{div-term}. This proves the first relation. In general dimensions, the Gauss equation
\eqref{Gauss-Weyl} for the Weyl tensor implies
\begin{align}\label{I5-bar}
   \bar{\Iv}_5 & = \Iv_5 + \tr (\lo^4) - |\lo|^4 + 2 (\lo^2,\JF) \notag \\
   & = \Iv_5 + \tr (\lo^4)  - |\lo|^4 + \frac{2}{n-2} \left(\tr(\lo^4) - \frac{1}{2(n-1)} |\lo|^4 + (\lo^2,\W)\right).
\end{align}
In particular, for $n=4$, we obtain
\begin{equation}\label{II5}
    \bar{\Iv}_5 = \Iv_5 - \frac{7}{6} \Iv_1 + 2 \Iv_2 + \Iv_6.
\end{equation}
Combining this identity with the first relation proves the second relation.
\end{proof}

The following result is a local version of the second relation in Proposition \ref{action}.

\begin{rem}\label{action-local} Let $n=4$. Then
$
    (D(\lo),\lo) = -I_5 - \frac{1}{2} I_7 - \delta (\lo,\overline{W}_{0}).
$
\end{rem}

\begin{proof} The arguments in the proof of Proposition \ref{action} lead to the additional divergence terms
\begin{align*}
     & 4 \delta (\lo dH) - \frac{4}{3} \delta (\lo \delta(\lo)) + \delta (L \overline{\Ric}_0)
     & \mbox{(in \eqref{DLL-ex})} \\
     & + \delta (\lo \bar{\Rho}_0) - \delta (H \overline{\Ric}_0) & \mbox{(by \eqref{curv-term})} \\
     & - \delta (\lo^{ij} \overline{W}_{\cdot ij 0}).
\end{align*}
By the Codazzi-Mainardi equation, all terms except the last one cancel.
\end{proof}

As a byproduct, Remark \ref{action-local} implies that $\delta (\lo,\overline{W}_{0})$ is a local conformal invariant of weight $-4$.

We finish this section with the 

\begin{proof}[Proof of Lemma \ref{Guven-Wm}] Proposition \ref{WM-deco} shows that for a flat background
$$
    \int_M {\text Wm} dvol_h = \frac{1}{2} \int_M  (D(\lo),\lo) dvol_h - 3 \int_M \Jv_1 dvol_h. 
$$
But $ \frac{1}{2} \int_M(D(\lo),\lo) dvol_h = \int_M (-\frac{7}{6} \Iv_1 + 2 \Iv_2) dvol_h$ by Proposition \ref{action}. 
Hence
$$
    \int_M {\text Wm} dvol_h = \int_M \left( \frac{5}{12} \Iv_1 - 2 \Iv_2 
   - 6 |dH|^2 + 6 H \tr(\lo^3) - 3 H^2 |\lo|^2 \right) dvol_h
$$
using the first relation in Lemma \ref{J1-int}. But Hopf's formula for $\chi(M)$ \cite[Theorem 5.7]{Gray} 
states that
$$
    24 \pi^2 \chi(M) = 18 \int_M \sigma_4 (L) dvol_h.
$$
Thus, Newton's identity for $\sigma_4(L)$ gives
$$
    24 \pi^2 \chi(M) 
   = \int_M \left(\frac{9}{4} \Iv_1 - \frac{9}{2} \Iv_2 + 6H \tr(\lo^3) - 9 H^2 |\lo|^2 + 18 H^4 \right) dvol_h.
$$
Therefore, we get
$$
    \int_M \left({\text Wm} + \frac{11}{6} \Iv_1 - \frac{5}{2} \Iv_2\right) dvol_h 
    - 24 \pi^2 \chi(M) 
   = \int_M \left(- 6 |dH|^2 + 6 H^2 |\lo|^2 - 18 H^4\right) dvol_h.
$$
The proof is complete.
\end{proof}

\subsection{Proof of the relation \eqref{JJD}}\label{diff-rel}

In the present section, we derive the relation \eqref{JJD}. We first prove the following result.

\begin{prop}\label{B-JJD} In the critical dimension $n=4$, it holds
\begin{align}\label{JJ-deco}
    \Jv_1 - 2 \Jv_2 = (D(\lo),\lo) + 3 (\lo^2,\W) + 2 \lo^{kl} \lo^{ij} W_{kijl} - \frac{4}{3} |\lo|^4
    + 3  \tr(\lo^4) + \delta (\lo,\overline{W}_0).
\end{align}
\end{prop}

Since all other terms on the right-hand side of \eqref{JJ-deco} are conformally invariant, this identity confirms
the conformal invariance of $\Jv_1-2\Jv_2$.

\begin{proof} In general dimensions, it holds
$$
   \bar{\nabla}^k(\overline{W})_{kij0} = (n-2) \bar{C}_{ij0} 
    = (n-2) \bar{\nabla}^0(\bar{\Rho})_{ij} - (n-2) \bar{\nabla}_j(\bar{\Rho})_{i0}.
$$
On the other hand, for tangential $\partial_k$, we find 
\begin{align}\label{n-reduce}
     \bar{\nabla}^k (\overline{W})_{kij0} & = \nabla^k \overline{W}_{kij0} - L^{kl} \overline{W}_{kijl}
     + L^k_i \overline{W}_{k0j0} + n H \overline{W}_{0ij0} \notag \\
     & = \nabla^k \overline{W}_{kij0} - L^{kl} \overline{W}_{kijl} - L^k_i \overline{W}_{0kj0} 
     + n H \overline{W}_{0ij0} \notag \\
     & = \nabla^k \overline{W}_{kij0} - \lo^{kl} \overline{W}_{kijl} + H \overline{W}_{0ij0} - \lo_i^k \overline{W}_{0kj0}
     - H \overline{W}_{0ij0} + n H \overline{W}_{0ij0} \notag \\
     & = \nabla^k \overline{W}_{kij0} - \lo^{kl} \overline{W}_{kijl} - \lo_i^k \overline{W}_{0kj0} 
    + n H \overline{W}_{0ij0}
\end{align}
using $\bar{\nabla}_i (\partial_j) = \nabla_i(\partial_j) - L_{ij} \partial_0$ and $\bar{\nabla}_k(\partial_0) 
= L_k^m \partial_m$. Hence
\begin{equation}\label{n-reduce-2}
    \bar{\nabla}^k(\overline{W})_{kij0} = \bar{\nabla}^0(\overline{W})_{0ij0} + 
    \nabla^k \overline{W}_{kij0} - \lo^{kl} \overline{W}_{kijl} - \lo_i^k \overline{W}_{0kj0} + n H \W_{ij}.
\end{equation}
Combining both results gives the relation
$$
    (n-2) \bar{\nabla}^0(\bar{\Rho})_{ij} - \bar{\nabla}^0(\overline{W})_{0ij0} 
   = \nabla^k \overline{W}_{kij0} + (n-2) \bar{\nabla}_j(\bar{\Rho})_{i0} 
    - \lo^{kl} \overline{W}_{kijl} - \lo_i^k \overline{W}_{0kj0} + n H \W_{ij}.
$$
Together with the trace-free Codazzi Mainardi equation \eqref{CM-TF-3}, we obtain
\begin{align*}
   &  (n-2) \bar{\nabla}^0(\bar{\Rho})_{ij} - \bar{\nabla}^0(\overline{W})_{0ij0} \\
   & = \nabla^k \nabla_i (\lo)_{kj} - \nabla^k \nabla_k (\lo)_{ij} + \frac{1}{n-1} \nabla^k \delta(\lo)_i h_{kj} 
   - \frac{1}{n-1} \nabla^k \delta(\lo)_k h_{ij} \\
    & + (n-2) \bar{\nabla}_j(\bar{\Rho})_{i0}  - \lo^{kl} \overline{W}_{kijl} - \lo_i^k \overline{W}_{0kj0} + n H \W_{ij}.
\end{align*}
Now we commute the covariant derivatives in the first term. Then
\begin{align*}
   &  (n-2) \bar{\nabla}^0(\bar{\Rho})_{ij} - \bar{\nabla}^0(\overline{W})_{0ij0} \\
   & = \nabla_i \delta(\lo)_j - R^k{}_{ijl} L^l_k - R^k{}_{ikl} L^l_j  - \Delta(\lo)_{ij} 
   + \frac{1}{n-1} \nabla_j \delta(\lo)_i - \frac{1}{n-1} \delta \delta (\lo) h_{ij} \\
   & + (n-2) \bar{\nabla}_j(\bar{\Rho})_{i0} - \lo^{kl} \overline{W}_{kijl} - \lo_i^k \overline{W}_{0kj0} + n H \W_{ij}.
\end{align*}
Next, we note that
\begin{align*}
   \bar{\nabla}_j(\bar{\Rho})_{i0} & = \partial_j (\bar{\Rho}_{0i}) - \bar{\Rho} (\bar{\nabla}_j(\partial_i),\partial_0) 
   - \bar{\Rho}(\partial_i, \bar{\nabla}_j(\partial_0)) \\
   & = \partial_j (\bar{\Rho}_{0i}) - \bar{\Rho}(\nabla_j(\partial_i),\partial_0) + L_{ij} \bar{\Rho}_{00}  
   -  L_j^l \bar{\Rho}(\partial_i,\partial_l) \\
   &=  \nabla_j (\bar{\Rho}_0)_{i} + L_{ij} \bar{\Rho}_{00} - L_j^l \bar{\Rho}_{il},
\end{align*}
and that the Codazzi-Mainardi equation implies
$$
    (n-1) \nabla_j(\bar{\Rho}_0)_i = \nabla_j \delta(\lo)_i - (n-1) \Hess_{ij}(H).
$$
These identities yield
\begin{align}\label{diff-main}
   &  (n-2) \bar{\nabla}^0(\bar{\Rho})_{ij} - \bar{\nabla}^0(\overline{W})_{0ij0} \notag \\
   & = - \Delta(\lo)_{ij} + \nabla_i \delta(\lo)_j + \nabla_j \delta(\lo)_i  - (n-2) \Hess_{ij}(H) 
   - \frac{1}{n-1} \delta \delta (\lo) h_{ij} \\
   & - (n-2) L_j^l \bar{\Rho}_{il} + (n-2) L_{ij} \bar{\Rho}_{00} 
   - R^k{}_{ijl} L^l_k - R^k{}_{ikl} L^l_j - \lo^{kl} \overline{W}_{kijl} - \lo_i^k \overline{W}_{0kj0} 
   + n H \W_{ij}. \notag 
\end{align}
The identity \eqref{diff-main} will also be important in Section \ref{BGW-1}. In the critical dimension $n=4$, we 
contract this identity with $\lo$ and obtain
\begin{align}\label{N-diff}
   & 2 \lo^{ij} \bar{\nabla}^0(\bar{\Rho})_{ij} - \lo^{ij} \bar{\nabla}^0(\overline{W})_{0ij0} \notag \\
   & = - (\lo,\Delta(\lo)) + 2 \lo^{ij} \nabla_i \delta(\lo)_j - 2 (\lo,\Hess(H)) - 2 (\lo^2,\bar{\Rho}) 
   - 2 H (\lo, \bar{\Rho}) + 2 \bar{\Rho}_{00} |\lo|^2 \notag \\
   &  -\lo^{ij} \lo^{kl} R_{kijl} + 2 (\lo^2,\Rho) + \J |\lo|^2 
   - \lo^{ij} \lo^{kl} \overline{W}_{kijl} - (\lo^2,\W) + 4 H (\lo,\W).
\end{align}
Now the definitions of $\Jv_1$ and $\Jv_2$ give
\begin{align*}
    & 2 \Jv_2 - \Jv_1 \notag \\
    & = 2 \lo^{ij} \bar{\nabla}_0 (\bar{\Rho})_{ij} - \lo^{ij} \bar{\nabla}_0(\overline{W})_{0ij0} \notag \\
    & + 2  (\lo^2,\Rho) + 2 H (\lo,\Rho) - 3  H (\lo,\W) + 2 (\lo,\Hess(H)) 
   - 2 \bar{\Rho}_{00} |\lo|^2 - 3 H^2 |\lo|^2 - H \tr(\lo^3) \notag \\
    & + \frac{2}{3} |\delta(\lo)|^2 - \delta \delta (\lo^2) + \frac{1}{2} \Delta (|\lo|^2).
\end{align*}
In this formula, we substitute the first line on the right-hand side by the sum displayed in \eqref{N-diff}. This yields
\begin{align*}
    & 2 \Jv_2 - \Jv_1 \notag \\
    & = - (\lo,\Delta(\lo)) + 2 \lo^{ij} \nabla_i \delta(\lo)_j  - 2 (\lo^2,\bar{\Rho}-\Rho) 
   - 2 H (\lo, \bar{\Rho} - \Rho) + 2 (\lo^2,\Rho) + \J |\lo|^2 \notag \\
    &  -\lo^{ij} \lo^{kl} R_{kijl} - \lo^{ij} \lo^{kl} \overline{W}_{kijl} - (\lo^2,\W) + H (\lo,\W) 
    - 3 H^2 |\lo|^2 - H \tr(\lo^3) \notag \\
    & + \frac{2}{3} |\delta(\lo)|^2 - \delta \delta (\lo^2) + \frac{1}{2} \Delta (|\lo|^2).
\end{align*}
In order to simplify that result, we apply the identities
\begin{align*}
    \delta (\lo \delta(\lo)) & = |\delta(\lo)|^2 + \lo^{ij} \nabla_i \delta(\lo)_j, \\
    \lo^{ij} \lo^{ks} R_{kijs} & = - 2 (\lo^2,\Rho) +  \lo^{ij} \lo^{kl} W_{kijl},
\end{align*}
together with the consequences 
\begin{align*}
    2(\lo,\bar{\Rho}-\Rho) & = - 2 H |\lo|^2 + \tr(\lo^3) + (\lo,\W), \\
    2(\lo^2,\bar{\Rho}-\Rho) & = - 2 H \tr(\lo^3) - H^2 |\lo|^2 + \tr(\lo^4) - \frac{1}{6} |\lo|^4 + (\lo^2,\W)
\end{align*}
of the Fialkow equation. Then we obtain
\begin{align*}
     2 \Jv_2 - \Jv_1 & = - (\lo,\Delta(\lo)) + 4 (\lo^2,\Rho) + \J |\lo|^2 \\
     & - \lo^{ij} \lo^{kl} W_{kijl} - \lo^{ij} \lo^{kl} \overline{W}_{kijl} - 2 (\lo^2,\W) - \tr(\lo^4) + \frac{1}{6} |\lo|^4 \\
     & - \frac{4}{3} |\delta(\lo)|^2 + 2  \delta(\lo \delta(\lo)) - \delta \delta (\lo^2) + \frac{1}{2} \Delta (|\lo|^2).
\end{align*}
Finally, \eqref{II5} shows that
$$
     \lo^{ij} \lo^{kl} \overline{W}_{kijl} =  \lo^{ij} \lo^{kl} W_{kijl} - \frac{7}{6} |\lo|^4 + 2 \tr(\lo^4) + (\lo^2,\W).
$$ 
Therefore, we get the final result
\begin{align*}
    2 \Jv_2 - \Jv_1 & = - (\lo,\Delta(\lo)) + 4 (\lo^2,\Rho) + \J |\lo|^2 \\
    & - 2 \lo^{ij} \lo^{kl} W_{kijl} - 3 (\lo^2,\W) + \frac{4}{3} |\lo|^4 - 3 \tr(\lo^4) \\
    & - \frac{4}{3} |\delta(\lo)|^2 + 2  \delta(\lo \delta(\lo))  - \delta \delta (\lo^2) + \frac{1}{2} \Delta (|\lo|^2).
\end{align*}
Combining this result with
$$
    (D(\lo),\lo) 
   = (\Delta(\lo),\lo) - 4 (\lo^2,\Rho) - \J |\lo|^2 + \frac{4}{3} |\delta (\lo)|^2 - \frac{4}{3} \delta (\lo \delta(\lo))
$$
(by \eqref{D-action}) finally yields
\begin{align*}
    \Jv_1 - 2 \Jv_2 & = (D(\lo),\lo) + 3 (\lo^2,\W) + 2 \lo^{kl} \lo^{ij} W_{kijl} - \frac{4}{3} |\lo|^4  + 3  \tr(\lo^4) \\
    & + \delta \delta (\lo^2) - \frac{1}{2} \Delta (|\lo|^2) - \frac{2}{3} \delta (\lo \delta(\lo)).
\end{align*}
Now Lemma \ref{van} completes the proof.
\end{proof}

\begin{cor}\label{J12-LC} In the critical dimension $n=4$, it holds
$$
   \Jv_1 - 2 \Jv_2 = -\frac{4}{3} \Iv_1 + 3 \Iv_2 + \Iv_5 + 3 \Iv_6 - \frac{1}{2} \Iv_7.
$$
\end{cor}

\begin{proof} Combine Proposition \ref{B-JJD} with Remark \ref{action-local}.
\end{proof}

This proves the relation \eqref{JJD}.

\subsection{Some comments on \cite{BGW-1}}\label{BGW-1}

Remark \ref{BGW-1-compare} shows that the formula for the second-order part of ${\bf P}_4$ 
given in \cite{BGW-1} is equivalent to our formula for this part. 

In the present section, we prove that the formula for ${\bf Q}_4$ in general dimensions displayed in 
Theorem \ref{main} is equivalent to the formula
\begin{align}
    {\bf Q}_4 & = Q_4 + \tfrac{2}{n} (\lo,\Delta(\lo)) \notag \\ 
    & + \tfrac{1}{n-3} \Big[ 2 (n-1) \delta \delta (\mathring{\JF}) 
    + \tfrac{3n^2-3n-2}{2n(n-1)} \Delta (|\lo|^2) + 4 \delta (\lo \delta(\lo)) \label{div-terms-2} \\
    & - 2(n-1) \lo^{ij} \overline{C}_{ij0} + \tfrac{6(n-1)}{n} \lo^{ij} \nabla^k \overline{W}_{kij0} \label{C-W-terms} \\
    & - 2 (n-1)(n-4) (\mathring{\JF},\Rho) - 4(n-5) (\lo^2,\Rho) \label{F-P-terms} \\
    & + 2 (n-1)(n-2) H (\lo,\mathring{\JF}) + 2(n-1)^2 (\mathring{\JF},\mathring{\JF}) 
   + 2(n-1) (\lo^2,\mathring{\JF}) \label{F-terms} \\
   & - \tfrac{n^3+5n^2-20n+20}{2n(n-1)} \J |\lo|^2 - \tfrac{2(n+3)}{n} \lo^{ij} \lo^{kl} W_{kijl} \Big]
   \label{J-W-terms}
\end{align}
displayed in \cite[Corollary 1.1]{BGW-1}. Here we omit the terms which are quartic in $L$. Accordingly, we shall 
omit the verification of the coincidence of the respective terms which are quartic in $L$.

In order to verify the equivalance of both formulas, we calculate the difference of the sum
$$
   \mbox{\eqref{div-terms} + \eqref{d-terms} + \eqref{W-terms} + \eqref{P-terms}}
$$
and the sum 
$$
   \mbox{\eqref{div-terms-2}} + \tfrac{2}{n} (\lo,\Delta(\lo)) + \mbox{\eqref{C-W-terms} + \eqref{F-P-terms} 
   + \eqref{F-terms} + \eqref{J-W-terms}},
$$
up to terms which are quartic in $L$. 

We first verify that the divergence terms in \eqref{div-terms} coincide with the divergence terms in 
\eqref{div-terms-2}. Indeed, \eqref{TF-F} implies
$$
   \delta \delta (\mathring{\JF}) = \tfrac{1}{n-2} \delta \delta (\lo^2) - \tfrac{1}{n (n-2)} \Delta (|\lo|^2 
   + \tfrac{1}{n-2} \delta \delta (\W).
$$ 
Hence \eqref{div-terms-2} reads
$$
     \tfrac{2(n-1)}{(n-3)(n-2)} \delta \delta (\W) +  \tfrac{2(n-1)}{(n-3)(n-2)} \delta \delta (\lo^2) 
     + \tfrac{3n-4}{2(n-1)(n-2)} \Delta (|\lo|^2) + \tfrac{4}{n-3} \delta (\lo \delta (\lo)).
$$
But this sum coincides with \eqref{div-terms}. 

In order to proceed, we use the identity \eqref{diff-main} to replace the normal derivative term 
$(\lo,\bar{\nabla}_0(\bar{\Rho}))$ in \eqref{d-terms} by $\lo^{ij} \bar{\nabla}_0(\overline{W})_{0ij0}$. 
We find that the sum \eqref{d-terms} equals
\begin{align*}
     & 2 (\lo,\Hess(H)) + 2 (\lo,\bar{\nabla}_0 (\bar{\Rho})) 
     - \tfrac{2}{n-2} \lo^{ij} \bar{\nabla}^0 (\overline{W})_{0ij0} 
     - \tfrac{2(n-1)}{(n-2)(n-3)} \lo^{ij} \bar{\nabla}^0(\overline{W})_{0ij0} \notag \\
     & = - \tfrac{2}{n-2} (\lo,\Delta(\lo)) + \tfrac{4}{n-2} \lo^{ij} \nabla_i \delta(\lo)_j 
     - \tfrac{2(n-1)}{(n-2)(n-3)} \lo^{ij} \bar{\nabla}^0(\overline{W})_{0ij0} \notag \\
     & - 2 \lo^{ij} L_j^l \bar{\Rho}_{il} + 2 |\lo|^2 \bar{\Rho}_{00} 
     - \tfrac{2}{n-2} \lo^{ij} L^l_k R^k {}_{ijl} - \tfrac{2}{n-2} \lo^{ij} L^l_j R^k {}_{ikl} \notag \\
     & - \tfrac{2}{n-2} \lo^{ij} \lo^{kl} \overline{W}_{kijl} 
     - \tfrac{2}{n-2} (\lo^2,\W) + \tfrac{2n}{n-2}  H (\lo,\W).
\end{align*}
The difference of this sum and the sum 
\begin{align*}
    & \tfrac{2}{n} (\lo,\Delta(\lo)) - \tfrac{2(n-1)}{n-3} \lo^{ij} \bar{C}_{ij0} 
    + \tfrac{6(n-1)}{n(n-3)} \lo^{ij} \nabla^k \overline{W}_{kij0} \notag \\
    & = \tfrac{2}{n} (\lo,\Delta(\lo)) - \tfrac{2(n-1)}{(n-2)(n-3)} \lo^{ij} \bar{\nabla}^k (\overline{W})_{kij0} 
   + \tfrac{6(n-1)}{n(n-3)} \lo^{ij} \nabla^k \overline{W}_{kij0}
\end{align*}
(see \eqref{C-W-terms}) equals
\begin{align*}
    & - \tfrac{4 (n-1)}{n(n-2)} (\lo,\Delta(\lo)) 
    + \tfrac{2(n-1)}{(n-2)(n-3)} (\lo^{ij} \bar{\nabla}^k (\overline{W})_{kij0} 
    -  \lo^{ij} \bar{\nabla}^0(\overline{W})_{0ij0}) \notag \\
    & + \tfrac{4}{n-2} \lo^{ij} \nabla_i \delta(\lo)_j 
    - \tfrac{6(n-1)}{n(n-3)} \lo^{ij} \nabla^k \overline{W}_{kij0} \notag \\
    & - 2 \lo^{ij} L_j^l \bar{\Rho}_{il} + 2 |\lo|^2 \bar{\Rho}_{00} 
     - \tfrac{2}{n-2} \lo^{ij} L^l_k R^k {}_{ijl} - \tfrac{2}{n-2} \lo^{ij} L^l_j R^k {}_{ikl} \notag \\
     & - \tfrac{2}{n-2} \lo^{ij} \lo^{kl} \overline{W}_{kijl} 
     - \tfrac{2}{n-2} (\lo^2,\W) + \tfrac{2n}{n-2}  H (\lo,\W).
\end{align*}
Now the identities \eqref{Simons-1} and \eqref{n-reduce-2} imply that this sum equals 
\begin{align}\label{T4}
     & - \tfrac{4 (n-1)}{n(n-2)} \Big( n (\lo,\Hess(H)) + L^{ij} \nabla_i (\overline{\Ric}_0)_j 
     + L^{ij} \nabla^k \bar{R}_{ikj0} 
     + (L^2)^{ij} \bar{R}_{ik}{}^k{}_j - L^{ij} L^{rs} \bar{R}_{rijs} \Big) \notag \\
     & + \tfrac{2(n-1)}{(n-2)(n-3)} \left( \lo^{ij} \nabla^k \overline{W}_{kij0} - \lo^{ij} \lo^{kl} \overline{W}_{kijl} 
      - (\lo^2,\W) + n H (\lo,\W) \right)  \notag \\
     & + \tfrac{4}{n-2} \lo^{ij} \nabla_i \delta(\lo)_j - \tfrac{6(n-1)}{n(n-3)} \lo^{ij} \nabla^k \overline{W}_{kij0} \notag \\
      & - 2 \lo^{ij} L_j^l \bar{\Rho}_{il} + 2 |\lo|^2 \bar{\Rho}_{00} 
     - \tfrac{2}{n-2} \lo^{ij} L^l_k R^k {}_{ijl} - \tfrac{2}{n-2} \lo^{ij} L^l_j R^k {}_{ikl} \notag \\
     & - \tfrac{2}{n-2} \lo^{ij} \lo^{kl} \overline{W}_{kijl} 
     - \tfrac{2}{n-2} (\lo^2,\W) + \tfrac{2n}{n-2}  H (\lo,\W),
\end{align}
up to terms which are quartic in $L$. The decomposition $\bar{R} = \overline{W} - \bar{\Rho} \owedge g$ yields
$$
    L^{ij} \nabla^k \bar{R}_{ikj0} = - L^{ij} \nabla^k \overline{W}_{kij0} + L^{ij} \nabla_j (\bar{\Rho}_0)_i 
    - n H \delta(\bar{\Rho}_0).
$$
Therefore, the sum \eqref{T4} further reduces to the sum of
\begin{align*}
    & \left(\tfrac{4(n-1)}{n(n-2)} + \tfrac{2(n-1)}{(n-2)(n-3)} - \tfrac{6(n-1)}{n(n-3)}\right)
    L^{ij} \nabla^k \overline{W}_{kij0}  \\
    & - \tfrac{4(n-1)}{n(n-2)} \left[ n (\lo,\Hess(H)) + \lo^{ij} \nabla_i (\overline{\Ric}_0)_j + H \delta (\overline{\Ric}_0) 
    + \lo^{ij} \nabla_j (\bar{\Rho}_0)_i + H \delta (\bar{\Rho}_0)  - n H \delta(\bar{\Rho}_0) \right] \\
    & + \tfrac{4}{n-2} \lo^{ij} \nabla_i \delta(\lo)_j,
\end{align*}
and
\begin{align*}
    & -\tfrac{4(n-1)}{n(n-2)} ( (L^2)^{ij} \bar{R}_{ik}{}^k{}_j - L^{ij} L^{rs} \bar{R}_{rijs}) \\
    & - \tfrac{2(n-1)}{(n-2)(n-3)} (\lo^{ij} \lo^{rs} \overline{W}_{rijs} + (\lo^2,\W) - n H (\lo,\W)) \\
    & - 2 \lo^{ij} L^l_j \bar{\Rho}_{il} + 2 |\lo|^2 \bar{\Rho}_{00} 
    - \tfrac{2}{n-2} \lo^{ij} L^{rs} R_{rijs} + \tfrac{2}{n-2} \lo^{ij} L^l_j \Ric_{il} \\
    & - \tfrac{2}{n-2} \lo^{ij} \lo^{rs} \overline{W}_{rijs} - \tfrac{2}{n-2} (\lo^2,\W) + \tfrac{2n}{n-2} H (\lo,\W),
\end{align*}
up to terms which are quartic in $L$. Now $\delta (\lo) = (n-1) dH + (n-1) \bar{\Rho}_0$ (Codazzi-Mainardi) implies 
that the first sum vanishes. In the second sum, we use the Gauss equations for $\bar{R}$, \eqref{I5-bar} and 
the Fialkow equation to replace curvature contributions of the background metric $g$ by curvature contributions of the 
induced metric $h$. Simplification gives
\begin{align}\label{T6}
    & -\tfrac{4(n-1)}{n(n-2)} ( (\lo^2,\Ric) - \lo^{ij} \lo^{rs} R_{rijs}) \notag \\ 
    & - (\tfrac{2(n-1)}{(n-2)(n-3)} + \tfrac{2}{n-2}) ( \lo^{ij} \lo^{rs} W_{rijs} + \tfrac{2}{n-2} (\lo^2,\W)) \notag \\
    & - (\tfrac{2(n-1)}{(n-2)(n-3)} + \tfrac{2}{n-2}) (\lo^2,\W) + (\tfrac{2n(n-1)}{(n-2)(n-3)} 
    + \tfrac{2n}{n-2}) H (\lo,\W) \notag \\
    & - 2 \lo^{ij} L^l_j (\tfrac{1}{n-2} \W_{il} + \Rho_{il}) + 2 |\lo|^2 \bar{\Rho}_{00} 
    - \tfrac{2}{n-2} \lo^{ij} \lo^{rs} R_{rijs} + \tfrac{2}{n-2} (\lo^2, \Ric) \notag \\
    & = -\tfrac{4(n-1)}{n(n-2)} ( (\lo^2,\Ric) - \lo^{ij} \lo^{rs} R_{rijs}) \notag \\
    & -\tfrac{4}{n-3} (\lo^{ij} \lo^{rs} W_{rijs} + \tfrac{2}{n-2} (\lo^2,\W)) 
    - \tfrac{4}{n-3} (\lo^2,\W) + \tfrac{4n}{n-3}  H (\lo,\W) \notag \\
    & - 2 \lo^{ij} L^l_j (\tfrac{1}{n-2} \W_{il} + \Rho_{il}) + 2 |\lo|^2 \bar{\Rho}_{00} 
    - \tfrac{2}{n-2} \lo^{ij} \lo^{rs} R_{rijs} + \tfrac{2}{n-2} (\lo^2, \Ric),
\end{align}
up to terms which are quartic in $L$. In the latter formula, we decompose $R = W - \Rho \owedge h$ and 
substitute $\Ric = (n-2) \Rho + \J h$. 

Finally, the terms in \eqref{F-P-terms}-\eqref{J-W-terms} yield
\begin{align}\label{T5}
   & - \tfrac{2(n-1)(n-4)}{(n-2)(n-3)} (\Rho,\W) + \tfrac{2(n-1)^2}{(n-2)^2(n-3)} (\W,\W) 
   + \big(\tfrac{4(n-1)^2}{(n-2)^2(n-3)} + \tfrac{2(n-1)}{(n-2)(n-3)}\big) (\lo^2,\W) \notag \\
   & - \tfrac{6n^2-38n+48}{(n-2)(n-3)} (\lo^2,\Rho) + \tfrac{2(n-1)(n-4)}{n(n-2)(n-3)} \J |\lo|^2
   + \tfrac{2(n-1)}{n-3} H (\lo,\W).
\end{align}

Now, we use these results to determine the remaining terms in the difference of both formulas for ${\bf Q}_4$.

\begin{itemize}
\item The term $(\Rho,\W)$. Its coefficient 
$$
    - \tfrac{2(n-1)(n-4)}{(n-2)(n-3)} 
$$
in \eqref{T5} coincides with the coefficient in \eqref{W-terms}. 

\item The term $(\W,\W)$. Its coefficient
$$
   \tfrac{2(n-1)^2}{(n-2)^2(n-3)}
$$
in  \eqref{T5} coincides  with its coefficient in \eqref{W-terms}. 

\item The term $H (\lo,\W)$. It contributes to \eqref{W-terms} and to \eqref{T6} with the respective coefficients
$$
   - \tfrac{2(n-1)^2}{(n-2)(n-3)} \quad \mbox{and} \quad \tfrac{4n}{n-3} - \tfrac{2}{n-2}. 
$$
The sum of these coefficients equals 
$$
    \tfrac{2(n-1)}{n-3}.
$$
On the other hand, it contributes to \eqref{T5} with the coefficient
$$
   \tfrac{2(n-1)}{n-3}.
$$
Thus, the term $H (\lo,\W)$ has the same coefficient in both formulas for ${\bf Q}_4$. The fact that, for $n=4$, this 
coefficient equals $-6$ is reflected by the contribution $-3\Jv_1$ in the decomposition of ${\bf Q}_4$ 
in Corollary \ref{Q4-crit}.

\item The term $(\lo^2,\W)$. It contributes to \eqref{W-terms} and \eqref{T6} by 
$$
   \tfrac{4(3n-5)(n-1)}{(n-2)^2(n-3)} \quad \mbox{and} 
   \quad - \tfrac{4}{n-3} \tfrac{2}{n-2} - \tfrac{4}{n-3} - \tfrac{2}{n-2} = - \tfrac{6(n-1)}{(n-2)(n-3)}.
$$
The sum of these coefficients equals
$$
   \tfrac{2(3n-4)(n-1)}{(n-2)^2(n-3)}.
$$
On the other hand, it contributes to \eqref{T5} by   
$$
    \tfrac{4(n-1)^2}{(n-2)^2(n-3)} + \tfrac{2(n-1)}{(n-2)(n-3)} =  \tfrac{2(3n-4)(n-1)}{(n-2)^2(n-3)}.
$$
Thus, the contributions of $(\lo^2,\W)$ to both formulas for ${\bf Q}_4$ coincide.
\item The term $(\lo^2,\Rho)$. It contributes to \eqref{P-terms} and \eqref{T6} with the respective coefficients
$$
    -\tfrac{2(n^2-9n+12)}{(n-2)(n-3)} \quad \mbox{and} \quad -4.                       
$$
The sum of these coefficients equals
$$
    -\tfrac{6n^2-38n+48}{(n-2)(n-3)}.
$$
This coefficient coincides with the coefficient of $(\lo^2,\Rho)$ in \eqref{T5}. 
\item The term $\J |\lo|^2$. It contributes to \eqref{T6} and \eqref{P-terms} with the respective coefficients
$$
      -\tfrac{4(n-1)}{n(n-2)} + \tfrac{2}{n-2} 
     = - \tfrac{2}{n} \quad \mbox{and} \quad - \tfrac{n^3-5n^2+18n-20}{2(n-3)(n-2)(n-1)}.
$$
The sum of these coefficients coincides with the sum 
$$
    \tfrac{2(n-1)(n-4)}{n(n-2)(n-3)} - \tfrac{n^3+5n^2-20n+20}{2n(n-1)} = -\tfrac{n^4 -n^3 - 6n^2 + 24n -24}{2(n(n-1)(n-2)(n-3)}
$$
of its contributions to \eqref{T5} and \eqref{J-W-terms}.

\item The term $\lo^{ij} \lo^{rs} W^{rijs}$. It contributes to \eqref{T6} with the coefficient 
$$
    \tfrac{4(n-1)}{n(n-2)} - \tfrac{4}{n-3} - \tfrac{2}{n-2} = - \tfrac{2(n+3)}{n(n-3)}.
$$
This coefficient coincides with its coefficient in \eqref{J-W-terms}.
\item The terms $|\lo|^2 \bar{\Rho}_{00}$ and $H (\lo,\Rho)$. Their contributions in \eqref{P-terms} cancel against
their respective contributions in \eqref{T6}.
\end{itemize}

In the critical dimension $n=4$, we may proceed more directly to compare the formulas for 
${\bf Q}_4$ displayed in Corollary \ref{Q4-crit} and \cite[Theorem 1.2]{BGW-1}. 
 
\cite[Theorem 1.2]{BGW-1} is the special case of \cite[Corollary 1.1]{BGW-1} for $n=4$. It gives 
the decomposition
\begin{equation}\label{Q4-BGW-1-0}
   {\bf Q}_4 = Q_4 + {\text Wm} + U + \mbox{divergence term}
\end{equation}
with ${\text Wm}$ as in \eqref{Wm-def}, the local conformal invariant
$$
   U \st 18 (\mathring{\JF},\mathring{\JF}) + 6 (\lo^2,\mathring{\JF})
   + \frac{49}{24} |\lo|^4 + \frac{9}{2} \lo^{ij} \nabla^k \overline{W}_{kij0},
   - \frac{7}{2} \lo^{ij} \lo^{kl} W_{iklj}
$$
and the divergence terms
\begin{align}\label{TD}
    \frac{8}{3} \delta (\lo \delta(\lo)) + 6 \delta \delta(\mathring{\JF}) - \frac{1}{12} \Delta (|\lo|^2)
   = \frac{8}{3} \delta (\lo \delta(\lo)) + 3 \delta \delta (\lo^2) -  \frac{5}{6} \Delta (|\lo|^2)  + 3 \delta \delta (\W).
\end{align}

We rewrite the sum in \eqref{Q4-BGW-1-0} in terms of the invariants $\Iv_j$, $\Jv_1$, and a divergence term.
Combining Proposition \ref{WM-deco}, Remark \ref{action-local} and \eqref{II5} gives
\begin{align*}
   {\text Wm} & = \left( -\frac{1}{2} \Iv_5 - \frac{1}{4} \Iv_7  - \frac{1}{2} \delta (\lo, \overline{W}_0) \right)  
   + \left(- 3 \Jv_1 + 3 \bar{\Iv}_5 + 3 \Iv_6 - \frac{3}{2} \Iv_7 - 6 \delta (\lo, \overline{W}_0) \right)  \\
   & = -\frac{7}{2} \Iv_1 + 6 \Iv_2  + \frac{5}{2} \Iv_5 + 6 \Iv_6 - \frac{7}{4} \Iv_7
   - 3 \Jv_1 - \frac{13}{2} \delta (\lo, \overline{W}_0)
\end{align*}
Moreover, we easily calculate
\begin{align*}
   U & = \frac{15}{2} \tr(\lo^4) + \frac{1}{6} |\lo|^4 + \frac{9}{2} |\W|^2 + 12 (\lo^2,\W)
   + \frac{9}{2} \lo^{ij} \nabla^k \overline{W}_{kij0} - \frac{7}{2} \lo^{ij} \lo^{kl} W_{iklj} \\
   & = \frac{1}{6} \Iv_1 + \frac{15}{2} \Iv_2 + \frac{9}{2} \Iv_4 - \frac{7}{2} \Iv_5 +12 \Iv_6 + \frac{9}{4} \Iv_7
    + \frac{9}{2} \delta (\lo,\overline{W}_{0})
\end{align*}
using Lemma \ref{div-term}. Finally, the sum \eqref{TD} coincides with the sum 
$$
     4 \delta (\lo \delta(\lo)) + \delta \delta(\lo^2)
    + \frac{1}{6} \Delta (|\lo|^2) + 3 \delta \delta (\W) + 2 \delta (\lo,\overline{W}_0)
$$
using Lemma \ref{van}. Note that the latter sum coincides with the second line of \eqref{Q4-ex2}, up to 
the last term (being a local conformal invariant).  

Summarizing these result, we find that \eqref{Q4-BGW-1-0} reads
\begin{align*}
    {\bf Q}_4 & = Q_4 - \frac{10}{3} \Iv_1 + \frac{27}{2} \Iv_2 + \frac{9}{2} \Iv_4 - \Iv_5
    + 18 \Iv_6 + \frac{1}{2} \Iv_7 - 3 \Jv_1 \\
    & + 4 \delta (\lo \delta(\lo)) + \delta \delta(\lo^2) + \frac{1}{6} \Delta (|\lo|^2) + 3 \delta \delta (\W).
\end{align*}
This shows that the formula for ${\bf Q}_4$ in Corollary \ref{Q4-crit} (or equivalently in Theorem \ref{alex}) coincides 
with the formula in \cite[Theorem 1.2]{BGW-1}.



\begin{thebibliography}{CHBRS21}

\bibitem [AGV81]{AGV}
E.~Abbena, A.~Gray and L.~Vanhecke, Steiner's formula for the volume of a parallel hypersurface in a
Riemannian manifold, {\em Annali Sc. Norm. Sup. Pisa} {\bf 8}, (3), (1981), 473--493.

\bibitem [AHR]{AHR}
J.~B.~Achour, E.~Huguet and J.~Renaud, {\em Conformally invariant wave equation for a symmetric second
rank tensor (spin-2) in $d$-dimensional curved background},  Physical Review D {\bf 89}, 064041 (2014).
\url{arXiv:1311.43124v3}

\bibitem [A12]{alex}
S.~Alexakis, {\em The decomposition of global conformal invariants}, Annals of Mathematics Studies
{\bf 182}, Princeton University Press, Princeton, NJ, 2012.

\bibitem [AGW21]{AGW}
C.~Arias, R.~Gover and A.~Waldron, Conformal geometry of embedded manifolds with boundary from universal
holographic formulae, {\em Advances in Math.} {\bf 384} (2021) 107700. \url{arXiv:1906.01731}

\bibitem [AS21]{AS}
A.~Astaneh and S.~Solodukhin, Boundary conformal invariants and the conformal anomaly in five dimensions.
{\em Physics Letters B} {\bf 816}, (2021), 136282. \url{arXiv:2102.07661}

\bibitem [AS22]{AS2}
A.~Astaneh and S.~Solodukhin, private communication.

\bibitem [BJ10]{BJ}
H.~Baum and A.~Juhl, {\em Conformal Differential Geometry: $Q$-Curvature and Conformal Holonomy}.
Oberwolfach Seminars {\bf 40}, 2010.

\bibitem [B87]{Besse}
A.~Besse,  {\em Einstein manifolds}, Ergebnisse der Mathematik und ihrer Grenzgebiete,
{\bf 10}, Springer-Verlag, (1987).

\bibitem [BGW21a]{BGW-2}
S.~Blitz, R.~Gover and A.~Waldron, Conformal fundamental forms and the asymptotically Poincar\'e-Einstein condition.
\url{arXiv:2107.10381v1}

\bibitem [BGW21b]{BGW-1}
S.~Blitz, R.~Gover and A.~Waldron, Generalized Willmore energies, $Q$-curvatures, extrinsic Paneitz
operators, and extrinsic Laplacian powers. \url{arXiv:2111.00179v1}

\bibitem [B95]{sharp}
T.~Branson, Sharp inequalities, the functional determinant, and the complementary series,
{\em Trans. Amer. Math. Soc.} {\bf 347}, (1995), 3671--3742.

\bibitem [B96]{br-NL}
T.~Branson, Nonlinear phenomena in the spectral theory of geometric linear differential operators,
{\em Proc. Symp. Pure Math.} {\bf 59}, (1996), 27-65.

\bibitem [B05]{br-last}
T.~Branson, $Q$-curvature and spectral invariants. {\em Rend. Circ. Mat. Palermo (2) Suppl.} {\bf 75},
(2005), 11--55.

\bibitem [BG94]{BG-BVP}
T.~Branson, The functional determinant of a four-dimensional boundary value problem,
{\em Trans. Amer. Math. Soc} {\bf 344}, (2), (1994), 479--531.

\bibitem [BGKV97]{BGKV}
T.~Branson, P.~Gilkey, K.~Kirsten and D.~Vassilevich,
Heat kernel asymptotics with mixed boundary conditions, {\em Nucl. Phys. B} {\bf 563}, (3), (1999), 603--626.

\bibitem [CHBRS21]{CHBRS}
A.~Chalabi, C.~Herzog, A.~O'Bannon, R.~Robinson and J.~Sisti, Weyl anomalies of four-dimensional
boundaries and defects, {\em J. High Energ. Phys.} {\bf 166}, (2022). \url{arXiv:2111.14713v1}

\bibitem [CG19]{CG-1}
S.-Y.~Alice Chang and Y.~Ge, Compactness of conformally compact Einstein manifolds in dimension $4$. 
{\em Advances in Math.} {\bf 340}, (2018), 588--652.
\url{arXiv:1809.05593v2}

\bibitem [CQ97]{CQ}
S.-Y.~Alice Chang and J.~Qing, The zeta functional determinants on manifolds with boundary,
{\em J. Funct. Anal.} {\bf 147}, (1997), 327--362.

\bibitem [CY95]{CY}
S.-Y.~Alice Chang and P.~Yang, Extremal metrics of zeta function determinant on $4$-manifolds,
{\em Ann. of Math.} (2) {\bf 142}, (1), (1995), 171--212.

\bibitem [CGY02]{CGY}
S.-Y.~Alice Chang, M.~Gursky and P.~Yang, An equation of Monge-Ampere type in conformal geometry
and $4$-manifolds of positive Ricci curvature, {\em Ann. of Math.} (2) {\bf 155}, (3), (2002), 709--787.

\bibitem [C05]{CBull}
S.-Y.~Alice Chang, Conformal invariants and partial differential equations,
{\em Bull. Amer. Math. Soc. (N.S.)}, {\bf 42}, (3), (2005), 365--393.

\bibitem [C18]{IMC}
S.-Y.~Alice Chang,  Conformal Geometry on Four Manifolds, {\em Proc. Int. Cong. of Math.} (2018), 
{\bf 1}, 119--146. \url{arXiv:1809.06339}

\bibitem [CMY21]{CMY}
S.-Y.~Alice Chang, S.~McKeown and P.~Yang, Scattering on singular Yamabe spaces. \url{arXiv:2109.02014}.

\bibitem [CK04]{CK}
B.~Chow and D.~Knopf, {\em The Ricci Flow: An Introduction}. Mathematical Surveys and Monographs
{\bf 110}, AMS (2004).

\bibitem [DS93]{DS}
S.~Deser and A.~Schwimmer, Geometric classification of conformal anomalies in arbitrary dimensions,
{\em Physics Letters B}, {\bf 309}, 279 (1993).

\bibitem[DGH08]{DGH}
Z.~Djadli, C.~Guillarmou and M.~Herzlich, {\em Op\'erateurs g\'eom\'etriques,
invariants conformes et vari\'et\'es asymptotiquement hyperboliques}, Panoramas et
Synth\`eses {\bf 26}, Soci\'et\'e Math\'ematique de France, 2008.

\bibitem[ES85]{ES}
M.~Eastwood and M.~Singer, A conformally invariant Maxwell gauge. {\em Physics Letters} 
{\bf 107A}, 2, (1985), 73--74.

\bibitem [EO]{EO}
J.~Erdmenger and H.~Osborn, {\em Conformally covariant differential operators: symmetric tensor fields},
Class. Quantum Grav. {\bf 15}, (1998), 273-280.

\bibitem [FG12]{FG-final}
C.~Fefferman and C.~R.~Graham, {\em The Ambient Metric}. Annals of Math. Studies {\bf 178},
Princeton University Press, 2012. \url{arXiv:0710.0919}

\bibitem[FG12]{FG-J}
C.~Feffermann and C.~R.~Graham, Juhl's formulae for GJMS-operators and $Q$-curvatures.
{\em Journal Amer. Math. Soc.} {\bf 26}, 4, (2013), 1191--1207. \url{arXiv:1203.0360}

\bibitem [Fu15]{Fu}
D.~V.~Fursaev, Conformal anomalies of CFTs with boundary, {\em J. High Energy Phys.} {\bf 2015},
1--10 (2015). \url{arXiv:1510.01427v2}

\bibitem [FV11]{FV-book}
D.~V.~Fursaev and D.~Vassilevich, {\em Operators, Geometry and Quanta. Methods of Spectral Geometry in
Quantum Field Theory}, Springer (2011).

\bibitem [FT82]{FT}
E.~S.~Fradkin and A.~A.~Tseytlin, Asymptotic freedom in extended conformal supergravity,
{\em Physics Letters} {\bf 110B}, 2, (1982), 117--121.

\bibitem [G84]{Gilkey-book}
P.~Gilkey, {\em Invariance Theory, the Heat Equation, and the Atiyah-Singer Index Theorem}. Publish or 
Perish Inc. (1984).

\bibitem[GGHW19]{GGHW}
M.~Glaros, R.~Gover, M.~Halbasch and A.~Waldron, Variational calculus for hypersurface functionals:
singular Yamabe problem Willmore energies, {\em J. Geom. Phys.} {\bf 138}, (2019), 168--193.
\url{arXiv:1508.01838v1}

\bibitem [GW15]{GW-LNY}
R.~Gover, A.~Waldron, Conformal hypersurface geometry via a boundary Loewner-Nirenberg-Yamabe problem,
{\em Comm. in Analysis and Geometry} {\bf 29}, (4), 779-- 836. \url{arXiv:150602723v3}

\bibitem [GW17]{GW-RV}
R.~Gover, A.~Waldron,  Renormalized volume, {\em Comm. in Math. Physics} {\bf 354}, (3), (2017),
1205--1244. \url{arXiv:1603.07367}

\bibitem [GJMS92]{GJMS}
C.~R.~Graham, R.~Jenne, L.~J.~Mason and G.~A.~J.~Sparling, Conformally invariant powers of the Laplacian. {I}.
Existence, {\em J. London Math. Soc.} (2) {\bf 46}, (3), (1992), 557--565.

\bibitem [GZ03]{GZ}
C.~R.~Graham and M.~Zworski, Scattering matrix in conformal geometry, {\em Inventiones math.}
{\bf 152}, (1), (2003), 89--118. \url{arXiv:math/0109089}

\bibitem [G17]{G-SY}
C.~R.~Graham, Volume renormalization for singular Yamabe metrics, {\em Proc. Amer. Math. Soc.} {\bf 145}, 
(2017), 1781--1792. \url{arXiv:1606.00069}

\bibitem [GR20]{GR}
C.~R.~Graham and N. Reichert, Higher-dimensional Willmore energies via minimal submanifold asymptotics,
{\em Asian J. Math.} {\bf 24}, (4), 571--610 (2020). \url{arXiv:1704.03852}

\bibitem [GG19]{GG}
C.~R.~Graham and M.~Gursky, Chern-Gauss-Bonnet formula for singular Yamabe metrics in dimension four. 
\url{arXiv:1902.01562}

\bibitem [GW99]{GWi}
C.~R.~Graham and E.~Witten, Conformal anomaly of submanifold observables in AdS/CFT correspondence, 
{\em Nuclear Phys. B} {\bf 546}, (1999), no. 1-22, 52--64.

\bibitem [G04]{Gray}
A.~Gray, {\em Tubes},  Birkhäuser, Progress in Mathematics {\bf 221}, 2004.

\bibitem [GZ20]{GZ1}
M.~Gursky and S.~Zhang, Rigidity for Bach-flat metrics on manifolds with boundary and applications. 
\url{arXiv:2007.04355v2}

\bibitem [G05]{Guven}
J.~Guven, Conformally invariant bending energy for hypersurfaces, {\em J. Phys. A}, {\bf 38}, 7943--7955 (2005).

\bibitem [HH17]{HH}
C.~P.~Herzog and Kuo-Wei Huang, Boundary conformal field theory and a boundary central charge. \url{arXiv:1707.06224.v4}

\bibitem [HP99]{HP}
G.~Huisken and A.~Polden, {\em Geometric Evolution Equations for Hypersurfaces}, in S. Hildebrand, M. Struwe (eds)
{\em Calculus of Variations and Geometric Evolution Problems}, Lecture Notes in Mathematics {\bf 1713}, 1999, 45--84.

\bibitem [J88]{Jenne}
R.~W.~Jenne, {\em A construction of conformally invariant differential operators}. Thesis. University of Washington. 1988.

\bibitem [J09]{J1}
A.~Juhl, {\em Families of conformally covariant differential operators, $Q$-curvature and holography}, Birkhäuser,
Progress in Mathematics {\bf 275}, 2009.

\bibitem[J13]{J-ex}
A.~Juhl, Explicit formulas for GJMS-operators and $Q$-curvatures. {\em Geom. Funct. Anal.} {\bf 23}, (2013),
1278--1370. \url{arXiv:1108.0273}

\bibitem [J16]{J-heat}
A.~Juhl, Heat kernels, ambient metrics and conformal invariants, {\em Advances in Math.} {\bf 286}, (2016), 545--682.
\url{arXiv:1411.7851}

\bibitem [JO20]{JO0}
A.~Juhl and B.~Orsted, Shift operators, residue families and degenerate Laplacians, {\em Pacific J. of Math.}
{\bf 308}, (1), (2020), 103--160. \url{arXiv:1806.02556}

\bibitem [JO21]{JO1}
A.~Juhl and B.~Orsted, Residue families, singular Yamabe problems and extrinsic conformal Laplacians, 
{\em Advances in Math.}  {\bf 409} (2022) 108634, 1--158. \url{arXiv:2101.09027v3}

\bibitem [JO22]{JO2}
A.~Juhl and B.~{\O}rsted, On singular Yamabe obstructions. {\em J. Geom. Anal.} {\bf 32}, 146, (2022). 
\url{arXiv:2103.01552v2}

\bibitem [J21]{J-announce}
A.~Juhl, Extrinsic Paneitz operators and $Q$-curvatures for hypersurfaces. \url{arXiv:2110.04838v2}

\bibitem [KS15]{KS}
T.~Kobayashi and B.~Speh, {\em Symmetry breaking for representations of rank
one orthogonal groups}, Memoirs of AMS  {\bf 238}, (2015). \url{arXiv:1310.3213}

\bibitem [M13]{Matsumoto}
Y.~Matsumoto, {A GJMS construction for 2-tensors and the second variation of the total $Q$-curvature},
Pacific J. Math. {\bf 262}, (2), (2013), 437--455. \url{arXiv:1202.3227}

\bibitem [MN18]{mondino}
A.~Mondino and Huy The Nguyen, Global conformal invariants of submanifolds,
{\em Ann. Inst. Fourier}, {\bf 68}, (6), (2018), 2663--2695. \url{arXiv:1501.07527v2}

\bibitem [P08]{pan}
S.~Paneitz, A quartic conformally covariant differential operator for arbitrary
pseudo-{R}iemannian manifolds (summary), {\em SIGMA Symmetry Integrability  Geom.
Methods Appl.} {\bf 4} (2008), paper 036, 3p. \url{arXiv:0803.4331}

\bibitem [RT17]{RT}
M.~Rangamani and T.~Takayanagi, {\em Holographic Entanglement Entropy},  
Lecture Notes in Physics {\bf 931} (2017)

\bibitem[R84]{Rieg}
R.~Riegert, A nonlocal action for the trace anomaly. {\em Physics Letters B} {\bf 134}, no. 1--2, 56--60 (1984).

\bibitem [S08]{Sol-entropy}
S.~N.~Solodukhin, Entanglement entropy, conformal invariance and extrinsic geometry, {\em Physics Letters B},
{\bf 665}, 305--309 (2008). \url{arXiv:0802.3117}

\bibitem [S16]{Sol}
S.~N.~Solodukhin, Boundary terms of conformal anomaly, {\em Physics Letters B}, {\bf 752}, 131--134 (2016).
\url{arXiv:1510.04566v4}

\bibitem [V13]{V}
Y.~Vyatkin, {\em Manufacturing conformal invariants of hypersurfaces}, PhD thesis,  University of Auckland, 2013.

\bibitem [W93]{Will}
T.~J.~Willmore, {\em Riemannian Geometry}, Oxford Science Publications, 1993.

\bibitem [W86]{Wunsch}
V.~Wünsch, {\em On conformally invariant differential operators}, Math. Nachr. {\bf 129}, 269--281 (1986).

\bibitem [Z21]{Zhang}
Y.~Zhang, Graham-Witten's conformal invariant for closed four dimensional submanifolds.
{\em J. Math. Study} {\bf 54}, 200--226 (2021). \url{arXiv:1703.08611}

\end{thebibliography}
\end{document}